\documentclass[11pt,a4paper,twoside]{article}

\usepackage[T1]{fontenc}
\usepackage{lmodern}
\usepackage{authblk}
\usepackage[top=3cm,bottom=3cm,left=3cm,right=3cm,footskip=1cm]{geometry}
\usepackage{mathtools,amsthm,amssymb}
\usepackage{bbm}
\usepackage{xcolor}
\usepackage{tikz-cd}
\usepackage{capt-of}
\usepackage{enumitem}
\usepackage{needspace}
\usepackage{aliascnt}

\setlist[itemize]{noitemsep}
\setlist[enumerate]{itemsep=0mm}

\usepackage[
defernumbers=true,maxnames=6,style=numeric,
giveninits=true,sortcites=true,sorting=nty,
doi=true,isbn=true,url=false,safeinputenc
]{biblatex}
\definecolor{persianblue}{rgb}{0.11,0.22,0.73}
\definecolor{antiquefuchsia}{rgb}{0.57,0.36,0.51}
\usepackage[bookmarks=false,hidelinks]{hyperref}
\hypersetup{
	colorlinks=true,
	urlcolor=antiquefuchsia,
	linkcolor=persianblue,
	citecolor=persianblue
}

\usepackage[nameinlink,capitalize]{cleveref}

\newcommand{\newaliastheorem}[4][theorem]{%
	\newaliascnt{#2}{#1}%
	\newtheorem{#2}[#2]{#3}%
	\aliascntresetthe{#2}%
	\crefname{#2}{#3}{#4}%
}

\theoremstyle{plain}
\newtheorem{theorem}{Theorem}[section]
\crefname{theorem}{Theorem}{Theorems}
\newaliastheorem{lemma}{Lemma}{Lemmas}
\newaliastheorem{proposition}{Proposition}{Propositions}
\newaliastheorem{corollary}{Corollary}{Corollaries}
\newaliastheorem{claim}{Claim}{Claims}
\newaliastheorem{fact}{Fact}{Facts}
\newaliastheorem{notation}{Notation}{Notations}

\theoremstyle{definition}
\newaliastheorem{definition}{Definition}{Definitions}
\newaliastheorem{remark}{Remark}{Remarks}
\crefname{subsection}{Subsection}{Subsections}

\newtheorem{innernamedtheorem}{}
\newenvironment{namedtheorem}[1]
{\renewcommand{\theinnernamedtheorem}{#1}\innernamedtheorem}
{\endinnernamedtheorem}
\crefformat{innernamedtheorem}{#2#1#3}

\newcounter{proofassertionsequence}
\newcounter{proofassertion}
\newcommand{\resetproofassertions}{%
	\stepcounter{proofassertionsequence}%
	\setcounter{proofassertion}{0}%
}

\newlist{proofassertionitems}{enumerate}{1}
\newenvironment{proofassertiongroup}[1][]
{\refstepcounter{proofassertion}%
	\begin{proofassertionitems}[
		label=(\arabic{proofassertion}.\arabic*),
		ref=\arabic{proofassertion}.\arabic*,#1]}
	{\end{proofassertionitems}}

\crefname{forcingproofstep}{Step}{Steps}
\newcounter{uniformizationproofstep}
\crefname{uniformizationproofstep}{Step}{Steps}
\newcommand{\proofstepheading}[3]{%
	\par\medskip
	\refstepcounter{#1}\label{#2}%
	\noindent\textbf{Step~\arabic{#1}. #3}%
	\par\nobreak\smallskip\noindent\ignorespaces
}
\newcommand{\forcingstep}[2]{\proofstepheading{forcingproofstep}{#1}{#2}}
\newcommand{\uniformizationstep}[2]{\proofstepheading{uniformizationproofstep}{#1}{#2}}

\newcommand{\axiom}[1]{\mathrm{#1}}
\newcommand{\ad}{\axiom{AD}}
\newcommand{\dc}{\axiom{DC}}
\newcommand{\ac}{\axiom{AC}}
\newcommand{\zfc}{\axiom{ZFC}}
\newcommand{\zf}{\axiom{ZF}}
\newcommand{\ch}{\axiom{CH}}

\newcommand{\gch}{\axiom{GCH}}
\newcommand{\hpc}{\axiom{HPC}}

\newcommand{\ntbase}{\axiom{NT}_{\mathrm{base}}}
\newcommand{\hH}{\mathsf{H}}
\newcommand{\refl}{\mathsf{ref}}
\newcommand{\ord}{\mathrm{ord}}
\newcommand{\hwo}{\mathrm{HWO}}
\newcommand{\ro}{\mathsf{RO}}
\newcommand{\Los}{\L o\'s}

\DeclareMathOperator{\hod}{HOD}
\DeclareMathOperator{\crit}{crit}
\DeclareMathOperator{\domOp}{dom}
\DeclareMathOperator{\rangOp}{rang}
\DeclareMathOperator{\tcOp}{tc}
\DeclareMathOperator{\cofOp}{cof}
\DeclareMathOperator{\otpOp}{ot}
\DeclareMathOperator{\UltOp}{Ult}
\DeclareMathOperator{\ColOp}{Col}
\newcommand{\dom}[1]{\domOp(#1)}
\newcommand{\rang}[1]{\rangOp(#1)}
\newcommand{\tc}[1]{\tcOp(#1)}
\newcommand{\cof}[1]{\cofOp(#1)}
\newcommand{\otp}[1]{\otpOp(#1)}
\newcommand{\ult}[2]{\UltOp(#1,#2)}
\newcommand{\ultstar}[2]{\UltOp^*(#1,#2)}
\newcommand{\col}[2]{\ColOp(#1,#2)}
\newcommand{\dotcol}[2]{\mathop{\dot{\mathrm{C}}\mathrm{ol}}\nolimits(#1,#2)}

\newcommand{\card}[1]{{\vert #1 \vert} }
\DeclareRobustCommand{\ogonekaccent}[1]{\csname T1\string\k\endcsname{#1}}
\DeclareTextSymbol{\textlstroke}{T1}{170}
\renewcommand{\k}{\kappa}
\renewcommand{\l}{\lambda}
\newcommand{\treg}{\Theta_{\mathsf{reg}}}
\newcommand{\powerset}{\mathcal{P}}
\newcommand{\rest}{\upharpoonright}
\newcommand{\mcl}[1]{\mathcal{#1}}
\newcommand{\mbb}[1]{\mathbb{#1}}
\newcommand{\R}{\mathbb{R}}
\newcommand{\id}[1]{\mathrm{id}_{#1}}
\newcommand{\df}{=_{\mathrm{def}}}
\newcommand{\mdls}[1]{\models \text{``}#1\text{''}}
\newcommand{\forces}[2]{\Vdash_{#1} \text{``}#2\text{''}}
\newcommand{\1}[1]{\mathbbm{1}_{{#1}}}

\newcommand{\hq}{{\mathfrak{q}}}
\newcommand{\hr}{{\mathfrak{r}}}
\newcommand{\hs}{{\mathfrak{s}}}
\newcommand{\hw}{{\mathfrak{w}}}
\newcommand{\hm}{{\mathfrak{m}}}

\makeatletter
\let\pendingproofname\@empty

\newcommand{\TrackProof}[2]{%
	\AddToHook{env/#1/end}{\protected@xdef\pendingproofname{#2}}}

\AddToHook{env/proof/begin}{%
	\let\thisproofname\pendingproofname
	\global\let\pendingproofname\@empty
	\renewcommand{\qedsymbol}{%
		\ensuremath{\square}%
		\ifx\thisproofname\@empty\else~\textnormal{(\thisproofname)}\fi}}
\makeatother

\TrackProof{lemma}{Lemma~\thelemma}
\TrackProof{theorem}{Theorem~\thetheorem}
\TrackProof{proposition}{Proposition~\theproposition}
\TrackProof{corollary}{Corollary~\thecorollary}
\TrackProof{claim}{Claim~\theclaim}
\TrackProof{fact}{Fact~\thefact}
\begin{document}
	
	\title{Dense ideals in forcing extensions of Nairian models}
	\author[1]{Rahman Mohammadpour\thanks{This research is part of project No. 2022/47/P/ST1/00705,
			co-funded by the National Science Centre and the European Union's Horizon 2020 research
			and innovation programme under the Marie Sk\textlstroke{}odowska-Curie grant agreement No. 945339.}}
	\author[2]{Grigor Sargsyan\thanks{The second author's work is funded by Poland's National Science Centre under the Weave-UNISONO call, registration number UMO-2023/05/Y/ST1/00194.}}
	\affil[1,2]{\small Institute of Mathematics of the Polish Academy of Sciences, Jana i J\ogonekaccent{e}drzeja \'{S}niadeckich 8, 00-656 Warsaw, Poland}
	\date{}
	\maketitle
	
	\begin{abstract}
		Complementing the work of Blue, Larson, and Sargsyan~\cite{NM,NM2},
		we obtain forcing extensions of minimal Nairian models in which the Axiom of Choice and
		the Generalized Continuum Hypothesis hold. In these extensions, every
		uncountable regular cardinal $\k$ carries a precipitous, uniform,
		$\k$-dense ideal $\mathcal{I}_\k$. In fact,
		$\powerset(\k)/\mathcal{I}_\k$ is isomorphic to the regular open algebra
		of $\col{\omega}{\k}$.
		This extends Woodin's theorem on the existence of $\omega_1$-dense ideals
		on $\omega_1$ in forcing extensions of models of $\treg$
		\cite[Section~2.3]{ASTWZ24}, where $\treg$ denotes the theory
		$\zf+\ad_{\R}+\text{``$\Theta$ is regular''}$. Moreover, our work places the consistency strength of the aforementioned theory on generic embeddings below a Woodin cardinal that is a limit of Woodin cardinals, and also complements the seminal work of Eskew and Hayut~\cite{EH24}.
	\end{abstract}
	
	\section{Introduction}\label{sec:introduction}
	
	Interest in Nairian models comes from the desire to unify useful foundational
	frameworks such as Martin's Maximum, generic embeddings, determinacy axioms and other such frameworks. The idea is that any sufficiently strong theory of infinity should
	identify a core whose theory is the same across all such frameworks. We do not
	attempt to make this philosophical idea precise here, but we require that the
	parent theory can be forced over this core.
	
	The papers~\cite{NM,NM2} propose Nairian models as determinacy models
	that can play this role. We write $\treg$ for the theory $\zf+\ad_{\R}+\text{``$\Theta$ is regular''}$, where $\Theta$ is the least nonzero ordinal that is not a surjective image of the reals. We work in a parent model
	$V\models\treg+V=L(\powerset(\R))$. For each ordinal $\xi\leq \Theta$, we set
	\[
	I_\xi\df\bigcup_{\beta<\xi}\beta^\omega
	\qquad\text{and}\qquad
	N_\xi\df\hod_{I_\xi}\cap V_\xi,
	\]
	where both constructions are computed in $V$. We call $N_\xi$
	a \emph{Nairian model} if $N_\xi\models\zf$. 
	
	Recall that $\kappa$ is called a \emph{Solovay cardinal} if, for every
	$\alpha<\Theta^N$, there is no ordinal-definable surjection
	$f:\powerset(\alpha)\longrightarrow\kappa$.  Woodin showed that, under $\treg$, the
	Solovay cardinals form a closed sequence in $\Theta^V$.
	We call a Nairian model $N=N_\xi$ \emph{minimal} if the following
	conditions hold:
	\begin{enumerate}
		\item\label{itm:intro-minimal-nairian-theta}
		$\Theta^N$ is a Solovay cardinal of $V$.
		\item\label{itm:intro-minimal-nairian-limits}
		If $\delta$ is the least Solovay cardinal of $V$ greater than
		$\Theta^N$, then in $\hod^V$, $\xi$ is both a limit of
		${<}\delta$-strong cardinals and a limit of Woodin cardinals.
		\item\label{itm:intro-minimal-nairian-minimality}
		For every $\gamma\in(\Theta^N,\xi)$, $N_\gamma$ is not a Nairian model.
	\end{enumerate}
	
	The papers~\cite{NM,NM2} investigate whether Martin's Maximum
	($\mathsf{MM}$) can be forced over Nairian models. Although they do not establish
	this, they show that the failure of Jensen's square at all regular cardinals can
	be forced over Nairian models. A recent preprint by the second author with Blue
	and Trang shows that Todor\v{c}evi\'c's Strong Reflection Principle
	($\mathsf{SRP}$) for $\omega_3$  can also be forced over Nairian
	models. This extends Woodin's result that $\mathsf{SRP}(\omega_2)$ holds in
	$\mathbb{P}_{\max}$ extensions of models of $\treg+V=L(\powerset(\R))$.
	
	The current paper experiments with a different foundational framework. The goal here is to force the existence of generic embeddings over Nairian models, extending the following seminal theorem of Woodin.
	
	\begin{theorem}[Woodin; {\cite[Section~2.3, Corollary~2.15]{ASTWZ24}}]
		\label{thm:woodin-ch-dense-ideal}
		Assume $\treg$ and $V=L(\powerset(\R))$.
		Let $G\subseteq\col{\omega_1}{\R}$ be $V$-generic,
		and let
		\[
		H\subseteq\bigl(\col{\Theta^V}{\powerset(\R)^V}\bigr)^{V[G]}
		\]
		be $V[G]$-generic. Then $W=V[G][H]$ satisfies $\zfc+\ch$.
		In $W$, there is a precipitous, normal, uniform ideal
		$\mcl I$ on $\omega_1$ that is $\omega_1$-dense. More precisely, $(\powerset(\omega_1)/\mcl I)^+$
		has a dense subset of cardinality $\omega_1$ in $W$.
	\end{theorem}
	
	If $W$ is as in Theorem~\ref{thm:woodin-ch-dense-ideal} and
	$g\subseteq\col{\omega}{\omega_1^W}$ is $W$-generic, then, in $W[g]$,
	there is an elementary embedding $j:W\longrightarrow M$ into a transitive
	class $M$ such that $\crit(j)=\omega_1^W$ and $M^\omega\subseteq M$;
	see \cite[Fact~2.4(iii)]{ASTWZ24}.
	Axioms postulating the existence of such embeddings in forcing extensions
	are called \emph{generic large cardinal axioms}.
	For applications of generic large cardinal axioms, we refer the reader to
	Foreman's survey~\cite{ForemanGLC} and the work of Eskew and Hayut~\cite{EH24}.
	
	\begin{definition}\label{def:generic-supercompactness}
		Suppose that $\k$ is a regular cardinal and $\l$ is an ordinal.
		We say that $\k$ is \emph{generically $\l$-supercompact} if there is
		a forcing $\mbb Q\in V$ such that, for every $V$-generic filter
		$H\subseteq\mbb Q$, in $V[H]$ there are a transitive class $M$ and
		an elementary embedding $j:V\longrightarrow M$ satisfying the following conditions:
		\begin{enumerate}
			\item\label{itm:generic-supercompactness-critical-point} $\crit(j)=\k$.
			\item\label{itm:generic-supercompactness-image} $j[\l]\in M$.
			\item\label{itm:generic-supercompactness-closure} $M^\omega\subseteq M$, where
			the sequences are computed in $V[H]$.
		\end{enumerate}
	\end{definition}
	
	The model $W$ from Theorem~\ref{thm:woodin-ch-dense-ideal} satisfies $\zfc+\ch$.
	In $W$, the forcing $\col{\omega}{\omega_1}$ witnesses that $\omega_1$ is
	generically $\omega_1$-supercompact in the sense of
	Definition~\ref{def:generic-supercompactness}; see \cite[Fact~2.4(iii)]{ASTWZ24}.
	\begin{definition}\label{def:theory-tis}
		The theory $T_{\mathrm{IS}}$ consists of $\zfc$ together with the following assertion:
		\begin{quote}
			For every regular cardinal $\kappa\geq\omega_1$, there is a
			precipitous, fine, normal, uniform, $|\kappa^\omega|$-dense ideal $\mcl I$ on
			$\mcl P_{\omega_1}(\kappa^\omega)$ such that $\mcl I^+$ witnesses that
			$\omega_1$ is generically $\kappa$-supercompact in the sense of
			Definition~\ref{def:generic-supercompactness}.\footnote{See \cref{subsec:supercompactness-measures} for relevant notions.}
		\end{quote}
	\end{definition}
	
	We write $\hpc$ for Steel's Hod Pair Capturing hypothesis;
	see \cite[Section~1.7]{Steel23} for its formulation.
	Our main result is the following.
	
	\begin{namedtheorem}{Main Theorem}\label{thm:main-introduction}
		Assume that $V\mdls{\treg+\hpc}$ and $V=L(\mcl P(\R))$.
		If $N$ is a minimal Nairian model in $V$, then $N$ has a
		class forcing extension ${N}[G]$ such that
		\[
		{N}[G]\mdls{T_{\mathrm{IS}}+\gch}.
		\]
	\end{namedtheorem}
	
	The second author has shown that the hypotheses of the
	\cref{thm:main-introduction} are consistent relative to a Woodin cardinal
	that is a limit of Woodin cardinals. Therefore, we obtain the following corollary.
	
	\begin{corollary}\label{cor:consistency-tis}
		If $\zfc$ together with the assertion that there is a Woodin cardinal
		that is a limit of Woodin cardinals is consistent, then
		$T_{\mathrm{IS}}$ is consistent.
	\end{corollary}
	
	Our proof of the \cref{thm:main-introduction} follows Woodin's approach to
	\cref{thm:woodin-ch-dense-ideal}. We generalize weak comeagerness for
	subsets of $\R^\omega$ to subsets of $(\kappa^\omega)^\omega$;
	see \cref{def:weak-comeagerness} in \cref{subsec:weak-comeagerness}.
	Given a supercompactness measure $\mu$ on
	$\mcl P_{\omega_1}(\kappa^\omega)$ and a generic filter
	$g\subseteq\col{\omega}{\kappa^\omega}$, the ultrafilter $\mu/g$
	consists of the ground-model sets $A\subseteq(\kappa^\omega)^\omega$
	that are weakly comeager with respect to $\mu$ below some $p\in g$,
	as computed in the ground model. We construct $\mu/g$ in
	\cref{subsec:generic-elementary-embedding} and use $\mu$-Uniformization
	to establish its \Los{} property in \cref{prop:los-property}.
	The uniformization results for $\textrm{EnSub}_\kappa$ and
	$\textrm{EnSub}_\kappa\times\R$ are
	\cref{thm:measure-uniformization,thm:measure-real-uniformization} in
	\cref{sec:iteration-uniformization}. We prove the result for weakly
	comeager subsets of $(\textrm{EnSub}_\kappa)^\omega$ in
	\cref{thm:weak-comeager-uniformization} (\cref{sec:weak-comeager}).
	The version for measures on $\mcl P_{\omega_1}(\kappa_\alpha^\omega)$
	needed for the generic ultrapower is stated in
	\cref{thm:uniformization_in_minimal_N_model}.
	
	\Needspace{8\baselineskip}
	\section{Basics}\label{sec:preliminaries}
	We work in $\zf$ throughout this section.
	We fix certain basic set-theoretic notation in \cref{subsec:notation}. Then in  \cref{subsec:supercompactness-measures}, we state our notation and terminology concerning filters and ideals. In particular, we recall the definitions of
	supercompactness measures.
	We then introduce determinacy, Dependent Choice, and the $\Theta$-sequence in \cref{subsec:determinacy}.
	In \cref{subsec:basic-nairian-theory}, we define the $\hH$-hierarchy and the basic Nairian theory,
	and recall the Choice-forcing iteration of Blue, Larson, and Sargsyan \cite{NM2},
	including their analyses of cardinal arithmetic, Dependent Choice at intermediate stages,
	and closure of the tail forcings. Finally, in \cref{ssec:minimal-Nairian-models}, we state certain theorems about minimal Nairian models.
	
	\subsection{Notation}\label{subsec:notation}
	Suppose that $X$ and $Y$ are sets and that $f$ is a function from $X$ into $Y$.
	We let $f:X\twoheadrightarrow Y$ indicate that $f$ is a surjection. We let $\rang{f}$ denote the range of $f$. If $Z\subseteq\dom{f}$,  $f[Z]$ denotes the pointwise image of $Z$ under $f$. We  let $\card{X} = \card{Y}$ indicate that  there is a bijection between $X$ and $Y$.
	We let $\tc{X}$ denote the transitive closure of $X$. 
	If $X$ is a set of ordinals, $\otp{X}$ denotes the order type of $(X,\in)$. 
	If $\kappa$ is a cardinal, we let
	\[
	\mcl P_\kappa(X)\df\{Z\subseteq X:\exists\alpha<\kappa\ \text{ such that }|Z|=|\alpha|\}.
	\]
	
	We let $X^Y$ denote the set of functions from $Y$ into $X$. For an ordinal $\beta$, we let
	\[
	X^{<\beta}\df\bigcup_{\alpha<\beta}X^\alpha.
	\]
	
	For $p,q\in X^{<\beta}$, $p\leq q$ denotes $p\supseteq q$.

	If $X\neq \emptyset$ and $\kappa$ is an infinite cardinal, we let $\col{\kappa}{X}\df (X^{<\kappa},\supseteq)$ denote the familiar collapse forcing.  
	A $V$-generic filter $G$ on $\col{\kappa}{X}$ induces a natural surjection from $\kappa$ onto $X$, which we often call a \emph{generic coding} of $X$. We use the convention that stronger conditions are smaller in the forcing order.
	For a poset $\mbb P$, $\ro(\mbb P)$ denotes the regular open completion of  $\mbb P$.

	The following is standard.
	\begin{lemma}\label{lem:zf-collapse-absorption}
		Let $X$ be an infinite set admitting a bijection $X\times2\to X$.
		Let $\mbb P$ be a forcing notion, and suppose
		that $c:X\twoheadrightarrow\mbb P$ is a surjection. Then $\mbb P\times\col{\omega}{X}$ and $\col{\omega}{X}$ have isomorphic dense suborders.
	\end{lemma}

	\subsection{Filters, ideals, and generic supercompactness}\label{subsec:supercompactness-measures}
	Let $\mcl F$ be a filter on a set $X$. We shall denote the collection of $\mcl F$-positive sets by $\mcl F^+$, where a set $A\subseteq X$ is $\mcl F$-positive if for every $B\in\mcl F$, $A\cap B\neq\emptyset$. We will use the standard notion of \emph{normality} in the absence of Choice, which states that whenever $U\in\mcl F^+$ and
	$F$ is a function on $U$ such that 
	$\emptyset\neq F(\sigma)\subseteq\sigma$ for every $\sigma\in U$, there are $U'\subseteq U$ with $U'\in\mcl F^+$ and an element $a$ such that
	$a\in F(\sigma)$ for every $\sigma\in U'$.
	Recall also that if $X=\mcl P_\kappa(Z)$ for some cardinal $\kappa$ and some set $Z$, we say $\mcl F$ is \emph{fine} if for every $z\in Z$, $\{A\in\mcl P_\kappa(Z):z\in A\}\in\mcl F$.
	
	The following easy lemma will be used in the proof of a uniformization theorem in \cref{sec:uniformization-in-nairian-models}.
	
	\begin{lemma}\label{lem:directed-supports}
		Let $X$ be an infinite set, and let $Y\subseteq X^\omega$. Assume that $\{\rang{y}:y\in Y\}$ is directed under inclusion.
		Let $\mcl F$ be a normal fine  filter on $\mcl P_{\omega_1}(Y)$.
		Then
		\[
		\left\{\sigma\in\mcl P_{\omega_1}(Y):
		\begin{gathered}
			\forall y,y'\in\sigma~\exists y''\in\sigma\\
			\rang{y}\cup\rang{y'}\subseteq\rang{y''}
		\end{gathered}\right\}\in\mcl F^+.
		\]
		
	\end{lemma}
	\begin{proof}
		Suppose, toward a contradiction, that the conclusion fails. By a double application
		of normality to the complement of the above set, there are $y,y'\in Y$ such that
		\[
		A\df\{\sigma\in\mcl P_{\omega_1}(Y):
		\forall y''\in\sigma[
		\rang{y}\cup\rang{y'}\nsubseteq\rang{y''}]
		\}\in\mcl F^+.
		\]
		By the directedness assumption, there is $z\in 
		Y$ with
		$\rang{y}\cup\rang{y'}\subseteq\rang{z}$, and by
		fineness,
		$\{\sigma\in\mcl P_{\omega_1}(Y):z\in\sigma\}\in\mcl F$. But $A\cap \{\sigma\in\mcl P_{\omega_1}(Y):z\in\sigma\}=\emptyset$, a contradiction!
	\end{proof}
	
	\begin{lemma}\label{lem:fine-and-complete}
		Let $X$ be an uncountable set.   Assume that $\mcl F$ is a fine, countably complete filter on $\mcl P_{\omega_1}(X)$. Let $a\in\mcl P_{\omega_1}(X)$. Then 
		$\{\sigma\in\mcl P_{\omega_1}(X): a\subseteq\sigma\}\in\mcl F$.
	\end{lemma}
	\begin{proof}
		We may assume that $a\neq\emptyset$. 
		Fix a surjection $s:\omega\twoheadrightarrow a$, and, for each $n<\omega$, let
		$C_n\df\{\sigma\in\mcl P_{\omega_1}(X):s(n)\in\sigma\}$.
		Notice that by fineness, each $C_n\in\mcl F$, and that
		\[
		\{\sigma\in\mcl P_{\omega_1}(X): a\subseteq\sigma\}=\bigcap_{n<\omega} C_n. 
		\]
		So by  countable completeness, we have
		$\{\sigma\in\mcl P_{\omega_1}(X): a\subseteq\sigma\}\in\mcl F$.
	\end{proof}
	\begin{definition}
		Let $\kappa$ be a cardinal and let $X$ be a set. A supercompactness measure on $\mcl P_\kappa(X)$ is a normal, fine, $\kappa$-complete, nonprincipal ultrafilter.
	\end{definition}
	
	\begin{remark}
		In applications, we may restrict a supercompactness measure on $\mcl P_{\omega_1}(X)$ to
		countably infinite $\sigma$'s, without further explanation.  
	\end{remark}
	\begin{definition}
		An uncountable cardinal  $\kappa$ is called \emph{supercompact} if, for every uncountable set $X$ with $X\notin\mcl P_\kappa(X)$, there is a \emph{supercompactness measure} on $\mcl P_\kappa(X)$.
	\end{definition}
	
	\begin{lemma}\label{lem:normality-fin}
		Let $\mu$ be a supercompactness measure on $\mcl P_{\omega_1}(X)$.
		Suppose that $F$ is a function on $U\in\mu$ such that
		$\emptyset\neq F(\sigma)\subseteq\sigma^{<\omega}$ for every $\sigma\in U$.
		Then there is $U'\subseteq U$ with $U'\in\mu$ and there is $q\in X^{<\omega}$
		such that $q\in F(\sigma)$ for every $\sigma\in U'$.
	\end{lemma}
	\begin{proof}
		Let us first prove the lemma under the additional assumption that there is a positive integer
		$\iota$ such that, for every $\sigma\in U$ and every $q\in F(\sigma)$,
		$|q|=\iota$. The case $\iota=1$ follows directly from the definition of normality. Suppose
		$\iota>1$. Let $F'$ be defined on $U$ by setting
		\[
		F'(\sigma)\df\{q\restriction(\iota-1):q\in F(\sigma)\}.
		\]
		By induction, there are $U'\subseteq U$ in $\mu$ and
		$q^*\in X^{\iota-1}$ such that $q^*\in F'(\sigma)$ for every
		$\sigma\in U'$. Now define  $F'':U'\to \mcl P(X)$ by
		\[
		F''(\sigma)\df\{q(\iota-1):q\in F(\sigma)\text{ and }q\supseteq q^*\}.
		\]
		Then $F''(\sigma)$ is a nonempty subset of $\sigma$
		for every $\sigma\in U'$. By normality, there are $U''\subseteq U'$ in
		$\mu$ and $a\in X$ such that $a\in F''(\sigma)$ for every
		$\sigma\in U''$. Let $q\in X^\iota$ be the unique extension of $q^*$
		such that $q(\iota-1)=a$. Then $q\in F(\sigma)$ for every
		$\sigma\in U''$, as required.
		
		For the general case, by countable completeness there are
		$U^*\subseteq U$ in $\mu$ and $\iota<\omega$ such that, for every
		$\sigma\in U^*$, $F(\sigma)$ contains an $\iota$-tuple. If $\iota=0$,
		then $\emptyset\in F(\sigma)$ for every $\sigma\in U^*$, and we are done. Otherwise, let $F^*(\sigma)\df\{q\in F(\sigma):|q|=\iota\}$.
		Applying the previous argument to $F^*$ and $U^*$ completes the proof.
	\end{proof}
	
	\begin{notation}\label{not:i_mu}
		Let $\mu$ be an ultrafilter. The canonical $\mu$-ultrapower embedding  $i_\mu:V\to\ult{V}{\mu}$ is defined by
		$x\mapsto[\lambda z.~x]_\mu$. Here $\lambda z.~x$ denotes the constant function with value $x$.
	\end{notation}

	\begin{definition}\label{def:los-property}
		Let $M$ be a transitive set or proper class satisfying $\zf$.
		Let $X\in M$, and let $\mu\subseteq \mcl P(X)\cap M$ be an ultrafilter.   We say that $\mu$ has the \emph{\Los{} property over $M$} if, for every formula $\varphi$ and all functions
		$f_0,\ldots,f_n:X\to M$ belonging to $M$,
		\[
		\ult{M}{\mu}\models{
			\varphi([f_0]_\mu,\ldots,[f_n]_\mu)}
		\quad\Longleftrightarrow\quad
		\{x\in X:M\models{
			\varphi(f_0(x),\ldots,f_n(x))}\}\in\mu.
		\]
	\end{definition}
	
	Note that under $\zfc$ every ultrafilter has the \Los{} property over the
	universe.

	For an ideal $\mcl I\subseteq\mcl P(X)$, we denote the collection of $\mcl I$-positive sets by $\mcl I^+$, where a set $A\subseteq X$ is $\mcl I$-positive if $A\notin\mcl I$. Notice that completeness, fineness, and normality for $\mcl I$ are naturally defined through its dual filter $\{X\setminus A:A\in\mcl I\}$. For $A,B\in\mcl I^+$, we let
	\[
	A\leq_{\mcl I}B\quad\Longleftrightarrow\quad A\setminus B\in\mcl I.
	\]
	Note that $(\mcl P(X)/\mcl I)\setminus\{0\}$ is the separative quotient  of the forcing  $(\mcl I^+,\leq_{\mcl I})$.
	Forcing with  $(\mcl I^+,\leq_{\mcl I})$ yields a generic ultrafilter $G$ on $\mcl P(X)^V$.  Then in $V[G]$, the ultrapower embedding is defined over $V$ as usual:
	\[
	i_G:V\to\ult{V}{G}.
	\] 
	Recall that $\mcl I$ is called \emph{precipitous} if for every $V$-generic filter $G$ on $\mcl I^+$,  $\ult{V}{G}$ is well-founded.
	
	\begin{definition}
		Assume $\ac$. Let $\kappa$ be an uncountable regular cardinal.  
		We say $\kappa$ is 
		\emph{generically  supercompact} if it is generically $\lambda$-supercompact for every ordinal $\lambda$.\footnote{See \cref{def:generic-supercompactness} for the definition of generic $\lambda$-supercompactness.}
	\end{definition}
	
	It is clear that $\kappa$ is generically supercompact  if it is generically $\lambda$-supercompact for an unbounded class of  ordinals $\lambda$.

	\begin{definition}\label{def:dense-ideal}
		Assume that $\mcl I$ is an ideal. Let $\kappa$ be a cardinal.
		We say
		that 
		$\mcl I$ is \emph{$\kappa$-dense} if
		$(\mcl I^+,\leq_{\mcl I})$ has a dense subset of cardinality at most
		$\kappa$.
	\end{definition}
	
	\begin{definition}\label{def:uniform-ideal}
		Let $\kappa$ be an infinite cardinal, and let $Z$ be a set with $|Z|=\kappa$.
		An ideal $\mcl I$ on  $Z$ is called \emph{uniform} if
		$\mcl P_\kappa(Z)\subseteq\mcl I$.
	\end{definition}
	Thus every uniform ideal is nonprincipal. When $Z=\kappa$ and $\kappa$ is regular, uniformity is
	equivalent to extending the ideal of bounded subsets.
	
	\begin{definition}[Wilson~\cite{Wilson-thesis}]\label{def:ordinal-covering}
		Let $W$ be a transitive model of $\zf$. Let  $\mcl I\in W$ be an ideal on $Z$. 
		Assume that $V\subseteq W$ is an inner model containing $Z$ as an element.
		We say that $\mcl I$
		has the \emph{ordinal covering property relative to $V$} if, for every
		ordinal $\eta$, every $F:Z\to\eta$, and every $S\in\mcl I^+$,
		there are $F_0:Z\to\eta$ in $V$ and $S'\in\mcl I^+$ with
		$S'\subseteq S$ such that $F\restriction S'=F_0\restriction S'$.
	\end{definition}
	
	\subsection{Determinacy}\label{subsec:determinacy}
	Let $X$ be a set and let $A\subseteq X^\omega$. Consider the two-player, perfect-information game of length $\omega$ with payoff set $A$, denoted by $\Game(A,X)$. The players $\mathrm{P}_{\mathrm{I}}$ and $\mathrm{P}_{\mathrm{II}}$ alternate playing elements of $X$, with $\mathrm{P}_{\mathrm{I}}$ moving first. At every stage, both players know all previous moves. A run of the game produces an element $x\in X^\omega$, and $\mathrm{P}_{\mathrm{I}}$ wins if and only if $x\in A$.
	
	A strategy for a player assigns a move to every finite position at which it is that player's turn. A strategy is winning if every run consistent with it is won by that player. The game $\Game(A,X)$ is called determined if one of the players has a winning strategy; in this case, we also say that $A$ is determined.
	
	The axiom $\ad_X$, called the \emph{Axiom of $X$-Determinacy}, asserts that every $A\subseteq X^\omega$ is determined. We refer to $\ad_{\R}$ as the \emph{Axiom of Real Determinacy}, while $\ad$ denotes the usual \emph{Axiom of Determinacy}, i.e., $\ad_\omega$. Here, $\R$ denotes the Baire space $\omega^\omega$.
	
	Let $\kappa$ be an infinite cardinal. Recall that the principle of \emph{$\kappa$-Dependent Choice for $X$}, denoted by
	$\dc_\kappa(X)$, asserts that every $\kappa$-closed transitive binary
	relation on $X$ has a $\kappa$-chain. The principle of
	\emph{$\kappa$-Dependent Choice}, denoted by $\dc_\kappa$, asserts that
	$\dc_\kappa(X)$ holds for every set $X$.\footnote{\label{fn:dependent-choice-chains}
		Suppose that $T$ is a transitive binary relation on a set $X$.  Given an ordinal $\eta$, an $\eta$-chain through $T$ is a function $f \colon \eta \to X$ such that $(f(\alpha), f(\beta)) \in T$ for all $\alpha < \beta < \eta$. 
		A $T$-upper bound for an $\eta$-chain $f$ is an element $u\in X$ such
		that $(f(\alpha),u)\in T$ for all $\alpha<\eta$. We say that $T$ is $\kappa$-closed if $f$ has a $T$-upper bound whenever $\eta<\kappa$ and $f\colon \eta\rightarrow X$ is an $\eta$-chain through $T$.} Finally, $\dc$ will stand for $\dc_{\omega}$. Notice that the Axiom of Choice is equivalent to the statement that for every infinite cardinal $\kappa$, $\dc_\kappa$ holds.
	
	For every transitive class $M$ and every infinite set $X\in M$, we let
	\[
	\Theta_X^M\df\sup\{\alpha \mid \exists f\in M, f: X\twoheadrightarrow \alpha \}. 
	\]
	We shall omit the superscript $M$ when $M=V$, unless emphasis is needed.
	For a transitive model $M$ of $\zf$, we shall denote $\Theta_{\R^M}^M$ by $\Theta^M$, following the standard convention, where $\R^M$ is the set of reals of $M$.
	
	Notice that for infinite $X$, every nonzero ordinal below $\Theta_X$ is a surjective image of $X$, whereas $\Theta_X$ is not.
	
	\begin{lemma}\label{lem:reg-cof}
		Let $X$ be an infinite set. Then
		$\Theta_X$ is regular if and only if there is no cofinal function $f: X \to \Theta_X$.
	\end{lemma}
	\begin{proof}
		$(\Rightarrow)$ Assume $\Theta_X$ is regular. Suppose, for a contradiction, that there is a cofinal function $f: X \to \Theta_X$. Let us set $\gamma \df \otp{\rang{f}}$.
		We must have $\gamma < \Theta_X$, for otherwise $f$ would yield a surjection from $X$ onto $\Theta_X$, which contradicts the definition of $\Theta_X$.
		Since $\rang{f}$ is cofinal in $\Theta_X$, we have $\cof{\Theta_X} \le \gamma < \Theta_X$. Thus $\Theta_X$ is singular, contradicting the assumption.
		
		$(\Leftarrow)$ Assume there is no cofinal function from $X$ into $\Theta_X$. Suppose, for a contradiction, that $\Theta_X$ is singular. Then $\kappa \df \cof{\Theta_X} < \Theta_X$. By the definition of $\Theta_X$, there is a surjection $\pi: X \twoheadrightarrow \kappa$. Let $c: \kappa \to \Theta_X$ be a cofinal function witnessing the singularity. Then the function $c \circ \pi: X \to \Theta_X$ is a cofinal function, a contradiction.
	\end{proof}

	\begin{definition}\label{def:theta-sequence}
		We define the $\Theta$-sequence $\langle \kappa_\alpha:\alpha\in\ord\rangle$ recursively. Set $\kappa_0\df\Theta$ and, for every ordinal $\alpha>0$, define
		\[
		\kappa_\alpha\df\sup\{\Theta_{\kappa_{\beta}^\omega}:\beta<\alpha\}.
		\]
		
		We also adopt the convention $\kappa_{-1}=\omega$.
		Note that for every ordinal $\alpha$, $\kappa_{\alpha+1}=\Theta_{\kappa_\alpha^\omega}$.
	\end{definition}
	
	\subsection{The basic Nairian theory and forcing Choice}\label{subsec:basic-nairian-theory}
	
	\begin{definition}\label{def:maximal-cardinality}
		Let $M$ be a nonempty set.  A set $X\in M$ is a \emph{maximal cardinality} of $M$ if for every nonempty $Y\in M$, there exists a surjection from $X$ to $Y$  in $M$. 
	\end{definition}
	\begin{definition}\label{def:has-maximal-successor} 
		Let $M$ be a nonempty set.
		We say $M$ \emph{has a maximal cardinality} if some $X\in M$ is a maximal cardinality of $M$. In this case, we say
		$M$ is of \emph{successor type}.
	\end{definition}

	\begin{definition}\label{def:full}
		A nonempty set $M$ is called \emph{full} if it satisfies the following three conditions.
		\begin{enumerate}
			\item Transitively closed: for every $Y\in M$, there is a transitive set $Z\in M$ such that $Z\subseteq M$ and $Y\in Z$.
			\item Closed under products: for all $X,Y\in M$, $X\times Y\in M$.
			\item Closed under subsets: for every $X\in M$, $\powerset(X)\subseteq M$.
		\end{enumerate}
	\end{definition}
	
	\begin{definition}\label{def:ordinal-fullness}
		Let $M$ and $X$ be sets. We say that $M$ is \emph{$X$-full} if
		$Y^X\in M$ for every $Y\in M$, where the function set is computed
		in the ambient universe. We say that $M$ is \emph{completely ordinal-full}
		if it is $\eta$-full whenever $\eta\in M$ is an ordinal below some
		cardinal as computed in $M$.
	\end{definition}
	
	\begin{definition}\label{def:strongly-regular}
		A nonempty set $M$ is called \emph{strongly regular} if, whenever $X\in M$ and $f$ is a function from $X$ to $M\cap\ord$, the range of $f$ is bounded in $M\cap\ord$.
	\end{definition}
	
	\begin{definition}\label{def:hierarchical}
		Let $M$ be a transitive set, and let $\kappa=M\cap\ord$. Then $M$ is called \emph{hierarchical} if there is a continuous, $\subseteq$-increasing sequence of transitive sets
		$\langle M_\alpha:\alpha<\kappa\rangle$
		such that $M=\bigcup_{\alpha<\kappa}M_\alpha$
		and $M_\alpha\in M$ for every $\alpha<\kappa$. 
	\end{definition}

	\begin{definition}\label{def:H-like}
		Suppose that $M$ is a transitive set. Then $M$ is \emph{$H$-like} if
		\begin{enumerate}
			\item\label{itm:h-like-full-and-hierarchical} $M$ is full and hierarchical, and
			\item\label{itm:h-like-successor-regularity} if $M$ is of successor type, then $M$ is strongly regular.
		\end{enumerate}
	\end{definition}
	We now introduce  the \emph{$\mcl H$-set} associated with a set $X$.
	
	\begin{notation}\label{not:h-set}
		Let $X$ be a set. We set
		\[
		\mcl H_X\df\{Y:\exists f\colon X\twoheadrightarrow\tc{\{Y\}}\}.
		\]
	\end{notation}
	
	For an infinite cardinal $\kappa$, $\mcl H_\kappa$ consists of sets of hereditary
	cardinality at most $\kappa$; under Choice, this is the usual $H_{\kappa^+}$.
	The sets $\mcl H_X$ will be used to define a hierarchy analogous to the usual $H_\kappa$-hierarchy under $\zfc$.
	
	\begin{proposition}[Blue, Larson, Sargsyan~{\cite{NM2}}]\label{prop:bls-h-like-preservation}
		Suppose that $Z$ is an $H$-like set, $X\in Z$ is transitive and a maximal
		cardinality of $Z$, and $\mathbb P\in Z$ is a forcing notion.
		If $G\subseteq\mathbb P$ is $V$-generic, then
		\[
		Z[G]=\mcl H_X^{V[G]},
		\]
		and $Z[G]$ is an $H$-like set of successor type in $V[G]$,
		with $Z[G]\cap\ord=Z\cap\ord$. In particular,
		$Z\cap\ord$ is regular in $V[G]$.
	\end{proposition}
	\begin{proof}
		See \cite[Theorem~4.13 and its proof]{NM2}.
	\end{proof}
	
	\begin{definition}[$\hH$-hierarchy]\label{def:nairian-hierarchy}
		The sequence $\langle\hH_\alpha:\alpha\in\ord\rangle$ is defined as follows.
		\begin{enumerate}
			\item\label{itm:nairian-hierarchy-zero} $\hH_0\df \mcl H_{\powerset(\omega)}$.
			
			\item\label{itm:nairian-hierarchy-one} $\hH_1\df \mcl H_{\hH_0}$.
			
			\item\label{itm:nairian-hierarchy-successor} For every ordinal $\alpha$,
			$	\hH_{\alpha+2}\df \mcl H_{\hH_{\alpha+1}}$.
			
			\item\label{itm:nairian-hierarchy-limit} For every nonzero limit ordinal $\alpha$,
			\[
			\hH_\alpha\df \bigcup_{\beta<\alpha}\hH_\beta
			\qquad\text{and}\qquad
			\hH_{\alpha+1}\df \mcl H_{\powerset(\hH_\alpha)}.
			\]
		\end{enumerate}
	\end{definition}
	
	It follows directly from the definition that each $\hH_\alpha$ is transitive and the sequence $\langle\hH_\alpha:\alpha\in\ord\rangle$ is $\subseteq$-increasing and $\subseteq$-continuous.  
	
	\begin{notation}
		Let $\alpha$ be an ordinal.
		We set $\eta_\alpha\df\hH_\alpha\cap\ord$.
		
	\end{notation}
	
	\begin{definition}[Axiom of Harmony]\label{ax:axiom-of-harmony}
		The \emph{Axiom of Harmony} asserts that $\hH_\alpha$ is $H$-like for every ordinal $\alpha$.
	\end{definition}

	\begin{definition}\label{def:x-closed}
		Let $X$ and $W$ be sets. We say that $W$ is \emph{$X$-closed} if
		$X\cup\{X\}\cup W^X\subseteq W$.
	\end{definition}
	
	\begin{definition}[Reflection]\label{def:reflection}
		For sets $X$, $Y$, and $Z$, the principle $\refl(X,Y,Z)$ is the statement that whenever $W$ is a transitive $X$-closed structure with $Y\in W$ and $a\in W$, there exist a transitive set $\overline{W}$ and an elementary embedding $\pi\colon \overline{W}\to W$ such that
		\begin{itemize}
			\item $X,Y\in \overline{W}$,
			\item $a\in\rang{\pi}$,
			\item $\overline{W}\in Z$,
			\item $\pi\rest\bigl(X\cup Y\cup\{X,Y\}\bigr)
			=\id{X\cup Y\cup\{X,Y\}}$,
			and
			\item $\overline{W}$ is $X$-closed.
		\end{itemize}
	\end{definition}
	
	\begin{definition}[Reflection Principle]\label{ax:reflection-principle}
		The \emph{Reflection Principle} asserts that for every nonlimit ordinal $\alpha$ and every maximal cardinality $X$ of $\hH_\alpha$, $\refl(X,\hH_\alpha,\hH_{\alpha+1})$ holds.
	\end{definition}
	
	\begin{definition}[$\ntbase$]\label{def:basic-nairian-theory}
		The \emph{basic Nairian theory}, denoted by $\ntbase$, consists of $\zf$ together with the Axiom of Harmony and the Reflection Principle.
	\end{definition}

	For the rest of this section, we work in a transitive model $\mathsf{M}$ of $\ntbase$, and write
	$\hH_\gamma$ for the levels of the $\hH$-hierarchy computed in $\mathsf{M}$.
	We now recall the forcing construction used by Blue, Larson, and Sargsyan \cite{NM2} to force Choice over
	$\mathsf{M}$.

	\begin{theorem}[Blue, Larson, Sargsyan~{\cite{NM2}}]\label{thm:forcing-choice}	Suppose that $\eta_\epsilon$ is singular in $\mathsf{M}$ for every nonzero
		limit ordinal $\epsilon$.
		Let $\alpha$ be $0$ or a successor ordinal, let $\mathbb P\in\hH_\alpha$
		be a forcing notion, and let $G\subseteq\mathbb P$ be $\mathsf{M}$-generic.
		Suppose that $\hH_\alpha[G]$ satisfies Choice and is completely
		ordinal-full in $\mathsf{M}[G]$.
		
		In $\mathsf{M}[G]$, let $\mathbb Q^{\textsf{ac}}\df\mathbb Q^{\textsf{ac}}(\alpha,G)$ be the full-support class iteration
		of length $\ord^{\mathsf{M}}$ with $\mathbb Q^{\textsf{ac}}_{<0}\df\{1\}$ and
		\[
		\mathbb Q^{\textsf{ac}}_{<\beta+1}\df\mathbb Q^{\textsf{ac}}_{<\beta}\ast\dot{\mathbb Q}^{\textsf{ac}}(\beta).
		\]
		The initial segment $\mathbb Q^{\textsf{ac}}_{<\beta}$ at a limit ordinal $\beta$
		is taken with full support. Set
		\[
		\mathbb Q^{\textsf{ac}}(0)\df\col{\eta_\alpha}{\hH_\alpha}^{\mathsf{M}[G]}.
		\]
		At each  stage $\beta>0$, let $H_\beta\subseteq\mathbb Q^{\textsf{ac}}_{<\beta}$
		be generic over $\mathsf{M}[G]$. In $\mathsf{M}[G][H_\beta]$, define the next forcing
		as follows:
		\begin{enumerate}
			\item If $\beta$ is a  limit, $\mathbb Q^{\textsf{ac}}(\beta)$ is trivial.
			\item If $\beta=1$ or $\beta$ is the successor of a successor ordinal, let
			\[
			\mathbb Q^{\textsf{ac}}(\beta)\df\col{\eta_{\alpha+\beta}}{\hH_{\alpha+\beta}}.
			\]
			\item If $\beta=\delta+1$ for a nonzero limit ordinal $\delta$, put
			$\upsilon=(\eta_{\alpha+\delta}^+)^{\mathsf{M}[G][H_\beta]}$ and let
			\[
			\mathbb Q^{\textsf{ac}}(\beta)\df
			\col{\upsilon}{\mcl P(\hH_{\alpha+\delta})^{\mathsf{M}}}
			\ast\dotcol{\eta_{\alpha+\delta+1}}{\hH_{\alpha+\delta+1}}.
			\]
		\end{enumerate}
		
		Then, for every $\mathsf{M}[G]$-generic filter $H\subseteq\mathbb Q^{\textsf{ac}}$,
		the extension $\mathsf{M}[G][H]$ satisfies $\zfc$, and in $\mathsf{M}[G][H]$, for every ordinal $\beta$,
		\begin{enumerate}
			\item\label{itm:forcing-choice-cardinals-gch} $\eta_{\alpha+\beta}$ remains a cardinal and
			$2^{\eta_{\alpha+\beta}}=\eta_{\alpha+\beta}^+$, and
			\item\label{itm:forcing-choice-successor-heights} we have
			\[
			\eta_{\alpha+\beta+1}=
			\begin{cases}
				\eta_{\alpha+\beta}^+,
				&\text{if $\beta=0$ or $\beta$ is a successor},\\
				\eta_{\alpha+\beta}^{++},
				&\text{if $\beta$ is a nonzero limit}.
			\end{cases}
			\]
		\end{enumerate}
		
	\end{theorem}
	\begin{proof}
		See~{\cite[Theorem~8.7]{NM2}}.
	\end{proof}
	For the rest of this subsection, assume the hypotheses and notation
	of \cref{thm:forcing-choice}. Let $H\subseteq\mathbb Q^{\textsf{ac}}$ be an $\mathsf{M}[G]$-generic filter,
	and write $H_\gamma$ for its restriction to $\mathbb Q^{\textsf{ac}}_{<\gamma}$.
	We write $H(\gamma)$ for the induced generic filter at coordinate $\gamma$.
	
	\begin{notation}\label{not:factors}
		Let $\beta$ be either $-1$ or an ordinal. Write
		\[
		\mathbb Q^{\textsf{ac}}(\beta+1)=\mathbb Q^{\textsf{ac}}(\beta+1)_0
		\ast\dot{\mathbb Q}^{\textsf{ac}}(\beta+1)_1,
		\]
		where the first factor, computed in $\mathsf{M}[G][H_{\beta+1}]$, is
		\[
		\mathbb Q^{\textsf{ac}}(\beta+1)_0\df
		\begin{cases}
			\{1\}, & \text{if $\beta=-1,0$ or $\beta$ is a successor},\\
			\col{\eta_{\alpha+\beta}^+}{\mcl P(\hH_{\alpha+\beta})^{\mathsf{M}}},
			& \text{if $\beta$ is a nonzero limit}.
		\end{cases}
		\]

		Let $H(\beta+1)_0$ and $H(\beta+1)_1$ be the generic filters naturally induced by $H(\beta+1)$ via this factorization. In
		$\mathsf{M}[G][H_{\beta+1}][H(\beta+1)_0]$, the second factor is
		\[
		\mathbb Q^{\textsf{ac}}(\beta+1)_1
		\df\col{\eta_{\alpha+(\beta+1)}}{\hH_{\alpha+(\beta+1)}}.
		\]
	\end{notation}
	Note that
	$\mathbb Q^{\textsf{ac}}(0)_1=\mathbb Q^{\textsf{ac}}(0)$.
	
	\begin{proposition}[Blue, Larson, Sargsyan~\cite{NM2}]\label{prop:bls-bounded-presentations}
		Let $\beta$ be an ordinal.
		\begin{enumerate}
			\item\label{itm:bls-successor-presentation}
			If $\beta$ is a successor or $\beta=0$, then $\mathbb P\ast\dot{\mathbb Q}^{\textsf{ac}}_{<\beta}\in\hH_{\alpha+\beta}$.
			\item\label{itm:bls-limit-presentation}
			If $\beta$ is a nonzero limit, then $\mathbb P\ast\dot{\mathbb Q}^{\textsf{ac}}_{<\beta+1}
			\ast\dot{\mathbb Q}^{\textsf{ac}}(\beta+1)_0\in\hH_{\alpha+\beta+1}$.
			
		\end{enumerate}
	\end{proposition}
	
	\begin{proof}
		See~{\cite[proof of Theorem~8.7]{NM2}}.
	\end{proof}
	
	\begin{proposition}[Blue, Larson, Sargsyan~{\cite{NM2}}]\label{prop:bls-stage-choice}
		Let $\theta$ be the largest cardinal of $\hH_\alpha[G]$. Then
		\begin{enumerate}
			\item\label{itm:bls-initial-dc}
			$\mathsf{M}[G]\models{\dc_\theta}$.
			Moreover, $\theta$ remains a cardinal in $\mathsf{M}[G][H]$, and
			\[
			\mathsf{M}[G][H]\mdls{2^\theta=\eta_\alpha=\theta^+}.
			\]
			\item\label{itm:bls-completed-stage}
			If $\beta=0$ or $\beta$ is a successor ordinal, then
			$\hH_{\alpha+\beta+1}[G][H_{\beta+1}]$ satisfies Choice and has
			largest cardinal $\eta_{\alpha+\beta}$, and
			$\mathsf{M}[G][H_{\beta+1}]\models{\dc_{\eta_{\alpha+\beta}}}$.
		\end{enumerate}
	\end{proposition}
	\begin{proof}
		See~{\cite[Claims~8.9--8.11 and~8.13]{NM2}}.
	\end{proof}
	\begin{proposition}[Blue, Larson, Sargsyan~{\cite{NM2}}]\label{prop:bls-limit-stages}
		Let $\delta$ be a nonzero limit ordinal.\footnote{Notice that since coordinate $\delta$ is trivial,
			$\mathsf{M}[G][H_\delta]=\mathsf{M}[G][H_{\delta+1}]$.
			Adjoining $H(\delta+1)_0$ and then $H(\delta+1)_1$ gives
			$\mathsf{M}[G][H_{\delta+2}]$.}
		\begin{enumerate}
			\item\label{itm:bls-limit-size}
			In $\mathsf{M}[G][H_\delta]$, $\eta_{\alpha+\delta}$ is a cardinal,
			$\dc_{\eta_{\alpha+\delta}}$ holds,
			and there is a bijection between $\eta_{\alpha+\delta}$ and
			$\hH_{\alpha+\delta}[G][H_\delta]$.
			Moreover, $(\eta_{\alpha+\delta}^+)^{\mathsf{M}[G][H_\delta]}$ is regular
			in $\mathsf{M}[G][H_\delta]$, and
			\[
			\begin{gathered}
				\mcl P(\eta_{\alpha+\delta})^{\mathsf{M}[G][H_\delta]}
				=\mcl P(\eta_{\alpha+\delta})^{\mathsf{M}[G][H]},\\
				(\eta_{\alpha+\delta}^+)^{\mathsf{M}[G][H_\delta]}
				=(\eta_{\alpha+\delta}^+)^{\mathsf{M}[G][H]}.
			\end{gathered}
			\]
			\item\label{itm:bls-first-limit-factor}
			The model $\hH_{\alpha+\delta+1}[G][H_{\delta+1}][H(\delta+1)_0]$
			satisfies Choice and has largest cardinal
			$(\eta_{\alpha+\delta}^+)^{\mathsf{M}[G][H_\delta]}$.
			In $\mathsf{M}[G][H_{\delta+1}][H(\delta+1)_0]$, we have $\dc_{\eta_{\alpha+\delta}^+}$ and
			\[
			2^{\eta_{\alpha+\delta}}=\eta_{\alpha+\delta}^+,\qquad\text{and}\qquad
			\eta_{\alpha+\delta}^{++}=\eta_{\alpha+\delta+1}.
			\]
			Moreover, $\eta_{\alpha+\delta}^+$ and $\eta_{\alpha+\delta+1}$ are regular in $\mathsf{M}[G][H_{\delta+1}][H(\delta+1)_0]$.
			\item\label{itm:bls-second-limit-factor}
			The model $\hH_{\alpha+\delta+2}[G][H_{\delta+2}]$
			satisfies Choice and has largest cardinal $\eta_{\alpha+\delta+1}$.
			In $\mathsf{M}[G][H_{\delta+2}]$, we have $\dc_{\eta_{\alpha+\delta+1}}$, and
			\[
			2^{\eta_{\alpha+\delta}}=\eta_{\alpha+\delta}^+,\qquad\text{and}\qquad
			2^{\eta_{\alpha+\delta}^+}=\eta_{\alpha+\delta}^{++}=\eta_{\alpha+\delta+1}.
			\]
			Moreover, both $\eta_{\alpha+\delta}^+$ and $\eta_{\alpha+\delta+1}$ are regular in $\mathsf{M}[G][H_{\delta+2}]$.
		\end{enumerate}
	\end{proposition}
	\begin{proof}
		See~\cite[Claim~8.12 and proof of Claim~8.13]{NM2}.
	\end{proof}
	For each ordinal $\beta$, write
	\[
	\mathbb Q^{\textsf{ac}}\simeq\mathbb Q^{\textsf{ac}}_{<\beta}\ast\dot{\mathbb Q}^{\textsf{ac},\beta},
	\]
	where the tail $\dot{\mathbb Q}^{\textsf{ac},\beta}$ starts at coordinate $\beta$.
	
	\begin{proposition}[Blue, Larson, Sargsyan~\cite{NM2}]\label{prop:bls-tail-closure}
		Let $\beta$ be a successor ordinal and let $H_\beta\subseteq\mathbb Q^{\textsf{ac}}_{<\beta}$
		be $\mathsf{M}[G]$-generic. Let
		$\mathbb Q^{\textsf{ac},\beta}\df(\dot{\mathbb Q}^{\textsf{ac},\beta})^{H_\beta}$.
		Suppose that $\kappa$ is a cardinal of $\mathsf{M}[G][H_\beta]$ with
		$\kappa\leq\eta_{\alpha+\beta-1}$ and that
		$\mathsf{M}[G][H_\beta]$ satisfies $\dc_\kappa$. Then $\mathbb Q^{\textsf{ac},\beta}$ is
		${<}(\kappa^+)^{\mathsf{M}[G][H_\beta]}$-closed in $\mathsf{M}[G][H_\beta]$.
	\end{proposition}
	
	\begin{proof}
		See~\cite[Claim~8.8]{NM2}.
	\end{proof}
	
	\begin{remark}
		In our application of \cref{thm:forcing-choice} and the subsequent analyses, we will initially take $\alpha=0$ and $\mathbb P=\col{\omega_1}{\R}$ although at some points we have to work with  intermediate stages.
	\end{remark}
	
	\subsection{Minimal Nairian models}\label{ssec:minimal-Nairian-models}
	Throughout this subsection, assume that the parent model $V$ satisfies $\hpc$.
	Assume that $N$ is a minimal Nairian model, and let $\xi\df N\cap\ord$.
	\begin{theorem}[Blue, Larson, Sargsyan~\cite{NM,NM2}]
		\label{thm:strong-cardinals-in-nairian-models}
		The $\Theta$-sequence of $N$ is the increasing enumeration of the ${<}\xi$-strong cardinals of $\hod$ below $\xi$ and their limits.
	\end{theorem}
	\begin{proof}
		See~{\cite[Theorem~19.1]{NM2}}.
	\end{proof}
	\begin{theorem}[Blue, Larson, Sargsyan~\cite{NM,NM2}]
		\label{thm:theta-sequence-properties}
		Let $\alpha\in\{-1\}\cup\xi$. Then the following hold in $N$.
		\begin{enumerate}
			\item\label{itm:theta-sequence-theta-and-regularity}
			$\kappa_{\alpha+1}$ is a regular cardinal.
			
			\item\label{itm:theta-sequence-subsets}
			$\mcl P(\kappa_\alpha^\omega)\cap N \subseteq N|\kappa_{\alpha+1}$.
		\end{enumerate}
	\end{theorem}
	
	\begin{proof}
		These are clauses~(2a) and~(2c) of \cite[Theorem~1.3]{NM2}.
	\end{proof}
	
	\begin{theorem}[Blue, Larson, Sargsyan~\cite{NM2}]
		\label{thm:consequences-of-nairian-models}
		The Nairian model  $N$ satisfies $\ntbase$, $\ad_{\R}$, and $\dc$. Moreover,
		$\omega_1$ is supercompact in $N$.
	\end{theorem}
	
	\begin{proof}
		By \cite[Theorems~1.17 and~1.3(8)]{NM2}, $N$ satisfies $\ntbase$
		and $\omega_1$ is supercompact in $N$.
		The model $N$ also satisfies $\ad_{\R}$ by \cite[Theorem~8.1]{NM} and the discussion following Definition~2.5 there. To prove $\dc$ in $N$, let $T\in N$ be a tree of height $\omega$ with no terminal nodes. Solovay's theorem \cite{Sol} gives $\dc$ in
		the ambient universe, so $T$ has a cofinal branch $b$.
		Choose $\lambda<N\cap\ord$ such that $T\subseteq N|\lambda$.
		By \cite[Theorem~14.12]{NM2}, $(N|\lambda)^\omega\subseteq N$.
		So  $b\in N$. Thus $\dc$ holds in $N$.  
	\end{proof}
	
	\begin{lemma}\label{lem:heights-and-theta-sequence}
		Let $N$ be a minimal Nairian model. Let   $\xi\df N\cap\ord$.
		Let $\eta_\gamma\df\hH^N_\gamma\cap\ord$ for $\gamma<\xi$.
		Then for every $\gamma<\xi$, $\hH^N_\gamma=N|\eta_\gamma$, and 
		\[
		\eta_\gamma=
		\begin{cases}	\kappa_\delta,&\gamma=\delta,\quad \text{where $\delta$ is zero or a limit ordinal}\\
			\kappa_{\delta+n},&\gamma=\delta+n,\quad 0<n<\omega,~\cof{\delta}^N=\omega \text{ or } \delta=0\\
			\kappa_{\delta+n+1},&\gamma=\delta+n,\quad 0<n<\omega,~\cof{\delta}^N>\omega.
		\end{cases}
		\]
		
	\end{lemma}
	\begin{proof}
		See~\cite{NM2}.
	\end{proof}
	
	\begin{lemma}\label{lem:singular-limit-heights}
		Let $N$ be a minimal Nairian model at $\xi$, and compute the
		$\hH$-hierarchy inside $N$. For every nonzero limit ordinal
		$\delta<\xi$, the height
		$\eta_\delta\df\hH^N_\delta\cap\ord$ is a singular cardinal in $N$.
	\end{lemma}
	\begin{proof}
		By \cref{lem:heights-and-theta-sequence}
		and continuity of the $\Theta$-sequence of $N$, $\alpha\mapsto \kappa_\alpha^N$, we have $\eta_\delta=\kappa^N_\delta$. So $\eta_\delta$ is a cardinal in $N$ by
		\cref{thm:theta-sequence-properties}\textup{(\ref{itm:theta-sequence-theta-and-regularity})}.
		By \cref{thm:strong-cardinals-in-nairian-models}, $\eta_\delta$
		is a limit of ${<}\xi$-strong cardinals of the ambient $\hod$.
		These cardinals are inaccessible in $\hod$, so $\eta_\delta$
		is strong limit there. By minimality of $N$,  $\eta_\delta$ must be singular in $N$.
	\end{proof}
	
	\section{Measures and embeddings under determinacy}
	\label{sec:measures-determinacy}
	The material in this section is influenced by work of Woodin; see
	\cite{KecWoo08inbookGene,Wilson-thesis}.
	We assume $\zf+\ad+\dc$ throughout this section.
	We define weak comeagerness and $\mu$-Uniformization in
	\cref{subsec:weak-comeagerness}.
	Then, in \cref{subsec:generic-elementary-embedding}, we construct the
	generic ultrafilter $\mu/g$ from a supercompactness measure $\mu$,
	prove that its ultrapower is well-founded, and show that $i_\mu$ and
	$i_{\mu/g}$ agree on ordinals.
	Under $\mu$-Uniformization, we also establish the \Los{} property.
	Finally,  we compute $i_\mu(\omega_1)$.

	\subsection{Weak comeagerness}\label{subsec:weak-comeagerness}
	Throughout this subsection, we fix an uncountable set $S$.
	Unless stated otherwise, for every nonempty set $\sigma$, we equip $\sigma^\omega$ with the product topology, where $\sigma$ carries the discrete topology. Note that the basic open sets are 
	\[
	\mcl N^\sigma_p\df\{x\in\sigma^\omega:x\supseteq p\},
	\qquad p\in\sigma^{<\omega}.
	\]
	This agrees with the standard notation for the Baire space $\omega^\omega$, which we shall denote by $\R$.
	It is provable in $\zf$ that for any countably infinite set $\sigma$, there is a homeomorphism between the product spaces $\sigma^\omega$ and $\R$. Thus, we can discuss the topological properties of $\sigma^\omega$ that follow from $\ad$ without explicitly fixing a homeomorphism. 
	
	Let us also recall that a subset $A$ of a topological space $T$ is called nowhere dense if the interior of the closure of $A$ is empty. The set $A$ is called meager if it is a countable union of nowhere dense subsets of $T$, and it is called comeager if $T\setminus A$ is meager. Finally, $A$ has the Baire property if $A\Delta U$ is meager for some open set $U\subseteq T$.
	It is well known that, under $\ad$, all sets of reals have the Baire property; hence all subsets of $\sigma^\omega$ have the Baire property whenever $\sigma$ is countably infinite.
	
	The following alternative notation will be useful. Although introduced for arbitrary infinite sets, it will most often be used when $R$ is uncountable.
	\begin{notation}\label{not:sequence-space}
		For an infinite set $R$,
		we shall often denote $R^\omega$ by $R_\infty$.
	\end{notation}
	Notice that $|R_\infty|= |(R_\infty)^\omega|$.  Note also that if $Q\subseteq R$ is infinite, then $Q_\infty$ is a topological subspace of $R_\infty$. In particular, for every countably infinite $\sigma\in\mcl P_{\omega_1}(S)$,  $\sigma^\omega$ is a topological subspace of $S_\infty$ and is homeomorphic to the Baire space.
	
	\begin{lemma}\label{lem:dense}
		Assume that $\mu$ is a supercompactness measure on $\mcl P_{\omega_1}(S)$.
		Let $A\subseteq S_\infty$.
		The following are equivalent.
		\begin{enumerate}
			\item\label{itm:density-in-sequence-space} $A$ is dense in $S_\infty$.
			\item\label{itm:density-on-almost-all-supports} $\{\sigma\in\mcl P_{\omega_1}(S): A\cap\sigma^\omega \text{ is dense in } \sigma^\omega\}\in\mu$.
		\end{enumerate}
	\end{lemma}
	\Needspace{3\baselineskip}
	\begin{proof}\leavevmode
		\begin{itemize}
			\item ($1\Rightarrow 2$) Suppose, toward a contradiction, that (2) fails. Since $\mu$ is an ultrafilter, 
			\[\mathcal Z\df \{\sigma\in\mcl P_{\omega_1}(S): A\cap\sigma^\omega \text{ is not dense in } \sigma^\omega\}\in\mu.\]
			So for every $\sigma\in\mathcal Z$, there is $q\in\sigma^{<\omega}$ such that $\mcl N^\sigma_{q}\cap A=\emptyset$.
			Define $F$ on $\mathcal Z$ by  $F(\sigma)\df
			\{q\in\sigma^{<\omega}:\mathcal N_q^\sigma\cap A=\emptyset\}$.
			Applying \cref{lem:normality-fin} to $F$ gives $q_*\in S^{<\omega}$ and a set $\mathcal Z_*\subseteq\mathcal Z$ in $\mu$ such that $\mathcal N^{\sigma}_{q_*}\cap A=\emptyset$ for every $\sigma\in\mathcal Z_*$. Since $A$ is dense in $S_\infty$, there is $x\in A\cap\mcl N^S_{q_*}$. By \cref{lem:fine-and-complete}, there is $\sigma\in \mcl Z_*$ with $\rang{x}\subseteq\sigma$, so
			$x\in A\cap\mathcal N_{q_*}^\sigma$, a contradiction.

			\item ($2\Rightarrow 1$) Pick an arbitrary basic open set  $\mcl N^S_{q}$ in $S_\infty$.
			By \cref{lem:fine-and-complete}, there is a nonempty $\sigma\in\mcl P_{\omega_1}(S)$ such that $\rang{q}\subseteq\sigma$ and $A\cap\sigma^\omega$ is dense in $\sigma^\omega$. Let $y\in \mathcal N_{q}^\sigma\cap A$. Then $y\in  \mcl N^S_{q}\cap A$. Therefore, $A$ is dense in $S_\infty$.\qedhere
		\end{itemize}
	\end{proof}
	\begin{lemma}\label{lem:nowheredense}
		Assume that $\mu$ is a supercompactness measure on $\mcl P_{\omega_1}(S)$. Let $A\subseteq S_\infty$ be nowhere dense.
		Then $\{\sigma\in\mcl P_{\omega_1}(S): A\cap\sigma^\omega\text{ is nowhere dense in }\sigma^\omega\}\in\mu$.
	\end{lemma}
	\begin{proof}
		For a topological space $T$ and a set $B\subseteq T$, let us
		write $\operatorname{cl}_T(B)$ for the topological closure of $B$ in $T$. Put $D\df S_\infty\setminus\operatorname{cl}_{S_\infty}(A)$. Since $A$ is
		nowhere dense, $D$ is dense in $S_\infty$. By \cref{lem:dense}, we have
		\[
		U\df\{\sigma\in\mcl P_{\omega_1}(S):
		D\cap\sigma^\omega\text{ is dense in }\sigma^\omega\}\in\mu.
		\]
		It is enough to show that $U$ is contained in $\{\sigma\in\mcl P_{\omega_1}(S): A\cap\sigma^\omega\text{ is nowhere dense in }\sigma^\omega\}$. Thus let us fix $\sigma\in U$. 
		The subspace topology gives $\operatorname{cl}_{\sigma^\omega}(A\cap\sigma^\omega)
		\subseteq
		\operatorname{cl}_{S_\infty}(A)\cap\sigma^\omega$.
		Consequently,
		$D\cap\sigma^\omega
		\subseteq
		\sigma^\omega\setminus
		\operatorname{cl}_{\sigma^\omega}(A\cap\sigma^\omega)$. Therefore, $\sigma^\omega\setminus
		\operatorname{cl}_{\sigma^\omega}(A\cap\sigma^\omega)$   is  dense in $\sigma^\omega$, and hence
		$A\cap\sigma^\omega$ is nowhere dense there, as required.
	\end{proof}
	
	\begin{definition}\label{def:relative-comeagerness}
		Let $A\subseteq S_\infty$. Let $\sigma\in\mcl P_{\omega_1}(S)$, and let $p\in \sigma^{<\omega}$. We say that
		$A$ is \emph{comeager on $\mcl N_p^\sigma$} if $A\cap\mcl N_p^\sigma$ is comeager in the relative topology of $\mcl N_p^\sigma$ induced by $\sigma^\omega$. We also say $A$ is \emph{meager on $\mcl N_p^\sigma$} if $\mcl N_p^\sigma\setminus A$ is comeager in $\mcl N_p^\sigma$.
	\end{definition}

	The following definition extends the notion of comeagerness to spaces
	$S^\omega$, where $S$ is uncountable, though only in a weak sense.
	
	\begin{definition}\label{def:weak-comeagerness}
		Assume that $\mu$ is a supercompactness measure on $\mcl P_{\omega_1}(S)$.
		Let $A\subseteq S_\infty$, and let $p\in S^{<\omega}$.
		We say that $A$ is \emph{weakly comeager below $p$ with respect to $\mu$} if
		\[
		\{\sigma\in\mcl P_{\omega_1}(S) \mid p\in\sigma^{<\omega} \land A \text{ is comeager on }  \mcl N_p^\sigma \}\in\mu.
		\]
		We also say that $A$ is \emph{weakly comeager with respect to $\mu$} if it is weakly comeager below $\emptyset$ with respect to $\mu$.
	\end{definition}
	
	\begin{lemma}\label{lem:basic-open-weak-comeagerness}
		Assume that $\mu$ is a supercompactness measure on $\mcl P_{\omega_1}(S)$.
		For every $p\in S^{<\omega}$, every comeager subset of $\mathcal N_p^S$
		is weakly comeager below $p$ with respect to $\mu$.
		In particular, every comeager subset of $S_\infty$ is weakly comeager
		with respect to $\mu$.
	\end{lemma}
	\begin{proof}
		Let $p\in S^{<\omega}$, and let $A\subseteq\mathcal N_p^S$ be
		comeager in $\mathcal N_p^S$. Since $\mathcal N_p^S$ is clopen in
		$S_\infty$, we can write $\mathcal N_p^S\setminus A$ as a countable union of  nowhere dense sets in $S_\infty$, say $\{A_n:n<\omega\}$.
		By \cref{lem:nowheredense}, fineness, and countable completeness, for
		$\mu$-almost every $\sigma\in\mcl P_{\omega_1}(S)$ we have
		$p\in\sigma^{<\omega}$ and $A_n\cap\sigma^\omega$ is nowhere dense
		in $\sigma^\omega$ for every $n<\omega$. For each such $\sigma$, we have
		\[
		\mathcal N_p^\sigma\setminus A
		=\bigcup_{n<\omega}(A_n\cap\mathcal N_p^\sigma).
		\]
		Thus $\mathcal N_p^\sigma\setminus A$ is meager in $\mathcal N_p^\sigma$, since $\mathcal N_p^\sigma$ is an open subspace
		of $\sigma^\omega$. Consequently, $A$ is weakly comeager below $p$ with
		respect to $\mu$. Taking $p=\emptyset$ gives the final assertion.
	\end{proof}
	
	\begin{lemma}\label{lem:generalized-baire-property}
		Assume that $\mu$ is a supercompactness measure on $\mcl P_{\omega_1}(S)$.
		Let $A\subseteq S_\infty$, and let $p\in S^{<\omega}$.
		Then either $A$ is weakly comeager below $p$ with respect to $\mu$ or $S_\infty\setminus A$ is weakly comeager below some $q\leq p$ with respect to $\mu$.
	\end{lemma}
	\begin{proof}
		Suppose that $A$ is not weakly comeager below $p$ with respect to $\mu$. Since $\mu$ is a fine ultrafilter, the set $U$ consisting of $\sigma\in\mcl P_{\omega_1}(S)$ such that $p\in\sigma^{<\omega}$ and $A\cap\sigma^\omega$ is not comeager on $\mcl N_p^\sigma$ belongs to $\mu$. 
		
		Since every set of reals has the Baire property under $\ad$, for every $\sigma\in U$ there is an extension $q\in\sigma^{<\omega}$ of $p$ (that is, $q\leq p$) such that $S_\infty\setminus A$ is comeager on $\mcl N_q^\sigma$.
		Define $F$ on $U$ by
		\[F(\sigma)\df
		\{q\in\sigma^{<\omega}:q\leq p
		\text{ and }S_\infty\setminus A
		\text{ is comeager on }\mcl N_q^\sigma\}.\]
		Applying \cref{lem:normality-fin} to $F$ gives a tuple $q$ and a set $U'\subseteq U$ in $\mu$ such that $S_\infty\setminus A$ is comeager on $\mcl N_q^\sigma$ for every $\sigma\in U'$. Consequently,
		\[
		\{\sigma \in \mcl P_{\omega_1}(S) \mid q\in\sigma^{<\omega}\land (S_\infty\setminus A)\cap\mcl N_q^\sigma \text{ is comeager in } \mcl N_q^\sigma \}\in\mu.
		\]
		Therefore, $S_\infty\setminus A$ is weakly comeager below $q\leq p$ with respect to $\mu$.
	\end{proof}
	
	\begin{definition}[$\mu$-Uniformization]\label{def:mu-uniformization}
		Let $\mu$ be a supercompactness measure on
		$\mcl P_{\omega_1}(S)$.
		We say that \emph{$\mu$-Uniformization} holds if, for every set
		$a$, every relation $R\subseteq S^\omega\times a$,
		and every $p\in S^{<\omega}$ such that $\dom{R}$ is weakly comeager
		below $p$ with respect to $\mu$, there exists a function $f\subseteq R$ whose domain is weakly comeager
		below $q$ with respect to $\mu$, for some $q\leq p$.
		
	\end{definition}
	\subsection{Generic elementary embeddings}\label{subsec:generic-elementary-embedding}
	
	We fix an infinite set $X$ and  a supercompactness measure $\mu$ on
	$\mcl P_{\omega_1}(X_\infty)$. We will use
	\cref{subsec:weak-comeagerness}, with $X_\infty$ in the role of $S$.

	\begin{lemma}\label{lem:theta-preserved}
		Assume $\Theta_{X_\infty}$ is regular. Let $g \subseteq \col{\omega}{X_\infty}$ be a $V$-generic filter. Then
		$\omega_1^{V[g]}=\Theta_{X_\infty}^V$. 
	\end{lemma}
	\begin{proof}
		Throughout the proof, write $\Theta_{X_\infty}=\Theta_{X_\infty}^{V}$. For every nonzero
		$\alpha<\Theta_{X_\infty}$, there is in $V$ a surjection from $X_\infty$ onto
		$\alpha$. Since forcing with $\col{\omega}{X_\infty}$ makes
		$X_\infty$ countable, every $\alpha<\Theta_{X_\infty}$ is countable in $V[g]$.
		It remains to show that $\Theta_{X_\infty}$ itself is not countable in $V[g]$.
		
		Let $\dot\pi$ be a $\col{\omega}{X_\infty}$-name for a function from
		$\omega$ into $\Theta_{X_\infty}$. Working in $V$, define the function $f:\col{\omega}{X_\infty}\longrightarrow\Theta_{X_\infty}$ by 
		\[
		f(p)\df
		\sup\{\alpha+1:\exists n<\omega\;
		p\forces{}{\dot\pi(n)=\check\alpha}\}.
		\]
		By the regularity and
		uncountability of $\Theta_{X_\infty}$, we must have $f(p)<\Theta_{X_\infty}$.
		
		Since $\col{\omega}{X_\infty}$ is  equinumerous with $X_\infty$, \cref{lem:reg-cof}
		implies that $f$ is bounded. Choose $\gamma<\Theta_{X_\infty}$ with
		$\rang{f}\subseteq\gamma$. So the maximal condition forces
		$\rang{\dot\pi}\subseteq\gamma$. In particular, no function from
		$\omega$ onto $\Theta_{X_\infty}$ is added. Thus $\Theta_{X_\infty}$ remains uncountable while
		every ordinal below it becomes countable, and hence 
		$\omega_1^{V[g]}=\Theta_{X_\infty}^{V}$.
	\end{proof}
	
	Let $g \subseteq \col{\omega}{X_\infty}$ be a $V$-generic filter. Working in
	$V[g]$ from now on, we define
	\[
	\mu/g \df \{ A \in \mathcal{P}(X_\infty^\omega)^V :
	A \text{ is weakly comeager below some } p \in g
	\text{ with respect to } \mu \},
	\]
	where weak comeagerness is computed in $V$.
	
	The notation $\mu/g$ should not be interpreted as denoting a quotient; as
	the next lemma shows, it denotes the generic ultrafilter derived from $\mu$.
	
	\begin{lemma}\label{lem:generic-ultrafilter}
		The set $\mu/g$ is a nonprincipal
		ultrafilter on $\mcl P(X_\infty^\omega)^V$.
	\end{lemma}
	\begin{proof}
		Throughout this proof we denote $(X_\infty^\omega)^V$ by $X_\infty^\omega$.
		It is easy to see that
		weak comeagerness is preserved under supersets. Suppose that
		$A,B\in\mu/g$, witnessed by $A^*,B^*\in\mu$ below $p,q\in g$, respectively. By directedness,
		choose $r\in g$ with $r\leq p,q$, and let $C_r\df\{\sigma\in\mcl P_{\omega_1}(X_\infty):
		r\in\sigma^{<\omega}\}$, 
		which belongs to $\mu$ by fineness.
		It is easily seen that $A^*\cap B^*\cap C_r$  witnesses that $A\cap B\in\mu/g$.  The whole space belongs to
		$\mu/g$ by \cref{lem:basic-open-weak-comeagerness}.
		Thus $\mu/g$ is a  filter on $\mcl P(X_\infty^\omega)^V$.
		
		To see that $\mu/g$ is an ultrafilter, fix $A\in\mcl P(X_\infty^\omega)^V$. 
		\[
		D_A\df\left\{p\in X_\infty^{<\omega}:
		\begin{aligned}
			&A\text{ or }X_\infty^\omega\setminus A\text{ is weakly comeager}\\
			&\text{below }p\text{ with respect to }\mu
		\end{aligned}\right\}.
		\]
		By \cref{lem:generalized-baire-property}, $D_A$ is dense in
		$\col{\omega}{X_\infty}$. Let $p\in g\cap D_A$. Then either $A$ or $X_\infty^\omega\setminus A$ belongs to    $\mu/g$. To see that $\mu/g$ is nonprincipal, it is enough to show that $\mu/g$ does not contain any singleton.  Thus let us fix $f\in X_\infty^\omega$. Then the set
		\[
		D_f\df\{p\in X_\infty^{<\omega}:p\text{ is incompatible with }f\}
		\]
		is dense in $\col{\omega}{X_\infty}$. Choose $p\in g\cap D_f$. By
		\cref{lem:basic-open-weak-comeagerness}, we have $\mcl N_p^{X_\infty}\in\mu/g$. Therefore, $\{f\}\notin\mu/g$.
	\end{proof}
	
	Our next goal is to compare the ultrapower by $\mu$ with the generic
	ultrapower by $\mu/g$. Since we do not yet assume that $\mu/g$ has the
	\Los{} property in the sense of \cref{def:los-property}, we first compare the parts of these ultrapowers
	represented by functions with \emph{hereditarily well-orderable} range. This
	will be sufficient for comparing their ordinal parts, and hence the
	action of the corresponding canonical ultrapower maps on ordinals.
	
	We need the following general definitions.  Let
	\[
	\hwo\df\{x:\tc{x}\text{ is well-orderable}\}.
	\]
	Thus $\hod\subseteq\hwo$.
	\begin{definition}
		Suppose that $\nu$ is an
		ultrafilter on an infinite set $Z$.
		Define the \emph{restricted ultrapower}
		\[
		\ultstar{\hwo}{\nu}
		\df
		\{f:\dom{f}=Z\text{ and }\rang{f}\in\hwo\}/\nu,
		\]
		with equality and membership induced by $\nu$.
	\end{definition}
	
	Let $\nu$ be as above, and let us define $j_\nu:\hwo\longrightarrow\ultstar{\hwo}{\nu}$ in a natural way by $x\mapsto[\lambda z.x]$ as in \cref{not:i_mu}, and let us call it
	the canonical ultrapower map. We have the  \Los{} property for
	$\Delta_0$-formulas: for every $\Delta_0$-formula $\varphi$ and representatives
	$f_0,\ldots,f_n$,
	\[
	\ultstar{\hwo}{\nu}\models{
		\varphi([f_0]_\nu,\ldots,[f_n]_\nu)}
	\Leftrightarrow
	\{z\in Z:\hwo\models{
		\varphi(f_0(z),\ldots,f_n(z))}\}\in\nu.
	\]
	We provide a brief argument using induction on the formulas. The atomic cases follow from the
	definition of the ultrapower, and the Boolean steps use the ultrafilter
	properties of $\nu$. For the bounded existential step, suppose that
	\[
	\hwo\models{\exists x\in h(z)~\psi(x,f_0(z),\ldots,f_n(z))}
	\]
	for $\nu$-almost every $z$. Since $\rang{h}\in\hwo$,
	the set $T\df\tc{\rang{h}}\cup\{\emptyset\}$ is well-orderable.
	Fix a well-order of $T$. For each such $z$, choose the least witness in
	this order, and put $\emptyset$ at the remaining inputs. This defines a
	function $h'$ whose range has well-orderable transitive closure, so $[h']_\nu$
	represents an element of the restricted ultrapower. By induction, its
	class witnesses the existential formula. For the reverse implication,
	apply the induction hypothesis to a representative of a witness in the
	ultrapower.
	
	Thus $j_\nu$ is a $\Delta_0$-elementary $\in$-embedding. 
	
	\begin{fact}\label{fact:extensionality-wellfounded-hwo}
		Suppose $\nu$ is as above.
		If $\nu$ is countably complete, then
		$\ultstar{\hwo}{\nu}$ is well-founded and extensional.
	\end{fact}
	\begin{proof}
		Extensionality follows from the \Los{} property for $\Delta_0$-formulas. The fact that the induced membership relation is well-founded follows from the countable completeness of $\nu$ via a routine argument using $\dc$.
	\end{proof}

	All ultrapowers by $\mu$ and $\mu/g$ below use representatives in $V$,
	with domains $\mathcal P_{\omega_1}(X_\infty)^V$ and $(X_\infty^\omega)^V$, respectively.
	For ease of notation, we omit the superscript $V$ and write
	$\mathcal P_{\omega_1}(X_\infty)$ and $X_\infty^\omega$.
	In $V[g]$, the   \Los{} property for $\Delta_0$-formulas holds for the    restricted ultrapowers $\ultstar{\hwo^V}{\mu}$ and $\ultstar{\hwo^V}{\mu/g}$.
	
	Now the map
	$f\mapsto\rang{f}$ from $X_\infty^\omega$ to
	$\mathcal{P}_{\omega_1}(X_\infty)$ induces a canonical map
	\[
	\begin{gathered}
		{\bullet}^*\colon V^{\mathcal P_{\omega_1}(X_\infty)}
		\longrightarrow V^{X_\infty^\omega},\\
		F\longmapsto F^*, 
	\end{gathered}
	\]
	where $F^*(f)\df F(\rang{f})$ for every $f\in X_\infty^\omega$.
	Since
	$\rang{F^*}\subseteq\rang{F}$, this map preserves $\hwo^V$, i.e., if
	$\rang{F}\in\hwo^V$, then $\rang{F^*}\in\hwo^V$. Thus, in $V[g]$, $\bullet^*$ induces an assignment
	\[
	\begin{gathered}
		\rho:\ultstar{\hwo^V}{\mu}
		\longrightarrow\ultstar{\hwo^V}{\mu/g},\\
		[F]_\mu\longmapsto[F^*]_{\mu/g}.
	\end{gathered}
	\]
	
	\begin{lemma}\label{lem:well-defined-embedding-hwo}
		Let $B\subseteq\mathcal P_{\omega_1}(X_\infty)$ belong to $V$. Then 
		\[
		B\in\mu
		\quad\Longleftrightarrow\quad
		\{f\in X_\infty^\omega:\rang{f}\in B\}\in\mu/g.
		\]
	\end{lemma}
	\begin{proof}
		Put $\widehat B\df\{f\in X_\infty^\omega:\rang{f}\in B\}$. For each countably infinite
		$\sigma\in\mathcal P_{\omega_1}(X_\infty)$, let
		$C_\sigma\df\{f\in\sigma^\omega:\rang{f}=\sigma\}$.
		Choose a surjection $e:\omega\twoheadrightarrow\sigma$.
		For each $n<\omega$, let $D_n\df\{f\in\sigma^\omega:e(n)\in\rang{f}\}$. Each $D_n$ is open and dense in $\sigma^\omega$, and
		$C_\sigma=\bigcap_{n<\omega}D_n$. Thus $C_\sigma$ is comeager.
		
		Now if $B\in\mu$, then $C_\sigma\subseteq\widehat B\cap\sigma^\omega$
		for every countably infinite $\sigma\in B$. Since $\mu$ concentrates
		on such $\sigma$'s, $\widehat B$ is weakly comeager with respect to $\mu$, and thus
		belongs to $\mu/g$. Conversely, if $B\notin\mu$, apply the
		forward implication to $B'\df \mathcal P_{\omega_1}(X_\infty)\setminus B$ to get
		\[
		X_\infty^\omega\setminus\widehat B
		=\widehat{B'}\in\mu/g.
		\]
		Since $\mu/g$ is an ultrafilter, $\widehat B\notin\mu/g$. 
		
	\end{proof}
	\begin{lemma}\label{lem:iso-hwo}
		The embedding	$\rho$ is an $\in$-isomorphism. Consequently, $\ultstar{\hwo^V}{\mu/g}$ is well-founded. 
	\end{lemma}
	
	\begin{proof}
		By \cref{lem:well-defined-embedding-hwo}, if
		$F,G\in V$ are functions on $\mathcal P_{\omega_1}(X_\infty)$ with
		$\rang{F},\rang{G}\in\hwo^V$, and  
		$R$ is either equality or membership, we have
		\[
		\{\sigma:F(\sigma)\mathrel{R}G(\sigma)\}\in\mu
		\quad\Longleftrightarrow\quad
		\{f:F^*(f)\mathrel{R}G^*(f)\}\in\mu/g.
		\]
		Therefore, $\rho$ is an $\in$-embedding.
		It remains to prove that it is surjective. So
		fix $G\in V$ with domain $X_\infty^\omega$ and
		$\rang{G}\in\hwo^V$, and put $R\df\rang{G}$. For each
		$F\in V$ with $F:\mathcal P_{\omega_1}(X_\infty)\to R$, write
		\[
		E_F\df\{f\in X_\infty^\omega:G(f)=F(\rang{f})\}.
		\]
		
		Notice that for $F$ as above, we have $\rang{F}\in \hwo^V$, and that $G\restriction E_F=F^*\restriction E_F$.
		
		Working in $V$,  let $D$ consist of  conditions $q\in\col{\omega}{X_\infty}$ such that for some  function $F:\mathcal P_{\omega_1}(X_\infty)\to R$ in $V$,  the set $E_F$
		is weakly comeager below $q$ with respect to $\mu$.  Assume for a moment that $D$ is dense in $\col{\omega}{X_\infty}$. Then, in $V[g]$, there is $q\in D\cap g$ and there is 
		$F\in V$ as above such that $E_F$ is weakly comeager below $q$ with respect to $\mu$.  So $E_F\in\mu/g$. Consequently,  $	\rho([F]_\mu)=[G]_{\mu/g}$, which proves surjectivity.
		
		We now prove that $D$ is dense. We work in $V$ and fix a well-order $\prec$ of $R$.
		Let $p$ be an arbitrary condition in $\col{\omega}{X_\infty}$. Let $U_p$ consist of the countably infinite
		$\sigma\in\mathcal P_{\omega_1}(X_\infty)$ containing $\rang{p}$ as a subset. By \cref{lem:fine-and-complete}, $U_p\in\mu$.
		For $\sigma\in U_p$, let $F_p(\sigma)$ be the $\prec$-least $y\in R$
		such that $G^{-1}(\{y\})$ is nonmeager on $\mathcal N_p^\sigma$.
		Such a $y$ exists because the fibers form a well-ordered cover of the
		nonempty Baire space $\mathcal N_p^\sigma$, and under $\ad$ a
		well-ordered union of meager sets is meager. For $\sigma\notin U_p$,
		put $F_p(\sigma)=\min_\prec R$. Then $F_p:\mcl P_{\omega_1}(X_\infty)\to \rang{F_p}\subseteq R$ is a function in $V$.  
		
		Define the function $\Phi$ on $U_p$ by letting
		\[
		\Phi(\sigma)\df\{r\in\sigma^{<\omega}:r\leq p\text{ and }
		G^{-1}(\{F_p(\sigma)\})\text{ is comeager on }\mathcal N_r^\sigma\}.
		\]
		The Baire property gives $\Phi(\sigma)\neq\emptyset$.
		By \cref{lem:normality-fin}, there are $U\subseteq U_p$ with $U\in\mu$ and
		$q\in X_\infty^{<\omega}$ such that $q\in\Phi(\sigma)$ for every
		$\sigma\in U$. In particular, $q\leq p$.
		
		For each $\sigma\in U$, let $C_\sigma$ be the comeager set of
		surjections $\omega\to\sigma$ from the proof of \cref{lem:well-defined-embedding-hwo}.
		For $f\in C_\sigma$, we have $F_p(\rang{f})=F_p(\sigma)$, and hence
		\[
		E_{F_p}\cap\mathcal N_q^\sigma
		\supseteq C_\sigma\cap G^{-1}(\{F_p(\sigma)\})\cap\mathcal N_q^\sigma.
		\]
		The right-hand side is comeager in $\mathcal N_q^\sigma$ for each $\sigma\in U$. Thus
		$E_{F_p}$ is weakly comeager below $q\leq p$ with respect to $\mu$.
		Since $p$ was arbitrary, this proves the required density.
		
		Well-foundedness follows from \cref{fact:extensionality-wellfounded-hwo} and the isomorphism $\rho$.
	\end{proof}
	
	\begin{lemma}\label{lem:well-founded-transitivized-ultrapowers}
		Let $\nu\in\{\mu,\mu/g\}$. Then $\ult{V}{\nu}$ is an end-extension of $\ultstar{\hwo^V}{\nu}$, and hence well-founded.
	\end{lemma}
	\begin{proof}	
		Work in $V[g]$, and let $\nu\in\{\mu,\mu/g\}$.
		Suppose that $[F]_\nu\in[G]_\nu$ and $\rang{G}\in\hwo^V$. Choose $F'\in V$ such that
		$F(z)=F'(z)\in G(z)$ on a $\nu$-measure-one set $A$ and $F'(z)=\emptyset$ outside $A$.
		Then $\rang{F'}\subseteq\bigcup\rang{G}\cup\{\emptyset\}$. Since $\rang{G}\in\hwo^V$, it follows that $\rang{F'}\in\hwo^V$.
		Thus $\ult{V}{\nu}$ is an end-extension of $\ultstar{\hwo^V}{\nu}$. By \cref{fact:extensionality-wellfounded-hwo,lem:iso-hwo}, $\ultstar{\hwo^V}{\nu}$ is well-founded and extensional and contains all ordinals. We easily conclude that $\ult{V}{\nu}$ is well-founded as well.
	\end{proof}
	
	\begin{remark}\label{rem:well-foundedness-and-extensionality}
		By \cref{lem:well-founded-transitivized-ultrapowers}, for each
		$\nu\in\{\mu,\mu/g\}$ there is a unique map $c_\nu$ on $\ult{V}{\nu}$ satisfying
		\[
		c_\nu(a)=\{c_\nu(b):b\in_\nu a\},
		\]
		where $\in_\nu$ is its membership relation. Put
		\[i_\nu(x)\df c_\nu([\lambda z.~x]_\nu).\]
		The range of $c_\nu$ is
		transitive, and $c_\nu$ is injective if and only if the    ultrapower
		is extensional. In that case, $c_\nu$ is the Mostowski collapse, and
		we identify the ultrapower with its transitive image. If $\nu$ has
		the  \Los{} property over $V$, extensionality holds and $i_\nu$
		is elementary.
	\end{remark}
	
	\begin{lemma}\label{lem:agreement-on-ordinals}
		We have $i_\mu\restriction\hwo^V=i_{\mu/g}\restriction\hwo^V$.
		In particular, $i_\mu\restriction\ord=i_{\mu/g}\restriction\ord$.
	\end{lemma}
	\begin{proof}
		For each $\nu\in\{\mu,\mu/g\}$, let $\pi_\nu$ denote the transitive
		collapse of $\ultstar{\hwo^V}{\nu}$. By
		\cref{lem:well-founded-transitivized-ultrapowers}, $\ult{V}{\nu}$ is
		an end-extension of $\ultstar{\hwo^V}{\nu}$. Uniqueness
		of well-founded recursion therefore gives
		\[
		c_\nu\restriction\ultstar{\hwo^V}{\nu}=\pi_\nu.
		\]
		By \cref{lem:iso-hwo}, the two restricted ultrapowers have the same
		transitive collapse $W$, and $\pi_{\mu/g}\circ\rho=\pi_\mu$.
		We also have $\rho\circ j_\mu=j_{\mu/g}$. Thus the following
		diagram commutes.
		\[\begin{tikzcd}
			&& {\ultstar{\hwo^V}{\mu}} \\
			\hwo^V &&&& W \\
			&& {\ultstar{\hwo^V}{\mu/g}}
			\arrow["{\pi_\mu}", from=1-3, to=2-5]
			\arrow["\rho", from=1-3, to=3-3]
			\arrow["{j_\mu}", from=2-1, to=1-3]
			\arrow["{j_{\mu/g}}"', from=2-1, to=3-3]
			\arrow["{\pi_{\mu/g}}"', from=3-3, to=2-5]
		\end{tikzcd}\]
		For $x\in\hwo^V$, it follows that
		\begin{equation*}
		\begin{aligned}[b]
			i_{\mu/g}(x)
			&=\pi_{\mu/g}(j_{\mu/g}(x))\\
			&=\pi_{\mu/g}(\rho(j_\mu(x)))\\
			&=\pi_\mu(j_\mu(x))=i_\mu(x).
		\end{aligned}\qedhere
		\end{equation*}
	\end{proof}
	
	\begin{lemma}\label{lem:reals-in-ult}
		Assume that $\ult{V}{\mu/g}$ is extensional.
		Let $M_{g}$ be its transitive collapse. Then
		$\R^{V[g]}\subseteq M_{g}$.
	\end{lemma}
	\begin{proof}
		Fix $x\in\R^{V[g]}$, and choose a $\col{\omega}{X_\infty}$-name
		$\dot x\in V$ and a condition $p\in g$ such that $p\forces{}{\dot x\in\R}$
		and $\dot x^g=x$. We are looking for a function
		$F:X_\infty^\omega\to\R$ in $V$ such that
		\[
		c_{\mu/g}([F]_{\mu/g})=x. 
		\]

		Working in $V$, we define
		\[
		D\df\{f\in\mathcal N_p^{X_\infty}:(\forall n<\omega)(\exists k\geq|p|)~
		f\restriction k\text{ decides }\dot x(n)\}.
		\]
		We first show that $D$ is comeager in $\mathcal N_p^{X_\infty}$.
		To see this, for each $n<\omega$, let
		\[
		D_n\df \bigcup\{\mathcal N_q^{X_\infty}:q\leq p
		\text{ and }q\text{ decides }\dot x(n)\}.
		\]
		Note that each $D_n$ is open and dense in $\mathcal N_p^{X_\infty}$ and that $D=\bigcap_{n<\omega}D_n$. Thus $D$ is comeager in $\mathcal N_p^{X_\infty}$.
		By \cref{lem:basic-open-weak-comeagerness}, $D$ is weakly comeager below $p$ with respect to $\mu$, so 
		$D\in\mu/g$.
		
		We now define the function $F$ as follows.
		For $f\in D$, first define $F_0(f)(n)$ to be the unique $m<\omega$ such that for some $k\geq|p|$,
		\[
		f\restriction k\forces{}{\dot x(n)=\check m}.
		\]
		Extend $F_0$ arbitrarily to a function $F:X_\infty^\omega\to\R$ in $V$, and set $z\df c_{\mu/g}([F]_{\mu/g})$. We show that $z=x$.
		By \cref{lem:agreement-on-ordinals}, $i_{\mu/g}(\omega)=\omega$, and thus
		$z\subseteq\omega\times\omega$. For $n,m<\omega$, by
		definition we have  
		\[
		(n,m)\in z
		\quad\Longleftrightarrow\quad
		\{f\in X_\infty^\omega:F(f)(n)=m\}\in\mu/g.
		\]
		Therefore $z$ is a
		partial function from $\omega$ to $\omega$.
		
		Fix $n<\omega$, put $m=x(n)$, and choose $q\in g$ with $q\leq p$
		forcing $\dot x(n)=\check m$. Therefore,
		\[
		D\cap\mathcal N_q^{X_\infty}
		\subseteq\{f\in X_\infty^\omega:F(f)(n)=m\}.
		\]
		The left-hand side belongs to $\mu/g$, so $(n,x(n))\in z$.
		Thus $x=z$.
	\end{proof}
	
	\begin{corollary}\label{cor:omega1-computation}
		Assume that $\Theta_{X_\infty}^V$ is regular and that $\mu/g$ has the
		\Los{} property over $V$. Let $M_{g}$ denote the transitive collapse of $\ult{V}{\mu/g}$. Then
		\[
		\omega_1^{M_{g}}=\Theta_{X_\infty}^V=\omega_1^{V[g]}.
		\]
	\end{corollary}
	\begin{proof}
		Note that $\ult{V}{\mu/g}$ is extensional, since $\mu/g$ has the
		\Los{} property over $V$.
		By \cref{lem:reals-in-ult}, $M_{g}$ and $V[g]$ have the same reals. As every countable ordinal can be coded by a real, they also have the same countable ordinals. Hence $\omega_1^{M_{g}}=\omega_1^{V[g]}$.
		On the other hand, by \cref{lem:theta-preserved} we have $\omega_1^{V[g]}=\Theta_{X_\infty}^V$.
	\end{proof}
	
	\begin{proposition}\label{prop:los-property}
		Assume that $\mu$-Uniformization holds. Then the ultrafilter $\mu/g$ has the
		\Los{} property over $V$.
	\end{proposition}
	\begin{proof}
		The ultrafilter property of $\mu/g$ was established in
		\cref{lem:generic-ultrafilter}. We prove the \Los{} property by induction on the complexity of formulas. Atomic formulas and Boolean
		connectives are handled by the usual induction. We verify the
		existential step.   So let $\phi(y,x_0,\ldots,x_{n-1})$ be a formula for which the induction
		hypothesis holds, and let $F_i:X_\infty^\omega\to V$, for $i<n$,  be functions in $V$. We want to show that
		\[
		\begin{gathered}
		\ult{V}{\mu/g}\models{\exists y~\phi(y,[F_0]_{\mu/g},\dots,[F_{n-1}]_{\mu/g})}\\
		\text{if and only if}\\
		U\df
		\{f\in X_\infty^\omega:V\models{
			\exists y~\phi(y,F_0(f),\ldots,F_{n-1}(f))}\}\in\mu/g.
		\end{gathered}
		\]
		The forward implication is trivial. Therefore, let us assume $U\in\mu/g$ and prove that there is $G:X_\infty^\omega\to V$ in $V$ such that
		\[
		\ult{V}{\mu/g}\models{\phi([G]_{\mu/g},[F_0]_{\mu/g},\dots,[F_{n-1}]_{\mu/g})}.
		\]
		
		Since $U\in\mu/g$, there exists
		$p\in g$ witnessing that $U$ is weakly comeager below $p$
		with respect to $\mu$. Fix such a condition $p$.
		We will apply $\mu$-Uniformization to $U$, conditions below $p$, and the following relation.
		We now work in $V$ and  fix a set $a\in V$ with the property that
		for every $f\in U$, there exists $y\in a$ such that $\phi(y,F_0(f),\ldots,F_{n-1}(f))$ holds. Then we define
		\[
		R\df\{(f,y)\in U\times a:
		\phi(y,F_0(f),\ldots,F_{n-1}(f))\}.
		\]
		Notice that $\dom{R}=U$. Let $q\leq p$ be an arbitrary condition. Then $U$ is also weakly comeager below $q$ with respect to $\mu$. By $\mu$-Uniformization, there exist a condition $r\leq q$, a set $U'\subseteq U$ which is weakly comeager below $r$ with respect to $\mu$, and a function $S\subseteq R$ with $\dom{S}=U'$.
		Therefore, the set of such conditions  is dense below $p$.
		
		In $V[g]$, since $p\in g$, choose $p'\in g$ below $p$, a set $U'\subseteq U$ in $V$, and a function $S\subseteq R$ in $V$ such that $U'$ is weakly comeager below $p'$ with respect to $\mu$ and $\dom{S}=U'$. Consequently, $U'\in\mu/g$.
		
		Now extend $S$ to  a function $G:X_\infty^\omega\longrightarrow V$ by setting
		\[
		G(f)\df
		\begin{cases}
			S(f),&f\in U',\\
			\emptyset,&f\notin U'.
		\end{cases}
		\]
		Note that $G\in V$. We have
		\[
		U'\subseteq \{f\in X_\infty^\omega:V\models{
			\phi(G(f),F_0(f),\ldots,F_{n-1}(f))}\}\in\mu/g.
		\]
		By the inductive hypothesis,
		\[
		\ult{V}{\mu/g}\models{
			\phi([G]_{\mu/g},[F_0]_{\mu/g},\ldots,[F_{n-1}]_{\mu/g})},
		\]
		and consequently
		\[
		\ult{V}{\mu/g}\models{
			\exists y~\phi(y,[F_0]_{\mu/g},\ldots,[F_{n-1}]_{\mu/g})}.
		\]
		
	\end{proof}

	\begin{proposition}\label{prop:countable-closure-of-ultrapower}
		Assume $\mu$-Uniformization holds.	The transitive collapse of $\ult{V}{\mu/g}$ is countably closed in $V[g]$.
	\end{proposition}
	\begin{proof}
		Let $M_g$ denote the transitive collapse of $\ult{V}{\mu/g}$, and identify each ultrapower class with its collapse. Let $x:\omega\to M_g$ be a sequence in $V[g]$. Fix a $\col{\omega}{X_\infty}$-name $\dot{x}$ with $\dot{x}^g=x$, and a condition $p\in g$ forcing that $\dot{x}$ is an $\omega$-sequence in $M_{\dot{g}}$. We will find a set $C$ which is countable in $M_g$ and contains $\rang{x}$.
		
		We now work in $V$. Since $p$ forces that each $\dot{x}(n)$ has a ground-model representative, Collection gives a set $a$ of functions with domain $X_\infty^\omega$ such that for every $q\leq p$ and every $n<\omega$, there are $F\in a$ and $r\leq q$ with $r\forces{}{[F]_{\mu/\dot{g}}=\dot{x}(n)}$.
		
		For each $q\leq p$ and every $n<\omega$, let 
		\[
		A_{q,n}\df
		\{F\in a:q\forces{}{[F]_{\mu/\dot g}=\dot{x}(n)}\}.
		\]
		For each $n<\omega$, let \[
		D_n\df\{q\leq p:A_{q,n}\neq\emptyset\}.
		\]
		Each $D_n$ is dense below $p$. Hence
		\[
			U\df\{f\in\mcl N_p^{X_\infty}:
			\forall n<\omega~\exists k\geq|p|~(f\restriction k\in D_n)\}
		\]
		is comeager in $\mcl N_p^{X_\infty}$ and hence weakly comeager below $p$ with respect to $\mu$, by \cref{lem:basic-open-weak-comeagerness}. We define a relation $R\subseteq U\times a^\omega$ by letting
		\[
		(f,\vec F)\in R
		\quad\Longleftrightarrow\quad
			\forall n<\omega~\exists k\geq|p|~
			(\vec F(n)\in A_{f\restriction k,n}).
		\]
		
		By $\dc$, for every $f\in U$, one can construct   a sequence $\vec F$ with $(f,\vec F)\in R$. Therefore,  $\dom{R}=U$.
		Let $q\leq p$ be arbitrary. Since $U$ is also weakly comeager below $q$ with respect to $\mu$, by $\mu$-Uniformization there exist a condition $r\leq q$ and a function $S_r\subseteq R$ whose domain $U_r$ is weakly comeager below $r$ with respect to $\mu$. Put
		\[
			b_r\df\{S_r(f)(n):f\in U_r,\ n<\omega\}\subseteq a.
		\]
		The set $b_r$ is a nonempty image of $U_r\times\omega$. Since $U_r\subseteq X_\infty^\omega$ and $|X_\infty^\omega\times\omega|=|X_\infty|$, there is a surjection $e_r:X_\infty\twoheadrightarrow b_r$.
		Since $q\leq p$ was arbitrary, the set of such conditions $r$ is dense below $p$.
		
		We now work in $V[g]$. By the above density, we can fix $r\in g$ extending $p$ and $S_r,U_r,b_r,e_r\in V$ as above.
		Let $C$ be represented by the function
		\[
			f\longmapsto\{e_r(y)(f):y\in\rang{f}\}.
		\]
		Each value of this function is countable, so $C$ is countable in $M_g$ by \cref{prop:los-property}. Fix $F\in b_r$ and choose $y\in X_\infty$ with $e_r(y)=F$. The set of $f$ with $y\in\rang{f}$ is dense open and hence belongs to $\mu/g$ by \cref{lem:basic-open-weak-comeagerness}. It follows that $[F]_{\mu/g}\in C$.

		We claim that $\rang{x}\subseteq C$. Fix $n<\omega$. For every $q\leq r$, $U_r$ is weakly comeager below $q$, so we can choose $f\in U_r$ extending $q$. Since $(f,S_r(f))\in R$, there is $k\geq |q|$ such that
		\[
			f\restriction k\forces{}{\dot{x}(n)=[S_r(f)(n)]_{\mu/\dot g}}.
		\]
		Thus the conditions forcing $\dot{x}(n)=[F]_{\mu/\dot g}$ for some $F\in b_r$ are dense below $r$. Since $r\in g$, we have $x(n)=[F]_{\mu/g}\in C$ for some $F\in b_r$.

		Fix an enumeration of $C$ in $M_g$. The function sending each $n$ to the least index of $x(n)$ in this enumeration is a real of $V[g]$, so it belongs to $M_g$ by \cref{lem:reals-in-ult}. Composing it with the enumeration gives $x\in M_g$.
	\end{proof}
	\begin{theorem}\label{thm:image-of-omega-one} 
		Assume that $\mu$-Uniformization holds and that $\Theta_{X_\infty}$ is regular. Then $i_\mu(\omega_1)=\Theta_{X_\infty}$.
	\end{theorem}
	\begin{proof}
		By \cref{prop:los-property}, $\mu/g$  has the \Los{} property over $V$. Consequently, the    ultrapower $\ult{V}{\mu/g}$ is extensional. Note that it is well-founded by
		\cref{lem:well-founded-transitivized-ultrapowers}. Let $M_g$
		denote its transitive collapse. Then
		$i_{\mu/g}:V\to M_g\subseteq V[g]$ is elementary. Since $\Theta_{X_\infty}$ is a regular cardinal, we can apply \cref{cor:omega1-computation,lem:theta-preserved}. In particular, by
		\cref{lem:agreement-on-ordinals} and elementarity, we have
		\[
		i_\mu(\omega_1)=i_{\mu/g}(\omega_1)=\omega_1^{M_g}=\omega_1^{V[g]}=\Theta_{X_\infty}.
		\]
	\end{proof}
	\section{Lifting supercompactness embeddings}\label{sec:lifting-supercompactness-embeddings}

	We assume $\zf+\ad+\dc$ throughout this section.
	Let $X\in\hwo$ be infinite, and let
	$\mu$ be a supercompactness measure on $\mcl P_{\omega_1}(X_\infty)$. Consider the following setup.
	\begin{enumerate}[label=(S\arabic*),ref=S\arabic*]
		\item\label{itm:lifting-ultrapower-hypotheses}We assume that $\Theta_{X_\infty}$ is regular.
		
		\item\label{itm:lifting-countable-target-closure} Whenever
		$g\subseteq\col{\omega}{X_\infty}$ is $V$-generic, $\mu/g$ has the
		\Los{} property over $V$, and the transitive collapse
		$M_g\df\ult{V}{\mu/g}$ is countably closed in $V[g]$.
		
		\item\label{itm:lifting-surjection-onto-forcing} Let $\mbb A$ be a forcing notion with a surjection
		$c:X_\infty\twoheadrightarrow\mbb A$.
		
		\item\label{itm:lifting-well-orderability-and-choice} The forcing $\mbb A$ forces that $X_\infty$ is well-orderable and 
		$\dc_{|X_\infty|}$ holds.
		
		\item\label{itm:lifting-countable-closure} We have an $\mbb A$-name $\dot{\mbb B}$
		for a forcing notion such that $\mbb A\ast\dot{\mbb B}$ is countably closed and $\mbb A$ forces that
		$\dot{\mbb B}$ is ${<}|X_\infty|^+$-closed.
		
		\item\label{itm:lifting-target-enumeration} For every $V$-generic filter
		$G_{\mathbb A}\ast G_{\dot{\mathbb B}}\subseteq\mbb A\ast\dot{\mbb B}$, there is 
		a surjection
		\[
		\tau:\Theta_{X_\infty}^V
		\twoheadrightarrow
		\mcl P(X_\infty)^V
		\]
		in  $V[G_{\mathbb A}][G_{\mathbb B}]$ such that $\tau\restriction\alpha\in V[G_{\mathbb A}]$ for every
		$\alpha<\Theta_{X_\infty}^V$. 
		
	\end{enumerate}
	Fix  such objects, and set 
	\[
	\lambda\df|X_\infty|^{V[G_{\mathbb A}]}\qquad\text{and}\qquad Z\df\mcl P_{\omega_1}(X_\infty)^V.
	\]

	\begin{remark}\label{rem:lifting-preliminaries}
		Let $G_{\mathbb A}\ast G_{\dot{\mathbb B}}$ be as in \textup{(\ref{itm:lifting-target-enumeration})}.
		\begin{enumerate}
			\item For every set $Y\in V$,
			$Y^\omega\cap V[G_{\mathbb A}][G_{\mathbb B}]=Y^\omega\cap V$. In particular, \[
			X^\omega\cap V[G_{\mathbb A}][G_{\mathbb B}]=X_\infty\cap V.\]
			We also have
			$\mcl P_{\omega_1}(X_\infty)^V=\mcl P_{\omega_1}(X_\infty)^{V[G_{\mathbb A}][G_{\mathbb B}]}$.
			
			\item\label{itm:lifting-no-new-subsets} 
			Notice that
			\textup{(\ref{itm:lifting-well-orderability-and-choice})} gives
			$|Z|=|X_\infty^\omega|=\lambda$ in $V[G_{\mathbb A}]$.
			Since $\dc_\lambda$ holds in $V[G_{\mathbb A}]$,  the ${<}\lambda^+$-closedness of $\dot{\mbb B}^{G_{\mathbb A}}$ in \textup{(\ref{itm:lifting-countable-closure})} gives
			\[
			\mcl P(Z)^{V[G_{\mathbb A}][G_{\mathbb B}]}=\mcl P(Z)^{V[G_{\mathbb A}]}.
			\]
		\end{enumerate}
	\end{remark}

	\begin{proposition}\label{lem:lifting}
		Forcing with $\col{\omega}{X_\infty}$ over $V[G_{\mathbb A}][G_{\mathbb B}]$ yields a
		$V$-generic filter $g^*$ on $\col{\omega}{X_\infty}$ such that $i_{\mu/g^*}$ lifts to
		an elementary embedding $i^*_{\mu/g^*}$ with domain $V[G_{\mathbb A}]$.
	\end{proposition}
	\begin{proof}
		Let $g\subseteq\col{\omega}{X_\infty}$ be generic over $V[G_{\mathbb A}][G_{\mathbb B}]$.
		By \cref{lem:zf-collapse-absorption} and
		\textup{(\ref{itm:lifting-surjection-onto-forcing})}, we can fix in $V$ a forcing
		equivalence $\varepsilon$ between
		$\mbb A\times\col{\omega}{X_\infty}$ and $\col{\omega}{X_\infty}$.
		Let $g^*\subseteq\col{\omega}{X_\infty}$ be the $V$-generic filter
		induced by $G_{\mathbb A}\times g$ under $\varepsilon$, so that $V[G_{\mathbb A}][g]=V[g^*]$.
		We will construct in $V[G_{\mathbb A}][G_{\mathbb B}][g]$ an $M_{g^*}$-generic filter
		$G_{\mathbb A}^*\subseteq i_{\mu/g^*}(\mbb A)$ containing $i_{\mu/g^*}[G_{\mathbb A}]$. By the standard lifting criterion, this gives an elementary embedding
		\[
		i^*_{\mu/g^*}:V[G_{\mathbb A}]\longrightarrow M_{g^*}[G_{\mathbb A}^*]
		\]
		defined by $i^*_{\mu/g^*}(\dot x^{G_{\mathbb A}})=i_{\mu/g^*}(\dot x)^{G_{\mathbb A}^*}$, with
		$i^*_{\mu/g^*}\restriction V=i_{\mu/g^*}$. Thus our goal is to construct such a filter $G^*_{\mathbb A}$ using a master condition argument.
		
		\begin{claim}\label{clm:lifting-countable-sequences-ordinals}
			$V[g^*]\cap\ord^\omega=V[G_{\mathbb A}][G_{\mathbb B}][g]\cap\ord^\omega$.
		\end{claim}
		\begin{proof}
			Working in $V[G_{\mathbb A}][G_{\mathbb B}]$,
			let $\dot \phi\in V[G_{\mathbb A}][G_{\mathbb B}]$ be a $\col{\omega}{X_\infty}$-name, and
			let $p_0\in g$ force that $\dot \phi:\omega\to\ord$ is a function.
			For each $n<\omega$, let $A_n\subseteq \col{\omega}{X_\infty}$ be a maximal antichain of conditions below $p_0$  deciding the value of $\dot \phi(n)$. Such maximal antichains exist since $\col{\omega}{X_\infty}$ is well-ordered in   $V[G_{\mathbb A}]$, in order type $\lambda$. Let 
			\[
			\dot{\phi}_*\df \{ (p,\check{\langle n,\alpha\rangle}): n<\omega,~p\in A_n,~p\forces{}{\dot{\phi}(n)=\alpha}\}. 
			\]
			By \cref{rem:lifting-preliminaries}~\textup{(\ref{itm:lifting-no-new-subsets})}, $ \dot{\phi}_*\in V[G_{\mathbb A}]$. Note that every $A_n$ is also a maximal antichain in $V[G_{\mathbb A}]$.
			It is clear that for every $n<\omega$ and every ordinal $\alpha$, \[
			\dot{\phi_*}^{g}(n)=\alpha\Longleftrightarrow \dot{\phi}^{g}(n)=\alpha.\] 
			Therefore, $\dot{\phi}^{g}\in V[G_\mbb A][g]\cap\ord^\omega=V[g^*]\cap\ord^\omega$.
		\end{proof}
		
		By \textup{(\ref{itm:lifting-ultrapower-hypotheses})} and
		\textup{(\ref{itm:lifting-countable-target-closure})}, together with
		\cref{clm:lifting-countable-sequences-ordinals,cor:omega1-computation}, we have
		\begin{equation}\label{eq:omega1}
			\omega_1^{M_{g^*}}=\Theta_{X_\infty}^V=\omega_1^{V[g^*]}=\omega_1^{V[G_{\mathbb A}][G_{\mathbb B}][g]}.
		\end{equation}
		
		Throughout this proof, we denote
		this ordinal simply by $\omega_1$.
		
		By \textup{(\ref{itm:lifting-countable-closure})}, $\mbb A$ is countably closed, so elementarity and \textup{(\ref{itm:lifting-countable-target-closure})}
		give the countable closure of $i_{\mu/g^*}(\mbb A)$ in $V[g^*]$.
		The forcing $\col{\omega}{X_\infty}$ and the surjection $c$ from
		\textup{(\ref{itm:lifting-surjection-onto-forcing})} make $G_{\mathbb A}$ countable there.
		Thus $i_{\mu/g^*}[G_{\mathbb A}]$ is a countable directed subset of
		$i_{\mu/g^*}(\mbb A)$, so countable closure gives a common lower bound for $i_{\mu/g^*}[G_{\mathbb A}]$.
		We use such a bound to start a sequence which is inductively constructed and meets the relevant dense subsets in $\mcl P(i_{\mu/g^*}(\mathbb A))$ belonging to $M_{g^*}$.
		
		\begin{claim}\label{clm:lifting-target-power-set-coding}
			In $V[g^*]$, there is a surjection $\pi:\mcl P(X_\infty)^V \twoheadrightarrow\mcl P(i_{\mu/g^*}(X_\infty))^{M_{g^*}}$.
		\end{claim}
		\begin{proof}
			We work in $V$ and fix a bijection $b:X_\infty^\omega\to X_\infty$.
			For $A\in\mcl P(X_\infty)$, define a map 
			$F_A:X_\infty^\omega\to\mcl P(X_\infty)$ in $V$ by
			\[
			F_A(f)=\{x\in X_\infty:  b(\langle x\rangle^\frown f)\in A\}.
			\]
			Now working in $V[g^*]$, we  define $\pi$ by 
			\[
			\pi(A)=[F_A]_{\mu/g^*},\]
			identifying ultrapower classes
			with their transitive collapses. Since $\mu/g^*$   has  the \Los{} property by  \textup{(\ref{itm:lifting-countable-target-closure})}, we have
			$\pi(A)\in\mcl P(i_{\mu/g^*}(X_\infty))^{M_{g^*}}$ for every
			$A\in\mcl P(X_\infty)^V$.
			
			To see that $\pi$ is surjective, 
			let $B\in\mcl P(i_{\mu/g^*}(X_\infty))^{M_{g^*}}$ and let 
			$F:(X_\infty^\omega)^V\to\mcl P(X_\infty)$ in $V$ be a representative in the $\mu/g^*$-ultrapower for $B$. Let
			\[
			A\df\{b(\langle x\rangle^\frown f):f\in(X_\infty^\omega)^V \land x\in F(f)\}.
			\]
			Then $A\in\mcl P(X_\infty)$. Since the map
			$(x,f)\mapsto b(\langle x\rangle^\frown f)$ is injective, we have $F_A=F$. Hence $\pi(A)=B$, as required.
		\end{proof}
		
		The surjection $c:X_\infty\twoheadrightarrow\mathbb A$ from \textup{(\ref{itm:lifting-surjection-onto-forcing})} induces in $V$ a surjection $\bar c:\mcl P(X_\infty)\twoheadrightarrow\mcl P(\mathbb A)$. By elementarity, $i_{\mu/g^*}(\bar c)$ maps $\mcl P(i_{\mu/g^*}(X_\infty))^{M_{g^*}}$ onto $\mcl P(i_{\mu/g^*}(\mathbb A))^{M_{g^*}}$. By \textup{(\ref{itm:lifting-target-enumeration})} and \cref{eq:omega1}, $\tau$ maps $\omega_1$ onto $\mcl P(X_\infty)^V$. Hence, by \cref{clm:lifting-target-power-set-coding}, we obtain the surjection
		\[
			s:\omega_1\twoheadrightarrow\mcl P(i_{\mu/g^*}(\mathbb A))^{M_{g^*}},
			\qquad s(\alpha)\df i_{\mu/g^*}(\bar c)\bigl(\pi(\tau(\alpha))\bigr).
		\]
		For every $\alpha<\omega_1$, \textup{(\ref{itm:lifting-target-enumeration})} gives $\tau\restriction\alpha\in V[G_{\mathbb A}]\subseteq V[g^*]$. Since $\pi$ and $i_{\mu/g^*}(\bar c)$ belong to $V[g^*]$, we also have $s\restriction\alpha\in V[g^*]$.

		Well-order the conditions of $i_{\mu/g^*}(\mathbb A)$ by the first occurrence of their singletons in $s$. Thus $p\lhd q$ if and only if
		\[
			\min\{\xi<\omega_1:s(\xi)=\{p\}\}<\min\{\xi<\omega_1:s(\xi)=\{q\}\}.
		\]
		Let $s^*$ enumerate the conditions in this order.

		\begin{claim}\label{claim:i(A)-countably-closed}
			The forcing  $i_{\mu/g^*}(\mbb A)$ is countably closed in $V[G_{\mathbb A}][G_{\mathbb B}][g]$.
		\end{claim}
		\begin{proof}
			By \cref{clm:lifting-countable-sequences-ordinals} and \textup{(\ref{itm:lifting-countable-target-closure})}, every $\omega$-sequence of ordinals in $V[G_{\mathbb A}][G_{\mathbb B}][g]$ belongs to $M_{g^*}$. Since $M_{g^*}$ satisfies $\dc$ and $\omega_1^{M_{g^*}}=\omega_1$ by \cref{eq:omega1}, every countable subset of $\omega_1$ in this extension is bounded.

			Fix $\alpha<\omega_1$. Choose $\beta_\alpha<\omega_1$ strictly above the first occurrences in $s$ of the singletons $\{p\}$ for $p\in\rang{s^*\restriction\alpha}$. Applying the same enumeration rule to $s\restriction\beta_\alpha$ and restricting the result to $\alpha$ gives $s^*\restriction\alpha$. Thus $s^*\restriction\alpha\in V[g^*]$.

			Now let $\langle p_n:n<\omega\rangle$ be a decreasing sequence in $i_{\mu/g^*}(\mathbb A)$ in $V[G_{\mathbb A}][G_{\mathbb B}][g]$. For each $n<\omega$, let $\gamma_n$ be the rank of $p_n$ in $\lhd$, and put $\gamma=\sup_{n<\omega}(\gamma_n+1)<\omega_1$. By \cref{clm:lifting-countable-sequences-ordinals}, $\langle\gamma_n:n<\omega\rangle\in V[g^*]$. Since $p_n=(s^*\restriction\gamma)(\gamma_n)$, the sequence $\langle p_n:n<\omega\rangle$ belongs to $V[g^*]$ and hence to $M_{g^*}$ by \textup{(\ref{itm:lifting-countable-target-closure})}. There $i_{\mu/g^*}(\mathbb A)$ is countably closed, so the sequence has a lower bound.
		\end{proof}
		We now construct a decreasing sequence $\langle p_\beta:\beta<\omega_1\rangle$ in $i_{\mu/g^*}(\mathbb A)$. Let $q_0$ be the $\lhd$-least lower bound of $i_{\mu/g^*}[G_{\mathbb A}]$. At each stage $0<\beta<\omega_1$, let $q_\beta$ be the $\lhd$-least lower bound of $\langle p_\alpha:\alpha<\beta\rangle$, which exists by \cref{claim:i(A)-countably-closed}. For every $\beta<\omega_1$, if $s(\beta)$ is dense in $i_{\mu/g^*}(\mathbb A)$, let $p_\beta$ be the $\lhd$-least member of $s(\beta)$ below $q_\beta$. Otherwise, let $p_\beta=q_\beta$.

		Let $G^*_{\mathbb A}$ be the filter generated by $\langle p_\beta:\beta<\omega_1\rangle$. Since $p_0\leq q_0$, this filter contains $i_{\mu/g^*}[G_{\mathbb A}]$. Every dense subset of $i_{\mu/g^*}(\mathbb A)$ in $M_{g^*}$ occurs as $s(\beta)$ for some $\beta<\omega_1$ and contains $p_\beta$. Thus $G^*_{\mathbb A}$ is $M_{g^*}$-generic, as required.
	\end{proof}

	Working  in $V[G_{\mathbb A}][G_{\mathbb B}]$, by  \cref{lem:lifting}, there is a
	$\col{\omega}{X_\infty}$-name $\dot{g}^*$ such that
	$\1{\col{\omega}{X_\infty}}$ forces that
	$\dot{g}^*\subseteq\col{\omega}{X_\infty}$ is $V$-generic and
	$i_{\mu/\dot{g}^*}$ lifts to an elementary embedding
	$i^*_{\mu/\dot{g}^*}$ with domain $V[G_{\mathbb A}]$. 
	\begin{definition}
		In $V$, define $\mathrm{r}:(X_\infty^\omega)\to Z$ by $\mathrm{r}(f)\df\rang{f}$. 
	\end{definition}
	
	\begin{lemma}\label{lem:seed}
		Fix a $V[G_\mbb A][G_{\mbb B}]$-generic filter $g$. Let $g^*\df (\dot{g}^*)^g$.
		Then $[\mathrm{r}]_{\mu/g^*}=i^*_{\mu/g^*}[X_\infty]$, where the equivalence classes are identified with classes in their transitive collapses.
	\end{lemma}
	\begin{proof}
		Let us set $D\df(X_\infty^\omega)^V$.
		Let  $e\df\bigcup g^*$ be the generic coding induced by $g^*$. Also set
		$s\df[\id{D}]_{\mu/g^*}$. By the \Los{} property, $s$ is a function with domain $\omega$.
		For each $n<\omega$, the basic open set determined by $e\restriction(n+1)$
		belongs to $\mu/g^*$, and every $f$ in this set satisfies $f(n)=e(n)$.
		Hence $s(n)=i^*_{\mu/g^*}(e(n))$. Since $e$ maps $\omega$ onto $X_\infty$, we obtain
		\[
		[\mathrm{r}]_{\mu/g^*}=\rang{s}=i^*_{\mu/g^*}[X_\infty].
		\]
	\end{proof}
	
	Define
	\[
	\Phi:\mcl P(\mcl P_{\omega_1}(X_\infty))\longrightarrow\ro(\col{\omega}{X_\infty})
	\]
	by
	\[
	\Phi(A)\df\left\|i^*_{\mu/\dot{g}^*}[X_\infty]\in
	i^*_{\mu/\dot{g}^*}(\check A)\right\|_{\ro(\col{\omega}{X_\infty})}.
	\]
	We now let
	\[
	\mcl I\df\ker(\Phi).
	\]
	Note that $\mcl I$ is a proper ideal.
	So the quotient homomorphism
	\[
\Phi^*:\mcl P(\mcl P_{\omega_1}(X_\infty))/\mcl I\to\ro(\col{\omega}{X_\infty})	\]
	is nontrivial and injective.
	\begin{lemma}\label{lem:completeness-phi}
		The map $\Phi^*$ is a $\lambda$-complete homomorphism.
	\end{lemma}
	\begin{proof}
		Note that $\Phi(\mcl P_{\omega_1}(X_\infty))=1$, since $\col{\omega}{X_\infty}$ makes $X_\infty$ countable and forces $M_{\dot{g}^*}$ to be countably closed.
	
		By elementarity, $\Phi$ is a Boolean homomorphism.
		By \cref{rem:lifting-preliminaries}, every sequence $\langle A_x:x\in X_\infty\rangle$ of subsets of $\mcl P_{\omega_1}(X_\infty)$ in $V[G_{\mbb A}][G_{\mbb B}]$ belongs to $V[G_{\mathbb A}]$. Note that $|X_\infty|=\lambda$. So we can apply  $i^*_{\mu/\dot{g}^*}$ to such sequences, so $\Phi^*(\nabla_{x\in X_\infty} [A_x]_{\mcl I})=\bigvee_{x\in X_\infty} \Phi^*([A_x]_\mcl I)$.
	\end{proof}

	\begin{lemma}\label{lem:basic-properties-of-I}
		The ideal  $\mcl I$ is countably complete, fine, normal, and uniform.
	\end{lemma}
	\begin{proof}
		Let $\langle A_n:n<\omega\rangle$ be a sequence of members of $\mcl I$ in $V[G_{\mathbb A}][G_{\mathbb B}]$. By \cref{rem:lifting-preliminaries}, the sequence belongs to $V[G_{\mbb A}]$. Then one can use the generic elementarity of $i^*_{\mu/\dot{g}^*}$ to conclude that \[
		\left\|i^*_{\mu/\dot{g}^*}[X_\infty]\notin i^*_{\mu/\dot{g}^*}(\bigcup_{n<\omega}A_n)\right\|=1.\]
		Hence  $\bigcup_{n<\omega}A_n\in\mcl I$.

		For fineness, fix $x\in X_\infty$ and put $C_x\df\{\sigma\in Z:x\in\sigma\}$.
		We can easily see that
		
		\[
		\left\|i^*_{\mu/\dot{g}^*}[X_\infty]\in i^*_{\mu/\dot{g}^*}(C_x)\right\|=1.\]
		Thus $Z\setminus C_x\in\mcl I$.
		
		To prove normality, let $U\in\mcl I^+$ and let $F$ be a function in $V[G_{\mbb A}][G_{\mbb B}]$
		with domain $U$ such that $\emptyset\ne F(\sigma)\subseteq\sigma$ for every $\sigma\in U$. As before, $F\in V[G_{\mbb A}]$.
		Write $s=i^*_{\mu/\dot{g}^*}[X_\infty]$. There is a condition $p\in\col{\omega}{X_\infty}$ forcing $s\in i^*_{\mu/\dot{g}^*}(U)$.
		By elementarity, $p$ forces that
		$\emptyset\ne i^*_{\mu/\dot{g}^*}(F)(s)\subseteq i^*_{\mu/\dot{g}^*}[X_\infty]$. Choose $q\leq p$ and $x\in X_\infty$ such that
		$q$ forces $i^*_{\mu/\dot{g}^*}(x)\in i^*_{\mu/\dot{g}^*}(F)(s)$. Set
		\[U_x\df\{\sigma\in U:x\in F(\sigma)\}.\]
		Then $q$ forces that $s\in i^*_{\mu/\dot{g}^*}(U_x)$. So $U_x\in\mcl I^+$, as required.
		
		To see that $\mcl I$ is uniform, suppose that $A\subseteq Z$ has cardinality $\nu<\lambda$ in $V[G_{\mbb A}][G_{\mbb B}]$.
		Fix in $V[G_{\mbb A}]$ a well-ordering of $X_\infty^\omega$, which has cardinality $\lambda$
		there by \cref{rem:lifting-preliminaries}. For each nonempty $\sigma\in A$,
		take its least enumeration in this well-ordering. Thus
		$|\bigcup A|<\lambda$, since $X_\infty$ is uncountable in $V[G_{\mbb A}][G_\mbb B]$.
		Choose $x\in X_\infty\setminus\bigcup A$. Then $A\subseteq Z\setminus C_x$,
		so $A\in\mcl I$ by fineness. 
	\end{proof}
	
	By \cref{lem:completeness-phi} and the fact that $\ro(\col{\omega}{X_\infty})$ has density at most $\lambda=|X_\infty|$, the homomorphism $\Phi^*$ is a complete Boolean homomorphism.
	
	Let $g\subseteq\col{\omega}{X_\infty}$ be generic over $V[G_{\mathbb A}][G_{\mathbb B}]$,
	and let $g^*=(\dot{g}^*)^g$ and
	$i^*_{\mu/g^*}$ be defined by the construction of
	$\Phi$ above.
	Let
	\[
	\nu_g\df
	\{A\in\mcl P(Z)^{V[G_{\mathbb A}][G_{\mathbb B}]}:
	i^*_{\mu/g^*}[X_\infty]\in i^*_{\mu/g^*}(A)\}.
	\]
	Then $\nu_g$ is a $V[G_{\mathbb A}][G_{\mathbb B}]$-generic filter on $\mcl I^+$, and hence an ultrafilter over $\mcl P(Z)^{V[G_{\mbb A}][G_{\mbb B}]}$.
	Moreover, \cref{lem:well-defined-embedding-hwo} gives
	\begin{equation}\label{eq:gen-filter-I}
	\nu_g\cap\mcl P(Z)^V=\mu
\qquad \text{ and } \qquad
\mcl I\cap\mcl P(Z)^V=\mcl P(Z)^V\setminus\mu.
	\end{equation}

	\begin{lemma}\label{lem:ideal-ultrapower-ord-cov}
		The ideal $\mcl I$   has the
		ordinal covering property relative to $V$.\footnote{See~\cref{def:ordinal-covering} for the definition of the ordinal covering property.}
	\end{lemma}
	\begin{proof}
		
		Fix an ordinal $\eta$, a function $F:Z\to\eta$ in $V[G_{\mbb A}][G_\mbb B]$, and a set  $S\in\mcl I^+$.
		Let $B\df\dot{B}^{G_\mbb A}$. By \cref{rem:lifting-preliminaries}, both $F$ and $S$ belong to $V[G_{\mbb A}]$.
		As $\Phi(S)>0$, there is a $\col{\omega}{X_\infty}$-generic filter $g$ over $V[G_{\mbb A}][G_\mbb B]$
		such that $i^*_{\mu/g^*}[X_\infty]\in i^*_{\mu/g^*}(S)$.

		Let $s\df [\mathrm{r}]_{\mu/g^*}$, and let  $\alpha\df i^*_{\mu/g^*}(F)(s)$. Since $i^*_{\mu/g^*}$ lifts $i_{\mu/g^*}$,
		by \cref{lem:agreement-on-ordinals} we have $\alpha<i^*_{\mu/g^*}(\eta)=i_\mu(\eta)$.
		Let $F_0:Z\to\eta$ in $V$ be a representation of $\alpha$ in the $\mu$-ultrapower.
		By \cref{lem:iso-hwo}, $[F_0]_\mu$ and $[F_0\circ \mathrm{r}]_{\mu/g^*}$ collapse
		to the same ordinal. Thus
		\[
		i^*_{\mu/g^*}(F_0)(s)=[F_0\circ \mathrm{r}]_{\mu/g^*}=\alpha=i^*_{\mu/g^*}(F)(s).
		\]
		Let $S'\df\{\sigma\in S:F(\sigma)=F_0(\sigma)\}$. Then $S'\in V[G_\mbb A][G_{\mbb B}]$,
		and by elementarity $s\in i^*_{\mu/g^*}(S')$. Hence $\Phi(S')>0$, so $S'\in\mcl I^+$, and $F\restriction S'=F_0\restriction S'$.
	\end{proof}
	\begin{lemma}\label{lem:ideal-ultrapower-comparison}
		We have
		$i_{\nu_g}\restriction\ord
		=i_\mu\restriction\ord
		=i_{\mu/g^*}\restriction\ord.$
	\end{lemma}
	\begin{proof}
		By \cref{lem:ideal-ultrapower-ord-cov} and the genericity of $\nu_g$, every ordinal-valued function on $Z$ in $V[G_{\mbb A}][G_{\mbb B}]$ agrees modulo $\nu_g$ with one in $V$.
		Since $\nu_g\cap\mcl P(Z)^V=\mu$, the predecessors of each  ordinal have the same order type in the two ultrapowers. This gives the first equality.
		The second equality follows from \cref{lem:agreement-on-ordinals} applied to $g^*$.
	\end{proof}
	
	\begin{lemma}\label{lem:ideal-ultrapower-precip}
		The ideal $\mcl I$ is precipitous.
	\end{lemma}
	
	\begin{proof} 
		Let $G_{\mcl I}$ be an $\mcl I$-generic ultrafilter over $V[G_{\mbb A}][G_{\mbb B}]$.
		By the second equality in \cref{eq:gen-filter-I}, we have $G_{\mcl I}\cap\mcl P(Z)^V=\mu$.
		As in the proof of \cref{lem:ideal-ultrapower-comparison}, ordinal covering and genericity give a representative in $V$ for every ordinal-valued $G_{\mcl I}$-class.
		Thus for every ordinal $\eta$, the classes of functions $Z\to\eta$ modulo $G_{\mcl I}$ have order type $i_\mu(\eta)$.
		The ordinals are therefore well-founded, and so is the generic $G_{\mcl I}$-ultrapower.
	\end{proof}
	\begin{corollary}\label{cor:embed-in-ro}
		There is an ideal $\mcl I$ on
		$\mcl P_{\omega_1}(X_\infty)$ in $V[G_{\mathbb A}][G_{\mathbb B}]$ such that
		\begin{enumerate}
			\item\label{itm:embed-in-ro-ideal} $\mcl I$ is a precipitous,
			$|X_\infty|$-dense, fine, normal, uniform ideal  with the ordinal covering property relative to $V$.
			\item\label{itm:embed-in-ro-embedding} The  Boolean algebra $\mcl P(\mcl P_{\omega_1}(X_\infty))/\mcl I$ completely
			embeds into $\ro(\col{\omega}{X_\infty})$.
			\item For every $V[G_{\mathbb A}][G_{\mathbb B}]$-generic ultrafilter $G_{\mcl I}\subseteq\mcl I^+$,
			the generic ultrapower embedding $i_{G_{\mcl I}}$
			is elementary and has critical point $\omega_1$. Moreover, the target of $i_{G_{\mcl I}}$ is countably closed in $V[G_{\mathbb A}][G_{\mathbb B}][G_{\mcl I}]$.
		\end{enumerate}
	\end{corollary}
	\begin{proof}
		By \cref{lem:completeness-phi} and the argument following \cref{lem:basic-properties-of-I}, $\Phi^*$ completely embeds $\mcl P(Z)/\mcl I$ into $\ro(\col{\omega}{X_\infty})$.
		 Therefore, $\mcl I$ is $|X_\infty|$-dense. 
		 The remaining ideal properties follow from \cref{lem:basic-properties-of-I,lem:ideal-ultrapower-ord-cov,lem:ideal-ultrapower-precip}.

		It remains to prove \textup{(3)}. By \textup{(\ref{itm:lifting-well-orderability-and-choice})} and \textup{(\ref{itm:lifting-countable-closure})},
		$V[G_{\mbb A}][G_{\mbb B}]\models\dc_\lambda$, hence $\ac_\lambda$.
		Let $G_{\mcl I}$ be $\mcl I$-generic over $V[G_{\mbb A}][G_{\mbb B}]$.
		Since $|Z|=\lambda$ by \cref{rem:lifting-preliminaries}, $\ac_\lambda$ gives the \Los{} property for $G_{\mcl I}$.
		By precipitousness, write
		\[
		i_{G_{\mcl I}}:V[G_{\mbb A}][G_{\mbb B}]\longrightarrow
		M_{G_{\mcl I}}\df\ult{V[G_{\mbb A}][G_{\mbb B}]}{G_{\mcl I}}
		\]
		for the collapsed ultrapower embedding.

		By \textup{(2)}, take $g\subseteq\col{\omega}{X_\infty}$ generic over $V[G_{\mbb A}][G_{\mbb B}]$ with $\nu_g=G_{\mcl I}$, and let $g^*=(\dot g^*)^g$.
		By \cref{lem:seed}, the factor map is
		\[
		\ult{V[G_{\mbb A}]}{G_{\mcl I}}\longrightarrow M_{g^*}[G_{\mbb A}^*],
		\qquad [F]_{G_{\mcl I}}\longmapsto i^*_{\mu/g^*}(F)([\mathrm r]_{\mu/g^*}).
		\]
		Its range contains $i^*_{\mu/g^*}[V[G_{\mbb A}]]$ and, by the proof of \cref{lem:ideal-ultrapower-ord-cov}, every ordinal.
		Since $(X_\infty^\omega)^V$ is well-orderable in $V[G_{\mbb A}]$, the range contains $[\id{(X_\infty^\omega)^V}]_{\mu/g^*}$, which generates $M_{g^*}$ over $i_{\mu/g^*}[V]$.
		It also contains $G_{\mbb A}^*=i^*_{\mu/g^*}(G_{\mbb A})$, so the factor map is onto.
		Since $|Z|=\lambda$ and $\mbb B$ adds no $\lambda$-sequences of elements of $V[G_{\mbb A}]$, the collapses agree on the smaller ultrapower. Thus
		\[
		M_{g^*}[G_{\mbb A}^*]=\ult{V[G_{\mbb A}]}{G_{\mcl I}}\subseteq M_{G_{\mcl I}},
		\qquad i_{G_{\mcl I}}\restriction V[G_{\mbb A}]=i^*_{\mu/g^*}.
		\]
		By \cref{rem:lifting-preliminaries,eq:omega1},
		\[i_{G_{\mcl I}}(\omega_1^V)=\omega_1^{V[g^*]}>\omega_1^V=\omega_1^{V[G_{\mbb A}][G_{\mbb B}]}.\]
		Hence $\crit(i_{G_{\mcl I}})=\omega_1^V$.

		Finally, let $x\in M_{G_{\mcl I}}^\omega\cap V[G_{\mbb A}][G_{\mbb B}][G_{\mcl I}]$.
		Working in $V[G_{\mbb A}][G_{\mbb B}]$, use \textup{(1)} and $\ac_\lambda$ to choose a family $\mcl F$ of at most $\lambda$ representatives on maximal antichains deciding the coordinates of a name for $x$.
		Since $|\mcl F\times Z|\leq\lambda$, fix a bijection
		\[
		e:\gamma\longrightarrow\{F(\sigma):F\in\mcl F,\ \sigma\in Z\}, 
		\]
		where $\gamma\leq\lambda$.
		Then $\rang{x}\subseteq\rang{i_{G_{\mcl I}}(e)}$, and by \cref{clm:lifting-countable-sequences-ordinals} and \textup{(\ref{itm:lifting-countable-target-closure})},
		\[
		i_{G_{\mcl I}}(e)^{-1}\circ x\in M_{g^*}\subseteq M_{G_{\mcl I}}.
		\]
		Hence $x\in M_{G_{\mcl I}}$.
	\end{proof}

	\section{Forcing generic supercompactness over minimal Nairian models}\label{sec:forcing-generic-supercompactness-and-gch}
	In this section we prove the following version of the main theorem under the assumption that $\mu$-Uniformization holds. We prove $\mu$-Uniformization in \cref{sec:uniformization-in-nairian-models}.
	
	Let us briefly recall that the theory $T_{\mathrm{IS}}$ is $\zfc$ together with the
	following assertion:
	\begin{quote}
		For every regular cardinal $\kappa\geq\omega_1$, there is a 
		precipitous, fine, normal, uniform, $|\kappa^\omega|$-dense ideal $\mcl I_\kappa$ on
		$\mcl P_{\omega_1}(\kappa^\omega)$ such that $\mcl I_\kappa^+$ witnesses that
		$\omega_1$ is generically $\kappa$-supercompact.
	\end{quote}
	
	In the following theorem, $\kappa_\alpha$ denotes the $\alpha$-th member of the $\Theta$-sequence computed in $N$ (see \cref{def:theta-sequence}).
	
	\begin{theorem}[$V=L(\mcl P(\R))+\hpc+\treg$]\label{thm:generic-supercompactness-and-gch}
		Suppose that $N$ is a minimal Nairian model in $V$. Assume that for every ordinal $\alpha\in N$ and every  supercompactness measure $\mu\in N$ on $\mcl P_{\omega_1}(\kappa_\alpha^\omega)^N$, $\mu$-Uniformization holds in $N$. Then $N$ has a class forcing extension in which $T_{\mathrm{IS}}+\gch$ holds.
	\end{theorem}
	
	We prove the theorem in three steps.
	Fix a minimal Nairian model $N$. This will be our ground model. Recall the notation
	\[
	\eta_\gamma\df\hH_\gamma^N\cap\ord.
	\]
	We work in $N$, and so we omit the superscript $N$ whenever it is clear. We shall apply \cref{thm:forcing-choice} with $\mathsf{M}=N$, $\alpha=0$, and
	the initial forcing $\mathbb P\df \col{\omega_1}{\R}$.
	Notice that by \cref{thm:consequences-of-nairian-models},
	$N\models{\ntbase}$, and by \cref{lem:singular-limit-heights},
	$\eta_\epsilon$ is singular in $N$ for every nonzero limit $\epsilon$.
	We will show that forcing with $\mathbb P$ followed by $\mathbb Q^{\textsf{ac}}$ as in
	\cref{thm:forcing-choice} yields an extension satisfying
	$T_{\mathrm{IS}}+\gch$.
	
	Fix an $N$-generic filter $G\subseteq\mathbb P$.
	\forcingstep{step:Preparation}{Preparation.}
	
	Our first task is to verify the hypotheses of \cref{thm:forcing-choice}
	for our initial forcing.
	\begin{lemma}\label{lem:colR-is-good}
		The structure $\hH_0[G]$ satisfies Choice and is completely ordinal-full\footnote{See~\cref{def:ordinal-fullness}.} in
		${N}[G]$. Moreover, $\omega_1$ is the largest cardinal of $\hH_0[G]$.
	\end{lemma}
	\begin{proof}
		In $N$, the forcing $\col{\omega_1}{\R}$ is countably closed, and $\dc$ holds  by
		\cref{thm:consequences-of-nairian-models}, so $\col{\omega_1}{\R}$
		preserves $\omega_1$ and adds no reals. Its conditions are hereditarily countable;
		hence fullness of $\hH_0$ (by \cref{ax:axiom-of-harmony}) gives
		$\col{\omega_1}{\R}\in\hH_0$.
		
		We now work in ${N}[G]$. The generic coding of 
		$\R$ induced by $\mbb P$ and the usual coding of
		subsets of $\omega$ by reals give
		$|\R|=|\powerset(\omega)^{N}|=\omega_1$.
		We then have, by \cref{prop:bls-h-like-preservation},
		\[
		\hH_0[G]=\mcl H_{\powerset(\omega)^{N}}=\mcl H_{\omega_1}.
		\]
		Every member of $\mcl H_{\omega_1}$ has a well-ordering in
		$\mcl H_{\omega_1}$, so $\hH_0[G]$ satisfies Choice and has largest
		cardinal $\omega_1$. We now show that $\hH_0[G]$ is completely ordinal-full in ${N}[G]$.
		Note that
		$|\omega_1^\omega|=\omega_1$ holds in ${N}[G]$.
		Thus, for $Y\in \mcl H_{\omega_1}$ and $\gamma<\omega_1$, the set
		$Y^\gamma$ has hereditary cardinality at most $\omega_1$ and hence
		belongs to $\mcl H_{\omega_1}$.  
	\end{proof}
	Now \cref{lem:colR-is-good} allows us to proceed with $\mathbb Q^{\textsf{ac}}\df\mathbb Q^{\textsf{ac}}(0,G)$ over ${N}[G]$.
	Let us recall the definition of $\mathbb Q^{\textsf{ac}}$ for completeness. In ${N}[G]$, $\mathbb Q^{\textsf{ac}}$ is the full-support class iteration
	of length $\ord$ with $\mathbb Q^{\textsf{ac}}_{<0}=\{1\}$ and
	\[
	\mathbb Q^{\textsf{ac}}_{<\beta+1}=\mathbb Q^{\textsf{ac}}_{<\beta}\ast\dot{\mathbb Q}^{\textsf{ac}}(\beta).
	\]
	The initial segment $\mathbb Q^{\textsf{ac}}_{<\epsilon}$ at a limit ordinal $\epsilon$
	is taken with full support. Set
	\[
	\mathbb Q^{\textsf{ac}}(0)=\col{\eta_0}{\hH_0}.
	\]
	At each positive stage $\beta$, let $H_\beta\subseteq\mathbb Q^{\textsf{ac}}_{<\beta}$
	be generic over ${N}[G]$. In ${N}[G][H_\beta]$, define the next forcing
	as follows:
	\begin{enumerate}
		\item If $\beta$ is a nonzero limit, $\mathbb Q^{\textsf{ac}}(\beta)$ is trivial.
		\item If $\beta=1$ or $\beta$ is the successor of a successor ordinal, let
		\[
		\mathbb Q^{\textsf{ac}}(\beta)=\col{\eta_{\beta}}{\hH_{\beta}}.
		\]
		\item If $\beta=\delta+1$ for a nonzero limit ordinal $\delta$, put
		$\upsilon\df(\eta_{\delta}^+)^{{N}[G][H_\beta]}$ and let
		\[
		\mathbb Q^{\textsf{ac}}(\beta)=
		\col{\upsilon}{\mcl P(\hH_{\delta})^N}
		\ast\dotcol{\eta_{\delta+1}}{\hH_{\delta+1}}.
		\]
		
	\end{enumerate}
	Let $H\subseteq\mathbb Q^{\textsf{ac}}$ be an ${N}[G]$-generic filter, and
	write $H_\gamma$ for its restriction to $\mathbb Q^{\textsf{ac}}_{<\gamma}$.
	By \cref{thm:forcing-choice}, 
	\[
	{N}[G][H]\models{\zfc}.
	\]
	\forcingstep{step:Cardinal arithmetic}{Cardinal arithmetic.}
	
	Our next task is to verify that $\gch$ holds in ${N}[G][H]$.
	For this, we again rely on the results of Blue, Larson, and Sargsyan
	\cite{NM2} recalled in \cref{subsec:basic-nairian-theory}.
	
	Since $\omega_1$ is the largest cardinal of $\hH_0[G]$ by \cref{lem:colR-is-good},
	\cref{prop:bls-stage-choice}\textup{(\ref{itm:bls-initial-dc})}
	shows that $\omega_1$ remains a cardinal in ${N}[G][H]$.
	\begin{lemma}\label{lem:bls-countable-closure}
		Let $\beta\in\ord$. Then
		$\mathbb P\ast\dot{\mathbb Q}^{\textsf{ac}}_{<\beta}$
		is countably closed.
	\end{lemma}
	\begin{proof}
		The forcing $\mathbb P$ is countably closed in $N$, so we may assume $\beta>0$.
		It is enough to show that the tail forcing $\mathbb Q^{\textsf{ac},1}$ is
		countably closed in ${N}[G][H_1]$.
		
		In ${N}[G]$, $\eta_0$ is regular by
		\cref{prop:bls-h-like-preservation}, applied to $\hH_0$.
		Thus $\mathbb Q^{\textsf{ac}}(0)$ is countably closed.
		By \cref{prop:bls-stage-choice}\textup{(\ref{itm:bls-completed-stage})},
		${N}[G][H_1]$ satisfies $\dc_{\eta_0}$, and
		$\eta_0$ is the largest cardinal of $\hH_1[G][H_1]$.
		By \cref{prop:bls-tail-closure}, applied with $\alpha=0,\beta=1$, and
		$\kappa=\eta_0$, the tail forcing $\mathbb Q^{\textsf{ac},1}$ is
		${<}\eta_0^+$-closed in ${N}[G][H_1]$.
	\end{proof}
	
	By the proof of \cref{lem:colR-is-good}, $\ch$ holds in $\hH_0[G]$.
	Thus, by \cref{lem:bls-countable-closure},
	\begin{equation}\label{eq:CH}
		{N}[G][H]\models{\ch}.
	\end{equation}
	By \cref{lem:bls-countable-closure} and $\dc$ in $N$, we have
	\[
	{N}[G][H]\cap\ord^\omega=N\cap \ord^\omega.
	\]
	Moreover, by \cref{prop:bls-stage-choice}\textup{(\ref{itm:bls-initial-dc})},
	$\omega_1$ remains a cardinal in ${N}[G][H]$, where $\eta_0=\omega_2$ and
	\begin{equation}\label{eq:bls-initial-cardinal}
		{N}[G][H]\mdls{2^{\omega_1}=\omega_2}.
	\end{equation}
	
	By induction using
	conclusion~\textup{(\ref{itm:forcing-choice-successor-heights})} of
	\cref{thm:forcing-choice} and the continuity of the height sequence $\delta\mapsto\eta_\delta$, we
	have, in ${N}[G][H]$,  for every $n<\omega$
	\[
	\eta_n=\aleph_{n+2}
	\]
	and, for every nonzero limit ordinal $\delta$ and every $n<\omega$,
	\[
	\eta_\delta=\aleph_\delta
	\qquad \text{ and } \qquad
	\eta_{\delta+n+1}=\aleph_{\delta+n+2}.
	\]
	In particular, ${N}[G][H]$ has no weakly inaccessible cardinals,
	and every singular cardinal of ${N}[G][H]$ is already singular in
	$N$, by \cref{lem:singular-limit-heights}.
	\begin{lemma}\label{lem:cardinals-in-extension}
		Let $\kappa\geq\omega_2$ be a regular cardinal in ${N}[G][H]$. Then exactly one of
		the following holds:
		\begin{enumerate}
			\item $\kappa=\eta_{\beta}$, where $\beta=0$ or $\beta$ is a
			successor ordinal;
			\item $\kappa=\eta_{\beta}^+$, where $\beta$ is a
			nonzero limit ordinal.
		\end{enumerate}
	\end{lemma}
	\begin{proof}
		There is a least $\gamma$ such that
		$\kappa<\eta_{\gamma}$. As $\eta_0=\omega_2\leq\kappa$, we have
		$\gamma>0$, and continuity of the  sequence $\delta\mapsto \eta_\delta$ rules out a limit $\gamma$. So let us write
		$\gamma=\beta+1$. We then have
		\[
		\eta_{\beta}\leq\kappa<\eta_{\beta+1}.
		\]
		If $\beta=0$ or $\beta$ is a successor, then by
		conclusion~\textup{(\ref{itm:forcing-choice-successor-heights})} of
		\cref{thm:forcing-choice}, we have
		$\eta_{\beta+1}=\eta_{\beta}^+$, and consequently $\kappa=\eta_{\beta}$. 
		If $\beta$ is a nonzero limit, we have
		$\eta_{\beta+1}=\eta_{\beta}^{++}$, so $\kappa$ is either
		$\eta_{\beta}$ or $\eta_{\beta}^+$. By
		\cref{lem:singular-limit-heights}, $\eta_{\beta}$ is singular in $N$,
		and by conclusion~\textup{(\ref{itm:forcing-choice-cardinals-gch})} of
		\cref{thm:forcing-choice}, it remains a cardinal in ${N}[G][H]$.
		Hence it remains singular there. Since $\kappa$ is regular,
		we must have $\kappa=\eta_{\beta}^+$.
	\end{proof}
	
	\begin{lemma}\label{lem:gch-holds}
		${N}[G][H]\models{\gch}$.
	\end{lemma}
	\begin{proof}
		By \cref{eq:CH} and \cref{eq:bls-initial-cardinal}, $\gch$ holds at
		$\omega$ and $\omega_1$, and
		conclusion~\textup{(\ref{itm:forcing-choice-cardinals-gch})} of
		\cref{thm:forcing-choice} gives
		$\gch$ at every $\eta_{\beta}$.
		By the preceding cardinal computations, the only remaining
		cardinals are $\eta_{\delta}^+$ for nonzero limit ordinals $\delta$.
		\Cref{prop:bls-limit-stages}\textup{(\ref{itm:bls-second-limit-factor})}
		gives $\gch$ at each such cardinal in ${N}[G][H_{\delta+2}]$.
		By \cref{prop:bls-tail-closure}, the tail $\mathbb Q^{\textsf{ac},\delta+2}$
		adds no subsets of this cardinal and preserves its successor,
		so $\gch$ at $\eta_\delta^+$ remains true in ${N}[G][H]$.
	\end{proof}
	\forcingstep{step:Getting the ideals}{Getting the ideals.}
	
	Our final task is to show that ${N}[G][H]$ satisfies $T_{\mathrm{IS}}$.
	We first prove two lemmas needed to apply the results of
	\cref{sec:lifting-supercompactness-embeddings}.
	
	\begin{notation}\label{not:beta_and_alpha}
		Let $\kappa\geq\omega_1$ be a regular cardinal in ${N}[G][H]$.
		
		\begin{enumerate}
			\item Let
			$\alpha(\kappa)$ be the least member of $\{-1\}\cup\ord$ such that
			$\kappa<\kappa_{\alpha(\kappa)+1}$.
			\item
			Set $\beta(\omega_1)\df-1$. For $\kappa\geq\omega_2$,
			let
			$\beta(\kappa)$ be determined by
			\[
			\begin{cases}
				\kappa=\eta_{\beta(\kappa)},
				& \text{if $\beta(\kappa)=0$ or $\beta(\kappa)$ is a successor},\\
				\kappa=(\eta_{\beta(\kappa)}^+)^{{N}[G][H]},
				& \text{if $\beta(\kappa)$ is a nonzero limit}.
			\end{cases}
			\]
		\end{enumerate}
	\end{notation}
	\begin{remark}
		Note that $\beta(\kappa)$ is unique by \cref{lem:cardinals-in-extension}.
	\end{remark}
	
	\begin{lemma}\label{lem:initial-forcing-coding}
		Let $\kappa\geq\omega_1$ be a regular cardinal in ${N}[G][H]$.
		Then $\eta_{\beta(\kappa)+1}=\kappa_{\alpha(\kappa)+1}$.
		
	\end{lemma}
	\begin{proof}
		Let us put $\beta\df\beta(\kappa)$ and $\alpha\df\alpha(\kappa)$.
		If $\beta=-1$, then $\alpha=-1$ and
		$\eta_0=\kappa_0$.
		If $\beta=0$ or $\beta$ is a successor, then
		$\eta_{\beta}$ and $\eta_{\beta+1}$ are consecutive
		terms of the $\Theta$-sequence by \cref{lem:heights-and-theta-sequence}.
		Since $\kappa=\eta_{\beta}$, the definition of $\alpha$
		gives $\eta_{\beta+1}=\kappa_{\alpha+1}$.
		
		Suppose now that $\beta$ is a nonzero limit.
		Using \cref{lem:heights-and-theta-sequence} and
		conclusion~\textup{(\ref{itm:forcing-choice-successor-heights})} of
		\cref{thm:forcing-choice}, we obtain
		\[
		\kappa_{\beta}=\eta_{\beta}
		<\kappa<\eta_{\beta+1}.
		\]
		If $\cof{\beta}=\omega$, then
		\cref{lem:heights-and-theta-sequence} gives
		$\eta_{\beta+1}=\kappa_{\beta+1}$,
		so $\alpha=\beta$.
		
		Assume now that $\cof{\beta}>\omega$.
		Since $\cof{\kappa_{\beta}}=\cof{\beta}$ is uncountable,
		every countable sequence in $\kappa_{\beta}$ is bounded.
		As in the proof of \cref{lem:heights-and-theta-sequence}, this gives
		$\kappa_{\beta}^\omega\subseteq\hH_{\beta}$.
		The bijection supplied by
		\cref{prop:bls-limit-stages}\textup{(\ref{itm:bls-limit-size})}
		therefore gives $|\kappa_{\beta}^\omega |
		=\kappa_{\beta}<\kappa$ in $N[G][H]$.
		Hence no surjection from $\kappa_{\beta}^\omega$ onto $\kappa$
		belongs to $N$.
		Since $\kappa_{\beta+1}=\Theta_{\kappa_{\beta}^\omega}$,
		we obtain $\kappa_{\beta+1}\leq\kappa$.
		Together with
		$\eta_{\beta+1}=\kappa_{\beta+2}$ from
		\cref{lem:heights-and-theta-sequence}, this yields
		$\alpha=\beta+1$.
		
		Thus $\eta_{\beta+1}=\kappa_{\alpha+1}$.  
	\end{proof}
	Recall $\dot{\mathbb Q}^{\textsf{ac}}(\beta+1)_0$ and $\dot{\mathbb Q}^{\textsf{ac}}(\beta+1)_1$ from \cref{not:factors}.
	\begin{lemma}\label{lem:initial-forcing-choice}
		Let $\kappa\geq\omega_1$ be a regular cardinal in ${N}[G][H]$. 
		Let $K$ be the generic filter on $\mathbb P\ast\dot{\mathbb Q}^{\textsf{ac}}_{<\beta(\kappa)+1}
		\ast\dot{\mathbb Q}^{\textsf{ac}}(\beta(\kappa)+1)_0$  induced by
		$G$, $H_{\beta(\kappa)+1}$, and $H(\beta(\kappa)+1)_0$.
		Then $\hH_{\beta(\kappa)+1}[K]$ satisfies Choice and $\dc_\kappa$
		and has $\kappa$ as its largest cardinal. Moreover,
		\[
		{N}[K]\mdls{
			|\kappa_{\alpha(\kappa)}^\omega|=\kappa\ \text{and}\ \dc_\kappa}.
		\]
		
	\end{lemma}
	\begin{proof}
		Put $\alpha\df\alpha(\kappa)$ and $\beta\df\beta(\kappa)$.
		For $\beta=-1$, \cref{lem:colR-is-good} gives Choice
		in $\hH_0[G]$, with largest cardinal $\omega_1$, and
		\cref{prop:bls-stage-choice}\textup{(\ref{itm:bls-initial-dc})}
		gives $\dc_{\omega_1}$ in ${N}[G]$.
		If $\beta=0$ or $\beta$ is a successor, the first factor
		$\mathbb Q^{\textsf{ac}}(\beta+1)_0$ is trivial, and the assertions about
		Choice, the largest cardinal, and Dependent Choice follow from
		\cref{prop:bls-stage-choice}\textup{(\ref{itm:bls-completed-stage})}.
		If $\beta$ is a nonzero limit, then
		\cref{prop:bls-limit-stages}\textup{(\ref{itm:bls-limit-size})}
		gives
		\[
		(\eta_{\beta}^+)^{{N}[G][H_{\beta+1}]}
		=(\eta_{\beta}^+)^{{N}[G][H]}=\kappa.
		\]
		The assertions now follow from
		\cref{prop:bls-limit-stages}\textup{(\ref{itm:bls-first-limit-factor})}.
		
		By \cref{lem:initial-forcing-coding} and
		\cref{thm:theta-sequence-properties}\textup{(\ref{itm:theta-sequence-subsets})},
		\[
		\kappa_\alpha^\omega\in\mcl P(\kappa_\alpha^\omega)
		\subseteq\hH|\kappa_{\alpha+1}
		=\hH_{\beta+1}.
		\]
		Choice and the maximality of $\kappa$ in
		$\hH_{\beta+1}[K]$ therefore give an injection
		$\kappa_\alpha^\omega\rightarrow\kappa$ in this model.
		Conversely, since
		$\kappa<\kappa_{\alpha+1}=\Theta_{\kappa_\alpha^\omega}$,
		there is a surjection $\kappa_\alpha^\omega\twoheadrightarrow\kappa$
		in $N$.
		As $\kappa$ is a cardinal in ${N}[K]$, these two maps
		show that $|\kappa_\alpha^\omega|=\kappa$ there.
		
	\end{proof}
	
	\begin{lemma}\label{lem:ideals-in-final-extension}
		${N}[G][H]\models{T_{\mathrm{IS}}}$.
	\end{lemma}
	\begin{proof}
		
		Fix a regular cardinal $\kappa\geq\omega_1$ in ${N}[G][H]$.
		By \cref{lem:gch-holds}, $|\kappa^\omega|=\kappa$ in this extension.
		By \cref{lem:bls-countable-closure}, the forcing adds no countable
		sequences of ordinals. Hence
		\[
		\omega_1=\omega_1^{{N}[G][H]}\qquad \text{and }\qquad
		\kappa^\omega=(\kappa^\omega)^{{N}[G][H]}.
		\]
		In particular, $\col{\omega}{\kappa^\omega}$ denotes the same forcing in
		these two models.
		Set $\alpha\df\alpha(\kappa)$ and $\beta\df\beta(\kappa)$, using
		\cref{not:beta_and_alpha}. Also set
		\[
		\begin{aligned}
			\mathbb A&\df \mathbb P\ast\dot{\mathbb Q}^{\textsf{ac}}_{<\beta+1}
			\ast\dot{\mathbb Q}^{\textsf{ac}}(\beta+1)_0,&\dot{\mathbb B}&\df\dot{\mathbb Q}^{\textsf{ac}}(\beta+1)_1,\\
			G_{\mathbb A}&\df G\ast\dot{H}_{\beta+1}\ast\dot{H}(\beta+1)_0,&G_{\mathbb B}&\df H(\beta+1)_1,\\
			W_1&\df{N}[G_{\mathbb A}],&W_2&\df{N}[G_{\mathbb A}][G_{\mathbb B}].
		\end{aligned}
		\]
		
		The forcing from $W_1$ to $W_2$ is the second factor at coordinate
		$\beta+1$. Evaluating its name in $W_1$, we obtain
		\[
		\mathbb B\df\dot{\mathbb B}^{G_{\mathbb A}}
		=\col{\eta_{\beta+1}}
		{\hH_{\beta+1}}^{W_1}.
		\]
		Forcing over $W_2$ with the
		tail $\mathbb Q^{\textsf{ac},\beta+2}$ yields ${N}[G][H]$.
		
		We will apply the results of \cref{sec:lifting-supercompactness-embeddings}. We take
		\[
		W\df N,\qquad X\df\kappa_\alpha,\qquad
		X_\infty\df\kappa_\alpha^\omega,\qquad
		Z\df\mcl P_{\omega_1}(X_\infty).
		\]
		By \cref{thm:consequences-of-nairian-models}, $N$ satisfies $\zf+\ad+\dc$ and
		contains a supercompactness measure $\mu$ on $Z$. Fix such a measure.
		Clearly, $X\in\hwo$.
		By \cref{thm:theta-sequence-properties}\textup{(\ref{itm:theta-sequence-theta-and-regularity})},
		$\Theta_{X_\infty}=\kappa_{\alpha+1}$.
		In particular, when $\kappa=\omega_1$, we have $\alpha=-1$,
		$X=\omega$, and $X_\infty=\R$, so $\Theta_{X_\infty}=\kappa_0=\Theta=\eta_0$.
		For $\alpha\geq0$, \cref{thm:image-of-omega-one} also gives
		$i_\mu(\omega_1^{N})=\kappa_{\alpha+1}$.
		
		We now verify assumptions
		\textup{(\ref{itm:lifting-ultrapower-hypotheses})}--\textup{(\ref{itm:lifting-target-enumeration})}
		of \cref{sec:lifting-supercompactness-embeddings}, using
		$\mathbb A$ followed by $\dot{\mathbb B}$.
		\begin{enumerate}[label=(S\arabic*),ref=S\arabic*]
			\item\label{itm:main-theta-regularity} $\Theta_{X_\infty}$ is regular
			by \cref{thm:theta-sequence-properties}\textup{(\ref{itm:theta-sequence-theta-and-regularity})}.
			\item\label{itm:main-ultrapower-properties}  
			For every
			$N$-generic filter $g$ on $\col{\omega}{X_\infty}$,
			$\mu/g$ has the \Los{} property (by \cref{prop:los-property}) and its transitive ultrapower
			$\ult{N}{\mu/g}$ is countably closed in ${N}[g]$
			(by \cref{prop:countable-closure-of-ultrapower}).
			
			\item\label{itm:main-forcing-coding} The forcing $\mathbb A$ is a surjective image of $X_\infty$.
			\begin{proof}
				It follows from \cref{prop:bls-bounded-presentations}\textup{(\ref{itm:bls-successor-presentation})},
				\textup{(\ref{itm:bls-limit-presentation})} when
				$\beta\geq0$, and from the proof of \cref{lem:colR-is-good} when
				$\beta=-1$, that $\mathbb A\in\hH_{\beta+1}$. By
				\cref{lem:initial-forcing-coding}, $\eta_{\beta+1}=\kappa_{\alpha+1}$.
				Therefore,
				\[
				\mathbb A\in\hH_{\beta+1}
				=\hH|\kappa_{\alpha+1}=\mcl H_{\kappa_{\alpha}^\omega}.
				\]
			\end{proof} 
			
			\item\label{itm:main-well-orderability}
			By \cref{lem:initial-forcing-choice},
			\[
			W_1\mdls{|X_\infty|=\kappa\ \text{and}\ \dc_\kappa}.
			\]
			
			\item\label{itm:main-forcing-closure}
			The forcing $\mathbb A\ast\dot{\mathbb B}$ is countably closed,
			and $\mathbb B$ is
			${<}(\kappa^+)^{W_1}$-closed in $W_1$.
			\begin{proof}
				The countable closedness of $\mbb A$ follows from \cref{lem:bls-countable-closure}.
				By conclusion~\textup{(\ref{itm:forcing-choice-successor-heights})} of
				\cref{thm:forcing-choice} and
				\cref{prop:bls-stage-choice}\textup{(\ref{itm:bls-initial-dc})},
				\[
				\eta_{\beta+1}=(\kappa^+)^{{N}[G][H]}.
				\]
				Note that $\eta_{\beta+1}$ is regular in ${N}[G][H]$ and
				hence also in $W_1$, so $(\kappa^+)^{W_1}\leq\eta_{\beta+1}$.
				Thus the collapse $\mathbb B$ is ${<}\eta_{\beta+1}$-closed in $W_1$, and therefore
				${<}(\kappa^+)^{W_1}$-closed.
			\end{proof}
			\item\label{itm:main-target-enumeration}
			The filter $G_{\mathbb A}\ast G_{\dot{\mathbb B}}$ is $N$-generic for
			$\mathbb A\ast\dot{\mathbb B}$, and there is a surjection
			\[
			\tau:\Theta_{X_\infty}\twoheadrightarrow\mcl P(X_\infty)
			\]
			in $W_2$  with
			$\tau\restriction\gamma\in W_1$ for every $\gamma<\Theta_{X_\infty}$.
			\begin{proof}
				The forcing $\mathbb B$ adds, over $W_1$, a bijection
				$f:\eta_{\beta+1}\to\hH_{\beta+1}$.
				By \cref{lem:initial-forcing-coding}, we have $\eta_{\beta+1}=\Theta_{X_\infty}$, and by
				\cref{thm:theta-sequence-properties}\textup{(\ref{itm:theta-sequence-subsets})},
				we have $\mcl P(X_\infty)\subseteq\hH_{\beta+1}$.
				Thus $f$ induces a surjection $\tau\in W_2$ as above.
				By \cref{lem:initial-forcing-choice}, $W_1$ satisfies $\dc_\kappa$,
				and every $\gamma<\Theta_{X_\infty}$ has cardinality at most $\kappa$.
				Since $\mathbb B$ is ${<}\Theta_{X_\infty}$-closed, it adds no such
				sequences of ground-model elements.
				Hence $\tau\restriction\gamma\in W_1$ for every
				$\gamma<\Theta_{X_\infty}$.
			\end{proof}
		\end{enumerate}

		By \cref{cor:embed-in-ro}, there is in
		$W_2$   a precipitous, fine, normal, uniform,
		$\kappa$-dense ideal $\mcl I_\kappa$ on $Z$, with the ordinal covering
		property relative to $N$, and a complete embedding
		\[
		\mcl P(Z)/\mcl I_\kappa\longrightarrow\ro(\col{\omega}{X_\infty}).
		\]
		By \cref{prop:bls-stage-choice}\textup{(\ref{itm:bls-completed-stage})},
		$W_2\models{\dc_{\kappa_{\alpha+1}}}$, and
		$\hH_{\beta+2}[G][H_{\beta+2}]$ satisfies Choice
		with largest cardinal $\kappa_{\alpha+1}$.
		By \cref{prop:bls-tail-closure}, applied at $\beta+2$
		with parameter $\kappa_{\alpha+1}$, the tail $\mathbb Q^{\textsf{ac},\beta+2}$ is
		${<}(\kappa_{\alpha+1}^+)^{W_2}$-closed. So $\mathbb Q^{\textsf{ac},\beta+2}$ adds no $\kappa$-sequences of ordinals. In particular, \[
		\mcl P(Z)^{N[G][H]}=\mcl P(Z)^{W_2}.
		\]
		In the final model, $\mcl I_\kappa$
		remains countably complete, fine, normal, uniform, and $\kappa$-dense and retains the ordinal covering property relative to $N$.
		We now show that $\mcl I_\kappa^+$ witnesses that $\omega_1$ is generically $\kappa$-supercompact.
		
		\setlength{\emergencystretch}{4em}
		For any $\mcl I_\kappa$-generic ultrafilter $G_{\mcl I_\kappa}$ over
		${N}[G][H]$, the arguments of
		\cref{lem:ideal-ultrapower-comparison,lem:ideal-ultrapower-precip} show that
		$\ult{{N}[G][H]}{G_{\mcl I_\kappa}}$ is well-founded and that its
		ultrapower map agrees with $i_\mu$ on the ordinals.
		In $N[G][H][G_{\mcl I_\kappa}]$, let $M_{G_{\mcl I_\kappa}}$ denote
		the transitive collapse of this ultrapower.
		Since ${N}[G][H]\models{\zfc}$, the \Los{} property gives an elementary
		embedding
		\[
		i_{G_{\mcl I_\kappa}}:{N}[G][H]\longrightarrow M_{G_{\mcl I_\kappa}}
		\]
		with critical point $\omega_1^{{N}[G][H]}=\omega_1$.
		Fineness, normality, and genericity give
		\[
		s\df[\sigma\mapsto\sigma]_{G_{\mcl I_\kappa}}
		=i_{G_{\mcl I_\kappa}}[X_\infty]\in M_{G_{\mcl I_\kappa}}.
		\]
		Fix in ${N}[G][H]$ a bijection $b:X_\infty\to\kappa$. Then
		\[
		i_{G_{\mcl I_\kappa}}[\kappa]=i_{G_{\mcl I_\kappa}}(b)[s]\in M_{G_{\mcl I_\kappa}}.
		\]
		The tail forcing $\mathbb Q^{\textsf{ac},\beta+2}$ from $W_2$ to $N[G][H]$
		is ${<}(\kappa_{\alpha+1}^+)^{W_2}$-closed, which suffices to preserve
		the $\kappa$-density of $\mcl I_\kappa^+$. As in the proof of
		\cref{cor:embed-in-ro}\textup{(3)}, the seed $s$ and $i_{G_{\mcl I_\kappa}}[\kappa]$
		identify the selected representatives on the deciding antichains, so
		$M_{G_{\mcl I_\kappa}}$ is countably closed in $N[G][H][G_{\mcl I_\kappa}]$.

		Thus $\mcl I_\kappa^+$ witnesses that $\omega_1$ is generically
		$\kappa$-supercompact. Since $\kappa$ was arbitrary, $\omega_1$ is
		generically supercompact in ${N}[G][H]$.
		Thus $\mcl I_\kappa$ remains precipitous, and ${N}[G][H]\models{T_{\mathrm{IS}}}$.
	\end{proof}
	
	Combining \cref{lem:gch-holds,lem:ideals-in-final-extension} proves
	\cref{thm:generic-supercompactness-and-gch}.
	\section{Uniformization in minimal Nairian models}\label{sec:uniformization-in-nairian-models}
	In this section, we will prove that under natural hypotheses from inner model theory, $\mu$-Uniformization (see~\cref{def:mu-uniformization}) holds in minimal Nairian models.
	Thus let us assume that $V=L(\mcl P(\R))$ is a model of $\treg+\hpc$. Suppose that there are Nairian models, and let us fix a minimal one, denoted by $N$. Let $\xi\df\ord\cap N$.

	For each uncountable $X\in N$, let $\widehat{\mu}_X$ be the normal fine ultrafilter on 
	$\mathcal P_{\omega_1}(X)$ in $V$, and let $\mu_X\df \widehat{\mu}_X\cap N$. 
	By Woodin's theorem, $\widehat{\mu}_X$ is the quasi-club filter on
	$\mathcal P_{\omega_1}(X)$ in the parent model $V$.
	Goldberg and Sargsyan \cite[Theorem~4.2]{GoldbergSargsyan2026Measures} showed 
	that in $N$, $\mu_X$ is the unique
	supercompactness measure on $\powerset_{\omega_1}(X)^{N}$.
	
	We work in $N$. 
	For an uncountable $X$, we associate to each $A\subseteq X^\omega$ the set
	\[
	A^*\df\{\rang{Z}:Z\in A\},
	\]
	and we naturally define  $\mathcal F_X\subseteq \mcl P(X^\omega)$  by letting
	\[
	A\in\mathcal F_X
	\quad\Longleftrightarrow\quad
	A^*\in\mu_X.
	\]
	Notice that $\mathcal F_X\in N$.
	
	Suppose now that $\kappa$ is a strong cardinal of $\hod\mid\xi$ or is a limit of such strong cardinals. 
	The countable elementary substructures and their enumerations form the sets
	\[
	\begin{aligned}
		\textrm{Sub}_\kappa &= \{Y\prec N\mid\kappa:Y\text{ is countable}\}
		\qquad\text{and}\\
		\textrm{EnSub}_\kappa &= \{Y\in(N\mid\kappa)^\omega:\rang{Y}\in \textrm{Sub}_\kappa\}.
	\end{aligned}
	\]
	We will prove two uniformization theorems, the first for $\textrm{EnSub}_\kappa$ and the second for weakly comeager subsets of $(\textrm{EnSub}_\kappa)^\omega$.
	We use $\mu_{N\mid\kappa}$ and $\mu_{\textrm{EnSub}_\kappa}$ for the
	corresponding measures.
	Since $\textrm{Sub}_\kappa$ is a club, we have
	$\textrm{Sub}_\kappa\in\mu_{N\mid\kappa}$.
	
	For $\nu<\xi$, we use the notation
	\[
	\begin{split}
		\operatorname{Its}(\nu)=\{Z:\omega\to\hod\mid\nu:\ &
		\text{there is a hod pair }\hq\text{ such that}\\
		&\nu_{\hq}\text{ is defined and }
		\pi_{\hq,\infty}[\mathcal M^{\hq}|\nu_{\hq}]
		=\rang{Z}\}.
	\end{split}
	\]
	This definition differs slightly from the one in
	\cite[Definition~12.1 and Notation~12.3]{NM2}.
	
	\subsection{Uniformization for \texorpdfstring{$\mathbb R$}{the reals}}\label{sec:real-uniformization}
	
	Since $N\models \ad_\R$,  we have that uniformization for relations on
	$\R\times\R$ holds in $N$.  In this section we conclude that the reflection theorems  \cite[Theorem~13.7 and the proof of
	Corollary~14.10]{NM2} imply uniformization for $\R\times A$ for any $A\in N$.
	
	\begin{theorem}\label{thm:real-uniformization}
		Let $A\in N$. Assume that $R\in N$ is a binary relation with $R\subseteq \R\times A$ and $\dom{R}=\R$.
		Then there is a function $F:\R\to A$ in $N$ such that $F\subseteq R$.
	\end{theorem}
	\begin{proof}
		By the Reflection Principle, there is $\eta<\xi$ such that $R,\R\in N\mid\eta$, and there is
		a transitive set $\overline{W}\in N$ and an elementary
		embedding
		\[
		\pi:\overline{W}\longrightarrow N\mid\eta
		\]
		such that
		\begin{enumerate}[label=(1.\arabic*),ref=1.\arabic*]
			\item\label{clause:real-uniformization:embedding-in-model}
			$\pi\in N$ and $R\in\rang{\pi}$;
			\item\label{clause:real-uniformization:fixed-reals}
			$\R\subseteq \overline{W}$ and $\pi\restriction\R=\mathrm{id}$;
			\item\label{clause:real-uniformization:real-surjection}
			there is a surjection $G:\R\longrightarrow \overline{W}$ with $G\in N$.
		\end{enumerate}
		We use $G$ as in clause~\textup{(\ref{clause:real-uniformization:real-surjection})}
		and write $\bar R\df \pi^{-1}(R)$. Clearly, $\bar R\in \overline{W}$.
		
		It follows by elementarity that for every $x\in\R$ there is $y\in \overline{W}$ such that $(x, y)\in \bar R$. Then let $F': \R\rightarrow \R$ be such that $F'\in N$ and for every $x\in\R$, $(x, G(F'(x)))\in \bar R$. It follows that $F=\pi\circ G\circ F'$ uniformizes $R$ in $N$. 
	\end{proof}
	
	Because there is a (natural) bijection between $\R$ and $\R^\omega$, we also have uniformization for weakly comeager sets. 
	
	\begin{theorem}\label{thm:real-weak-comeager-uniformization}
		Let $p\in\R^{<\omega}$, and suppose that
		$A\subseteq\R^\omega$ is weakly comeager below $p$ with respect to
		$\mu_\R$.
		Suppose $F:A\longrightarrow N$ is a function in $N$ such that for each $x\in A$, $F(x)\not =\emptyset$.
		Then, in $N$, there is a function $G:A\longrightarrow N$ such that
		\[
		(\forall x\in A)~G(x)\in F(x).
		\]
	\end{theorem}
	
	Uniformization for $\textrm{EnSub}_\kappa$ requires shrinking the given
	measure-one set to a measure-one subset.
	
	\subsection{Uniformization for \texorpdfstring{$\textrm{EnSub}_\kappa$}{EnSub kappa}}\label{sec:iteration-uniformization}
	
	The following are our basic uniformization theorems. We prove only the second theorem, since it implies the first.
	
	While our uniformizing function is defined on a smaller domain, the
	retained ranges still form a measure-one set. We retain all the
	enumerations of these ranges that belong to the original set $A$.
	See clause~\textup{(\ref{clause:measure-uniformization:measure-one})} of
	\cref{thm:measure-uniformization} and
	clause~\textup{(\ref{clause:measure-real-uniformization:measure-one})} of
	\cref{thm:measure-real-uniformization}.
	
	\begin{theorem}\label{thm:measure-uniformization}
		Suppose that
		\begin{itemize}
			\item $\kappa$ is a strong cardinal of $\hod\mid\xi$ or a limit of strong cardinals;
			\item $A\subseteq \textrm{EnSub}_\kappa$ and $A\in\mathcal F_{N\mid\kappa}$;
			\item $R\in N$ and, for some $\kappa'$,
			\[
			R\subseteq A\times(N|\kappa').
			\]
		\end{itemize}
		Suppose that for every $Y\in A$ there is
		$Y'\in N|\kappa'$ with $(Y,Y')\in R$.
		Then there are a set $U\subseteq A$ and a function $F:U\to N|\kappa'$ such that
		\begin{enumerate}[label=(\arabic*),ref=\arabic*]
			\item\label{clause:measure-uniformization:measure-one}
			$U^*\in\mu_{N\mid\kappa}$ and $U=\{Y\in A:\rang{Y}\in U^*\}$;
			
			\item\label{clause:measure-uniformization:function-in-model}
			$F\in N$;
			\item\label{clause:measure-uniformization:uniformizes}
			$F\subseteq R$.
		\end{enumerate}
		
	\end{theorem}
	
	The proof of \cref{thm:measure-real-uniformization} is motivated by
	the reflection theorem of \cite[Theorem~13.7]{NM2}.
	For $A\subseteq \textrm{EnSub}_\kappa$ with $A\in\mathcal F_{N\mid\kappa}$, we restrict
	to $U\subseteq A$ with $U^*\in\mu_{N\mid\kappa}$.
	For $Y\in U$, the set $X_Y=\rang{Y}\cap(\hod\mid\kappa)$ extends to
	the iteration set $Z[X_Y]$ constructed in
	\cref{step:measure-real-uniformization:parameter}.
	If $P_Y$ is the transitive collapse of $\rang{Y}$, then $Y$ induces
	an enumeration $\bar Y:\omega\to P_Y$.
	The HOD predicate $\bar H_Y$ of $P_Y$ satisfies
	$\bar H_Y=\mathcal M^{\hq_Y}\mid\kappa_Y$, where
	$\hq_Y=\hq_{Z[X_Y]}$, as in
	\cref{step:measure-real-uniformization:notation}.
	
	For each $b\in\mathbb R$, we make a real code for $P_Y$, $\bar Y$,
	and the auxiliary data generic over a complete iterate $\hr_{Y,b}$
	of $\hq_Y$. We then construct the local Nairian model $N_{Y,b}$
	and embed $P_Y$ into $N_{Y,b}\mid\kappa_{Y,b}$, as in
	\cref{step:measure-real-uniformization:local-model,step:measure-real-uniformization:map-k}.
	The argument of \cite[Theorem~11.3(2)--(4) and the proof of
	Claim~13.8]{NM2} gives an elementary embedding
	$j_{Y,b}:N_{Y,b}\to N\mid\zeta^0$.
	We arrange that $Y,b,B\in\rang{j_{Y,b}}$ and use $j_{Y,b}$ to send a
	chosen witness in $N_{Y,b}$ to an element of the $(Y,b)$-section of $B$,
	as in \cref{step:measure-real-uniformization:function}.
	The given enumeration $Y$, together with our fixed parameters and $b$,
	determines the model $M_{Y,b}$ of
	\cref{step:measure-real-uniformization:coding}, whose constructibility
	order we use to choose the iterate and the generic filter.
	\Cref{fig:measure-real-uniformization:models} shows the models
	used in this construction.
	
	\begin{theorem}\label{thm:measure-real-uniformization}
		Suppose that
		\begin{itemize}[beginpenalty=10000]
			\item $\kappa$ is a strong cardinal of $\hod\mid\xi$ or a limit of strong cardinals;
			\item $A\subseteq \textrm{EnSub}_\kappa$ and $A\in\mathcal F_{N\mid\kappa}$;
			\item $B\in N$ and, for some $\kappa'$,
			\[
			B\subseteq(A\times\mathbb R)\times(N|\kappa').
			\]
		\end{itemize}
		Suppose that for every $(Y,a)\in A\times\mathbb R$ there is
		$Y'\in N|\kappa'$ such that
		$((Y,a),Y')\in B$.
		Then there are a set $U\subseteq A$ and a function
		$F:U\times\mathbb R\longrightarrow N|\kappa'$ such that
		\begin{enumerate}[label=(\arabic*),ref=\arabic*]
			\item\label{clause:measure-real-uniformization:measure-one}
			$U^*\in\mu_{N\mid\kappa}$ and
			\[
			U=\{Y\in A:\rang{Y}\in U^*\};
			\]
			\item\label{clause:measure-real-uniformization:function-in-model}
			$F\in N$;
			\item\label{clause:measure-real-uniformization:uniformizes}
			$(\forall(Y,a)\in U\times\mathbb R)~((Y,a),F(Y,a))\in B$.
		\end{enumerate}
		
	\end{theorem}
	
	\begin{proof}[Proof of \cref{thm:measure-real-uniformization}]
		\resetproofassertions
		\setcounter{uniformizationproofstep}{0}
		If $P\prec N\mid\gamma$, we informally call
		$P\cap(\hod\mid\gamma)$ the ``HOD predicate'' of $P$.
		If $\sigma:P\to N\mid\gamma$ is an elementary embedding, we use
		this terminology for $\sigma^{-1}[\hod\mid\gamma]$.
		Throughout this proof, we consider these models together with their
		HOD predicates when asserting elementarity.
		
		\uniformizationstep{step:measure-real-uniformization:parameter}{Choosing the parameter $Z$.}
		For $A$ and $B$ as in the statement, we may increase $\kappa'$ if
		necessary so that $\kappa'\geq\kappa$ and $\kappa'$
		is a strong cardinal of $\hod\mid\xi$. We let $\zeta^0$ be the least
		strong cardinal of $\hod\mid\xi$ above $\kappa'$. Thus
		$\zeta^0=(\kappa'^+)^N$. Since $B$ is total on $A\times\mathbb R$,
		the parameters defining $B$ also define $A$.
		By \cite[Theorems~10.6 and~10.13, Lemma~10.16]{NM}
		and \cite[Theorem~12.11]{NM2}, we have $A,B\in N\mid\zeta^0$.
		
		We choose a strong cardinal $\zeta^1$ of $\hod\mid\xi$, an
		inaccessible proper cutpoint $\zeta^2$ of $\hod$, and a hod pair
		$\hq$ catching $\kappa$, $\kappa'$, $\zeta^0$, $\zeta^1$, and
		$\zeta^2$, such that
		\[
		\zeta^0<\zeta^1<\zeta^2<\xi
		\quad\text{and}\quad
		\operatorname{RefGen}_3
		(\hq,\zeta^0_{\hq},\zeta^1_{\hq},\zeta^2_{\hq})
		\]
		holds in the sense of \cite[Definition~9.4]{NM2}.
		We also require that the least Woodin cardinal of $\hod$ above
		$\zeta^0$ is below $\zeta^1$.
		We set
		\[
		Z\df\pi_{\hq,\infty}
		[\mathcal M^{\hq}\mid\zeta^2_{\hq}].
		\]
		We choose enumerations $Z^1$ of $Z$ and $Z^0$ of
		$Z\cap(\hod\mid\zeta^0)$. We require that $\hq$ and $Z^0$
		be chosen so that $A$ and $B$ are definable over $N\mid\zeta^0$
		from $Z^0$. The real $a_0$ codes a formula
		$\psi(v,z)$ such that
		\[
		N\mid\zeta^0\models\text{``$B$ is the unique $v$ such that
			$\psi(v,Z^0)$.''}
		\]
		We let $a_1$ code $a_0$ and a function $\ell:\omega\to\omega$
		such that $Z^0(n)=Z^1(\ell(n))$ for every $n<\omega$.
		
		For $X\subseteq\hod\mid\kappa$, let
		\[
		Z[X]=\{f(x):f\in Z\text{ is a function and }
		x\in X\cap\dom{f}\}.
		\]
		This is the operation of \cite[Definition~12.5]{NM2}.
		The set $D$ consists of the countable $X\prec\hod\mid\kappa$ such that
		\begin{proofassertiongroup}[series=measure-real-trace]
			\item\label{clause:measure-real-uniformization:trace-iteration}
			$Z[X]$ is an iteration set at $\zeta^2$;
			\item\label{clause:measure-real-uniformization:trace-intersection}
			$Z[X]\cap(\hod\mid\kappa)=X$.
		\end{proofassertiongroup}
		In the terminology of \cite[Definition~12.5]{NM2},
		clause~\textup{(\ref{clause:measure-real-uniformization:trace-iteration})} says
		that $X\in\operatorname{Its}(Z)$. Each such $X$ also satisfies the following.
		\begin{proofassertionitems}[resume*=measure-real-trace]
			\item\label{clause:measure-real-uniformization:trace-realization}
			There are hod pairs $\hs$ and $\hw$ such that $\hs$ certifies that
			$Z$ is an iteration set at $\zeta^2$, $\hw$ is a $\zeta^2$-bounded
			complete iterate of $\hs$, and
			\[
			\hq_{Z[X]}=\hw\mid\zeta^2_{\hw},\qquad
			\tau_{Z[X]}=\pi_{\hw,\infty}
			\restriction\mathcal M^{\hq_{Z[X]}}.
			\]
		\end{proofassertionitems}
		Here the tuple $(\hq_{Z[X]},\tau_{Z[X]})$ is defined in
		\cite[Notation~12.3]{NM2}.
		In particular, $\hq_Z=\hq\mid\zeta^2_{\hq}$.
		The set $D$ contains a quasi-club in $V$, so
		\[
		D\in\widehat{\mu}_{\hod\mid\kappa}.
		\]
		Since $\operatorname{Its}(Z)\in N$ by
		\cite[Corollary~12.12]{NM2}, we have $D\in N$.
		Thus $D\in\mu_{\hod\mid\kappa}$.
		
		\uniformizationstep{step:measure-real-uniformization:notation}{Setting up notation.}
		For $Y\in \textrm{EnSub}_\kappa$, we use the notation
		\[
		H_Y=\rang{Y},\qquad X_Y=H_Y\cap(\hod\mid\kappa).
		\]
		Let
		\[
		\begin{split}
			C=\{M\in\mathcal P_{\omega_1}(N\mid\kappa):\ &
			M\prec(N\mid\kappa,\hod\mid\kappa)\text{ and}\\
			&M\cap(\hod\mid\kappa)\in D\}
		\end{split}
		\]
		and
		\[
		U=\{Y\in A:H_Y\in C\}.
		\]
		It follows that $C\in\mu_{N\mid\kappa}$. Hence
		$U^*=A^*\cap C\in\mu_{N\mid\kappa}$ and
		$U=\{Y\in A:\rang{Y}\in U^*\}$, which proves
		clause~\textup{(\ref{clause:measure-real-uniformization:measure-one})}.
		
		For $Y\in U$, $P_Y$ denotes the transitive collapse of $H_Y$, and
		\[
		\sigma_Y:P_Y\longrightarrow N\mid\kappa
		\]
		is the inverse of the transitive collapse of $H_Y$, with range $H_Y$.
		The corresponding enumeration is
		\[
		\bar Y=\sigma_Y^{-1}\circ Y:\omega\longrightarrow P_Y.
		\]
		Thus $\bar Y$ enumerates all of $P_Y$.
		We let $\bar H_Y$ be the image of $X_Y$ under the transitive collapse
		of $H_Y$. More precisely, we define
		\[
		\bar H_Y\df\sigma_Y^{-1}[X_Y]
		=\{u\in P_Y:\sigma_Y(u)\in\hod\mid\kappa\}.
		\]
		This is the HOD predicate of $P_Y$.
		
		In the notation of \cite[Notation~12.3]{NM2}, we set
		\[
		\hq_Y=\hq_{Z[X_Y]},\qquad
		\tau_Y=\tau_{Z[X_Y]}:\mathcal M^{\hq_Y}\longrightarrow\hod\mid\zeta^2.
		\]
		Here $\tau_Y$ is the inverse of the transitive collapse of $Z[X_Y]$
		and has range $Z[X_Y]$.
		Thus,
		\begin{proofassertiongroup}
			\item\label{clause:measure-real-uniformization:base-in-trace}
			$Z\subseteq Z[X_Y]$. 
			\item\label{eq:measure-real-uniformization:collapsed-hod}
			We set $\kappa_Y=\tau_Y^{-1}(\kappa)$,
			$\zeta^0_Y=\tau_Y^{-1}(\zeta^0)$, and $\zeta^1_Y=\tau_Y^{-1}(\zeta^1)$.
			Then
			\[\displaystyle
			\bar H_Y=\mathcal M^{\hq_Y}\mid\kappa_Y,
			\qquad
			\sigma_Y\restriction\bar H_Y
			=\tau_Y\restriction\bar H_Y.
			\]
		\end{proofassertiongroup}
		\uniformizationstep{step:measure-real-uniformization:coding}{Coding, the model $M_{Y,b}$, and the hod pair $\hr_{Y,b}$.}
		Fixing a recursive pairing function $\langle\cdot,\cdot\rangle$ on $\omega$,
		we define the enumeration $m_Y:\omega\to\mathcal M^{\hq_Y}$ by
		\[
		m_Y(\langle s,t\rangle)=
		\begin{cases}
			\tau_Y^{-1}(Z^1(s)(Y(t))),&
			\begin{array}{l}
				Y(t)\in X_Y,\ Z^1(s)\text{ is a function, and}\\
				Y(t)\in\dom{Z^1(s)},
			\end{array}\\[2pt]
			0,&\text{otherwise}.
		\end{cases}
		\]
		It follows that $Z^0$ is coded by $t_Y$, where
		\[
		t_Y:\omega\longrightarrow\mathcal M^{\hq_Y}\mid\zeta^0_Y,
		\qquad t_Y(n)=\tau_Y^{-1}(Z^0(n)).
		\]
		We use a fixed recursive coding of countable structures by reals.
		Let $c_Y$ be the code of $(P_Y,\bar H_Y,\bar Y,t_Y)$ uniquely
		determined by our recursive coding of countable structures and the
		tuple $(m_Y,\bar Y)$. Thus $c_Y$ is uniquely determined by
		$(Z^1,a_1,Y)$.
		
		For $b\in\mathbb R$, we set
		\[
		M_{Y,b}=L^{\Sigma_Y}(m_Y,c_Y,a_1,b),
		\qquad \Sigma_Y=\Sigma^{\hq_Y}.
		\]
		Notice that $M_{Y, b}$ is a model of $\zfc$ and has a canonical well-ordering given by the constructibility hierarchy. 
		We let $\delta_Y$ be the least Woodin cardinal of $\mathcal M^{\hq_Y}$
		above $\zeta^0_Y$. Our choice of $\zeta^1$ gives
		$\zeta^0_Y<\delta_Y<\zeta^1_Y$.
		Working in $M_{Y,b}$, by the proof of
		\cite[Theorem~10.8]{NM}, there is a complete iterate
		$\hr$ of $\hq_Y$ such that
		\begin{enumerate}[label=(\roman*),ref=\roman*]
			\item\label{clause:measure-real-uniformization:iterate-based}
			$\mathcal T_{\hq_Y,\hr}$ is based on $(\zeta^0_Y,\delta_Y)$;
			\item\label{clause:measure-real-uniformization:iterate-genericity}
			$(c_Y,a_1,b)$ is generic for the extender algebra of
			$\mathcal M^{\hr}$ at $\pi_{\hq_Y,\hr}(\delta_Y)$.
		\end{enumerate}
		Then let $\hr_{Y,b}$ be the least such iterate of $\hq_Y$ in the constructibility
		order of $M_{Y,b}$ and set
		\[
		i_{Y,b}=\pi_{\hq_Y, \hr_{Y, b}}:\mathcal M^{\hq_Y}\longrightarrow\mathcal M^{\hr_{Y,b}}.
		\]
		We write
		\[
		\zeta^0_{\hr_{Y,b}}=i_{Y,b}(\zeta^0_Y),\qquad
		\zeta^1_{\hr_{Y,b}}=i_{Y,b}(\zeta^1_Y).
		\]
		We thus have
		\begin{proofassertiongroup}
			\item\label{clause:measure-real-uniformization:iteration-in-ambient-model}
			$\mathcal T_{\hq_Y,\hr_{Y,b}}\in M_{Y,b}$;
			\item\label{eq:measure-real-uniformization:fixed-segment}
			$i_{Y,b}\restriction(\mathcal M^{\hq_Y}\mid\zeta^0_Y)=\mathrm{id}$;
			\item\label{clause:measure-real-uniformization:decoded-collapse}
			$P_Y,\bar H_Y,\bar Y,t_Y\in\mathcal M^{\hr_{Y,b}}[c_Y]$.
		\end{proofassertiongroup}
		
		\uniformizationstep{step:measure-real-uniformization:local-model}{Construction of the local Nairian model and the embedding $j_{Y,b}$.}
		Using clause~\textup{(\ref{clause:measure-real-uniformization:iterate-genericity})}
		and the countability of
		$\mathcal M^{\hr_{Y,b}}$ in $M_{Y,b}$, we have, in $M_{Y,b}$, filters
		\[
		g\subseteq
		\col{\omega}{{<}\zeta^1_{\hr_{Y,b}}}
		\]
		generic over $\mathcal M^{\hr_{Y,b}}$, such that
		\[
		c_Y,a_1,b\in\mathcal M^{\hr_{Y,b}}[g].
		\]
		Let $e_{Y,b}$ be the least such filter in the constructibility
		order of $M_{Y,b}$ and set
		\[
		W_{Y,b}=\mathcal M^{\hr_{Y,b}}[e_{Y,b}].
		\]
		
		Working in $W_{Y,b}$, we form the local hod limit
		\[
		\hm_\infty^{Y,b}
		\df\hm(\hr_{Y,b},\zeta^1_{\hr_{Y,b}})
		\]
		of \cite[Notation~9.2(5)]{NM2}.
		Here $\mathcal M_\infty^{Y,b}$ is the direct limit, under the iteration
		embeddings, of the models $\mathcal M^{\hs}$ for all complete iterates
		$\hs$ of $\hr_{Y,b}$ such that
		\begin{enumerate}[label=(\roman*),ref=\roman*]
			\item\label{clause:measure-real-uniformization:local-system-tree}
			$\mathcal T_{\hr_{Y,b},\hs}\in W_{Y,b}$;
			\item\label{clause:measure-real-uniformization:local-system-length}
			$\operatorname{lh}(\mathcal T_{\hr_{Y,b},\hs})
			<\zeta^1_{\hr_{Y,b}}$;
			\item\label{clause:measure-real-uniformization:local-system-based}
			$\mathcal T_{\hr_{Y,b},\hs}$ is based on
			$\hr_{Y,b}\mid\zeta^1_{\hr_{Y,b}}$.
		\end{enumerate}
		The strategy of $\hm_\infty^{Y,b}$ is induced by $\Sigma^{\hr_{Y,b}}$.
		The direct-limit embedding is
		\[
		d_{Y,b}:\mathcal M^{\hr_{Y,b}}\longrightarrow
		\mathcal M_\infty^{Y,b}.
		\]
		We set
		\[
		\gamma_{Y,b}=d_{Y,b}(\zeta^0_{\hr_{Y,b}}),\qquad
		\kappa_{Y,b}=d_{Y,b}(\kappa_Y).
		\]
		We use the following initial segment of the local Nairian model at
		$\zeta^1_{\hr_{Y,b}}$:
		\[
		N_{Y,b}\df
		\left(N^\ell(\hr_{Y,b},\zeta^1_{\hr_{Y,b}})
		\mid\gamma_{Y,b}\right)^{W_{Y,b}}.
		\]
		Its HOD predicate is
		$\mathcal M_\infty^{Y,b}\mid\gamma_{Y,b}$.
		
		\medskip
		\noindent\begin{minipage}{\linewidth}
			The dashed arrows below indicate the constructions in
			\cref{step:measure-real-uniformization:coding,step:measure-real-uniformization:local-model}.
			The genericity iteration is carried out in $M_{Y,b}$, and
			$e_{Y,b}$ is chosen there. Both local models are formed in $W_{Y,b}$;
			the local Nairian model is then restricted to $\gamma_{Y,b}$ to obtain $N_{Y,b}$.
			\[
			\begin{tikzcd}[column sep=7em,row sep=huge]
				\mathcal M^{\hq_Y}
				\arrow[r,dashed,"\text{genericity iteration}"]
				& \mathcal M^{\hr_{Y,b}}
				\arrow[r,dashed,"\text{adjoin }e_{Y,b}"]
				& W_{Y,b}
				\arrow[dl,dashed,"\text{local hod limit}"']
				\arrow[d,dashed,"\text{local Nairian model}"] \\
				& \mathcal M_\infty^{Y,b} & N_{Y,b}
			\end{tikzcd}
			\]
			\captionof{figure}{The models constructed in
				\cref{step:measure-real-uniformization:coding,step:measure-real-uniformization:local-model}.}
			\label{fig:measure-real-uniformization:models}
		\end{minipage}
		
		The local hod limit $\hm_\infty^{Y,b}$ is a complete iterate of $\hq_Y$.
		By \cite[Theorem~12.6 and Definition~12.8]{NM2}, there is
		a unique map
		\[
		\tau_{Y,b}:\mathcal M_\infty^{Y,b}
		\longrightarrow\hod\mid\zeta^2
		\quad\text{such that}\quad
		(\hm_\infty^{Y,b},\tau_{Y,b})\in\operatorname{rel}(Z[X_Y]).
		\]
		By \cite[Theorem~11.3(2)--(4), Theorem~13.7, and the proof of
		Claim~13.8]{NM2}, the restriction
		\[
		\tau_{Y,b}\restriction
		(\mathcal M_\infty^{Y,b}\mid\gamma_{Y,b}):
		\mathcal M_\infty^{Y,b}\mid\gamma_{Y,b}
		\longrightarrow\hod\mid\zeta^0
		\]
		induces an elementary embedding
		\[
		j_{Y,b}:N_{Y,b}\longrightarrow N\mid\zeta^0
		\]
		in the sense of \cite[Definition~13.1]{NM2}.
		These results give the following properties:
		\begin{proofassertiongroup}
			\item\label{clause:measure-real-uniformization:fixed-reals}
			$j_{Y,b}\restriction\mathbb R^{N_{Y,b}}=\mathrm{id}$;
			\item\label{clause:measure-real-uniformization:hod-predicate}
			The restriction of $j_{Y,b}$ to the HOD predicate of $N_{Y,b}$ is
			$\tau_{Y,b}\restriction
			(\mathcal M_\infty^{Y,b}\mid\gamma_{Y,b})$, and
			\[
			j_{Y,b}^{-1}[\hod\mid\zeta^0]
			=\mathcal M_\infty^{Y,b}\mid\gamma_{Y,b};
			\]
			\item\label{eq:measure-real-uniformization:realization}
			$\displaystyle j_{Y,b}\circ d_{Y,b}\circ i_{Y,b}
			\restriction(\mathcal M^{\hq_Y}\mid\zeta^0_Y)
			=\tau_Y\restriction(\mathcal M^{\hq_Y}\mid\zeta^0_Y)$;
			\item\label{clause:measure-real-uniformization:kappa-image}
			$j_{Y,b}(\kappa_{Y,b})=\kappa$;
			\item\label{clause:measure-real-uniformization:bounded-local-image}
			There is $\eta_{Y,b}<\gamma_{Y,b}$ such that
			\[
			d_{Y,b}[\mathcal M^{\hr_{Y,b}}\mid\zeta^0_{\hr_{Y,b}}]
			\subseteq\mathcal M_\infty^{Y,b}\mid\eta_{Y,b}.
			\]
		\end{proofassertiongroup}
		The following diagram commutes by
		\textup{(\ref{eq:measure-real-uniformization:realization})}. Each arrow is
		the restriction of the indicated map to the displayed domain.
		\[
		\begin{tikzcd}[column sep=large,row sep=large]
			\mathcal M^{\hq_Y}\mid\zeta^0_Y
			\arrow[r,"i_{Y,b}"] \arrow[d,"\tau_Y"']
			& \mathcal M^{\hr_{Y,b}}\mid\zeta^0_Y
			\arrow[d,"d_{Y,b}"] \\
			\hod\mid\zeta^0
			& \mathcal M_\infty^{Y,b}\mid\eta_{Y,b}
			\arrow[l,"j_{Y,b}"]
		\end{tikzcd}
		\]
		
		\uniformizationstep{step:measure-real-uniformization:map-k}{The map $k_{Y,b}:P_Y\to D_{Y,b}$.}
		We extend the map on the HOD predicate of
		$P_Y$ to all of $P_Y$.
		For this extension, we write
		\[
		D_{Y,b}=N_{Y,b}\mid\kappa_{Y,b},\qquad
		h_{Y,b}=d_{Y,b}\restriction\bar H_Y.
		\]
		In particular,
		$h_{Y,b}:\bar H_Y\to\mathcal M_\infty^{Y,b}\mid\kappa_{Y,b}$.
		We have the following conditions:
		\begin{proofassertiongroup}
			\item\label{eq:measure-real-uniformization:hod-composition}
			$j_{Y,b}(D_{Y,b})=N\mid\kappa$;
			\item\label{clause:measure-real-uniformization:hod-map-composition}
			$j_{Y,b}\circ h_{Y,b}=\sigma_Y\restriction\bar H_Y$;
			\item\label{clause:measure-real-uniformization:code-in-local-model}
			$c_Y\in N_{Y,b}$;
			\item\label{clause:measure-real-uniformization:collapse-data-in-local-model}
			$P_Y,\bar Y,t_Y\in N_{Y,b}$;
			\item\label{clause:measure-real-uniformization:hod-map-in-local-model}
			$h_{Y,b}\in N_{Y,b}$.
		\end{proofassertiongroup}
		
		In the following definition and claim, the HOD predicates of
		$P_Y$ and $D_{Y,b}$ are $\bar H_Y$ and
		$\mathcal M_\infty^{Y,b}\mid\kappa_{Y,b}$, respectively.
		For $u\in P_Y$, we choose a sequence of ordinals $t\in P_Y$
		with domain $\omega$ and a formula $\varphi$ defining $u$ uniquely
		over $P_Y$ from $t$. We define
		\[
		v\df\langle h_{Y,b}(t(n)):n<\omega\rangle.
		\]
		Then $v\in D_{Y,b}$. \textup{(\ref{clause:measure-real-uniformization:hod-map-composition})} implies
		\[
		j_{Y,b}(v)=\sigma_Y(t).
		\]
		We define $k_{Y,b}:P_Y\longrightarrow D_{Y,b}$ by letting, for $u$ and $v$ as above,
		$k_{Y,b}(u)$ be the unique $u'\in D_{Y,b}$ such that
		\[
		D_{Y,b}\models\varphi(u',v).
		\]
		
		\begin{claim}\label{claim:measure-real-uniformization:elementarity}
			The map $k_{Y,b}:P_Y\longrightarrow D_{Y,b}$ is well-defined and elementary.
		\end{claim}
		
		\begin{proof}
			By \textup{(\ref{eq:measure-real-uniformization:hod-composition})}, the restriction
			\[
			j_{Y,b}\restriction D_{Y,b}:D_{Y,b}\longrightarrow N\mid\kappa
			\]
			is elementary. Elementarity of $\sigma_Y$ and of $j_{Y,b}\restriction D_{Y,b}$,
			together with $j_{Y,b}(v)=\sigma_Y(t)$, gives a unique
			$u'\in D_{Y,b}$ such that $D_{Y,b}\models\varphi(u',v)$, and
			$j_{Y,b}(u')=\sigma_Y(u)$.
			Since $j_{Y,b}$ is injective, this equality determines $u'$
			independently of $t$ and $\varphi$. Thus $k_{Y,b}$ is well-defined, and
			\begin{proofassertiongroup}
				\item\label{eq:measure-real-uniformization:collapse-composition}
				$k_{Y,b}:P_Y\longrightarrow D_{Y,b}$ and
				$j_{Y,b}\circ k_{Y,b}=\sigma_Y$.
			\end{proofassertiongroup}
			For every formula $\psi$ and every finite tuple $\bar a$ from $P_Y$,
			elementarity of $\sigma_Y$ and $j_{Y,b}\restriction D_{Y,b}$,
			together with \textup{(\ref{eq:measure-real-uniformization:collapse-composition})}, gives
			\[
			\begin{aligned}
				P_Y\models\psi(\bar a)
				&\quad\Longleftrightarrow\quad
				N\mid\kappa\models\psi(\sigma_Y(\bar a))\\
				&\quad\Longleftrightarrow\quad
				N\mid\kappa\models\psi(j_{Y,b}(k_{Y,b}(\bar a)))\\
				&\quad\Longleftrightarrow\quad
				D_{Y,b}\models\psi(k_{Y,b}(\bar a)).
			\end{aligned}
			\]
			Thus $k_{Y,b}$ is elementary.
		\end{proof}
		
		\medskip
		\noindent\begin{minipage}{\linewidth}
			The next diagram shows how $h_{Y,b}$ extends to $k_{Y,b}$.
			The vertical arrows are inclusions.
			By \textup{(\ref{eq:measure-real-uniformization:collapse-composition})} and
			\textup{(\ref{clause:measure-real-uniformization:hod-map-composition})},
			$j_{Y,b}\circ k_{Y,b}=\sigma_Y$ and
			$j_{Y,b}\circ h_{Y,b}=\sigma_Y\restriction\bar H_Y$.
			Since $j_{Y,b}$ is injective, $k_{Y,b}\restriction\bar H_Y=h_{Y,b}$,
			so the diagram commutes.
			
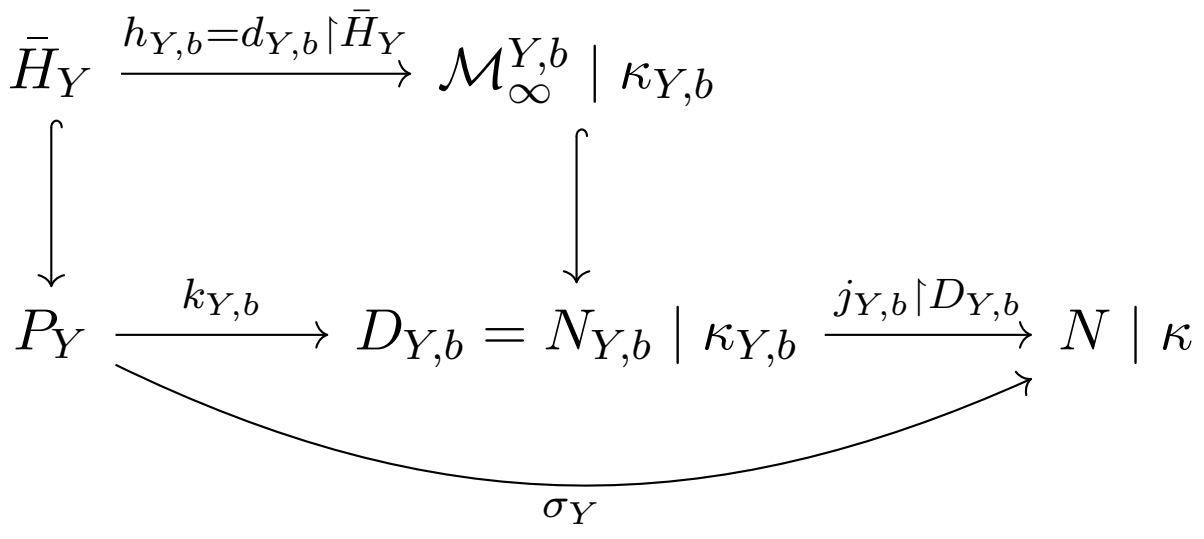
\[
			\begin{tikzcd}[column sep=large,row sep=large]
				\bar H_Y \arrow[r,"h_{Y,b}=d_{Y,b}\restriction\bar H_Y"]
				\arrow[d,hook]
				& \mathcal M_\infty^{Y,b}\mid\kappa_{Y,b} \arrow[d,hook] & \\
				P_Y \arrow[r,"k_{Y,b}"]
				\arrow[rr,bend right=25,"\sigma_Y"']
				& D_{Y,b}=N_{Y,b}\mid\kappa_{Y,b}
				\arrow[r,"j_{Y,b}\restriction D_{Y,b}"]
				& N\mid\kappa
			\end{tikzcd}
			\]
			\bigskip
		\end{minipage}
		
		The structures $(P_Y,\bar H_Y)$ and
		$(D_{Y,b},\mathcal M_\infty^{Y,b}\mid\kappa_{Y,b})$, together with
		$h_{Y,b}$, belong to $N_{Y,b}$. Hence the definition gives
		$k_{Y,b}\in N_{Y,b}$.
		
		We define
		\[
		S_{Y,b}=k_{Y,b}\circ\bar Y,\qquad
		\bar Z^0_{Y,b}(n)=d_{Y,b}(t_Y(n))\quad(n<\omega).
		\]
		Both sequences belong to $N_{Y,b}$. For $\bar Z^0_{Y,b}$, we use
		\textup{(\ref{eq:measure-real-uniformization:fixed-segment})},
		\textup{(\ref{clause:measure-real-uniformization:bounded-local-image})},
		and $t_Y\in W_{Y,b}$. \textup{(\ref{eq:measure-real-uniformization:realization})}
		and \textup{(\ref{eq:measure-real-uniformization:collapse-composition})} give
		\begin{proofassertiongroup}
			\item\label{eq:measure-real-uniformization:enumerations}
			$j_{Y,b}(S_{Y,b})=Y$;
			\item\label{clause:measure-real-uniformization:parameter-enumeration-image}
			$j_{Y,b}(\bar Z^0_{Y,b})=Z^0$.
		\end{proofassertiongroup}
		\uniformizationstep{step:measure-real-uniformization:function}{Construction of the uniformizing function.}
		We have proved the four properties needed for uniformization:
		\begin{proofassertiongroup}
			\item\label{clause:thm:measure-real-uniformization:upsilon-rge-j-upsilon}
			$Y\in\rang{j_{Y,b}}$;
			\item\label{clause:thm:measure-real-uniformization:z-rge-j-upsilon}
			$Z^0\in\rang{j_{Y,b}}$;
			\item\label{clause:thm:measure-real-uniformization:n-upsilon}
			$a_0\in N_{Y,b}$ and $j_{Y,b}(a_0)=a_0$;
			\item\label{clause:thm:measure-real-uniformization:b-n-upsilon}
			$b\in N_{Y,b}$ and $j_{Y,b}(b)=b$.
		\end{proofassertiongroup}
		
		By elementarity of $j_{Y,b}$ and
		\textup{(\ref{clause:measure-real-uniformization:parameter-enumeration-image})},
		there is $\bar B_{Y,b}\in N_{Y,b}$ with $j_{Y,b}(\bar B_{Y,b})=B$.
		Since $B$ is total on $A\times\mathbb R$, elementarity of $j_{Y,b}$
		gives $T\in N_{Y,b}$ such that $((S_{Y,b},b),T)\in\bar B_{Y,b}$.
		Since $W_{Y,b}=\mathcal M^{\hr_{Y,b}}[e_{Y,b}]$, it inherits a
		well-order from the constructibility order of $\mathcal M^{\hr_{Y,b}}$
		by ordering its elements by their least names evaluated at $e_{Y,b}$.
		We let $T_{Y,b}$ be the least such $T$ in this order and set
		\[
		\begin{split}
			F:U\times\mathbb R&\longrightarrow N\mid\kappa',\\
			F(Y,b)&\df j_{Y,b}(T_{Y,b}).
		\end{split}
		\]
		Then $((Y,b),F(Y,b))\in B$ for every $(Y,b)\in U\times\mathbb R$.
		The closure of $N$ under ordinal definability used in
		\cite[Proposition~12.2]{NM2} gives $F\in N$.
		Together with clause~\textup{(\ref{clause:measure-real-uniformization:measure-one})},
		this proves the theorem.
	\end{proof}
	
	\subsection{Uniformization for weakly comeager sets}\label{sec:weak-comeager}
	
	Our goal in this subsection is to prove \cref{thm:weak-comeager-uniformization}, a uniformization theorem for weakly comeager sets.
	
	We let $\kappa$ be a strong cardinal of $\hod\mid\xi$ or a limit
	of such cardinals, and let $\mu_{\textrm{EnSub}_\kappa}$ be the supercompactness
	measure in $N$ on $\mcl P_{\omega_1}(\textrm{EnSub}_\kappa)$.
	We use \cref{def:weak-comeagerness} with $X=\textrm{EnSub}_\kappa$.
	For every nonempty countable $\sigma\subseteq \textrm{EnSub}_\kappa$, the set
	$\{f\in\sigma^\omega:\rang{f}=\sigma\}$ is comeager in $\sigma^\omega$.

	Let $p_0,p_1,\ldots$ be the increasing enumeration of prime numbers.
	For $f\in(\textrm{EnSub}_\kappa)^\omega$, $\operatorname{Code}(f)$ is the sequence
	in $(N\mid\kappa)^\omega$ defined by
	\[
	\operatorname{Code}(f)(i)=
	\begin{cases}
		0, & \text{if }i\ne p_k^{j+1}\text{ for every }k,j<\omega,\\
		f(k)(j), & \text{if }i=p_k^{j+1}.
	\end{cases}
	\]
	This defines an injection
	\[
	\operatorname{Code}\colon(\textrm{EnSub}_\kappa)^\omega\longrightarrow(N\mid\kappa)^\omega.
	\]
	By \cref{lem:directed-supports} applied to $\mu_{\textrm{EnSub}_\kappa}$, if $A$ is
	$\mu_{\textrm{EnSub}_\kappa}$-weakly comeager below $p$, then there is $A'\subseteq A$
	such that $A'$ is $\mu_{\textrm{EnSub}_\kappa}$-weakly comeager below $p$ and
	\[
	\operatorname{Code}[A']\subseteq \textrm{EnSub}_\kappa.
	\]
	
	We need the following lemma.
	
	\begin{lemma}\label{lem:measure-projection}
		Suppose that $\mu_{\textrm{EnSub}_\kappa},\mu_{N\mid\kappa}\in N$ are the
		supercompactness measures on
		$\mcl P_{\omega_1}(\textrm{EnSub}_\kappa)$ and
		$\mcl P_{\omega_1}(N\mid\kappa)$, respectively.
		Suppose $\mathcal A\in\mu_{\textrm{EnSub}_\kappa}$ and $\mathcal B\in\mu_{N\mid\kappa}$.
		Then, for $\mu_{\textrm{EnSub}_\kappa}$-almost all $\sigma\in\mathcal A$,
		$\bigcup\{\rang{u}:u\in\sigma\}\in\mathcal B$.
	\end{lemma}
	\begin{proof}
		Set
		\[
		\mathcal D=\left\{\sigma\in\mathcal A:
		\bigcup\{\rang{u}:u\in\sigma\}\in\mathcal B\right\}.
		\]
		The measures $\mu_{\textrm{EnSub}_\kappa}$ and $\mu_{N\mid\kappa}$ are restrictions to $N$ of
		supercompactness measures in $V$. Since the latter measures are the
		quasi-club filters, we can find $C\subseteq\R$, $\alpha<\Theta^V$, and
		$F\in L_\alpha(C,\R)$ such that
		$\textrm{EnSub}_\kappa,N\mid\kappa\in L_\alpha(C,\R)$ and, setting
		\[
		W=(L_\alpha(C,\R),F,\in),
		\]
		every $Y\prec W$ that is countable in $V$ and contains $\textrm{EnSub}_\kappa$
		satisfies
		\begin{itemize}
			\item $Y\cap \textrm{EnSub}_\kappa\in\mathcal A$,
			\item $Y\cap(N\mid\kappa)\in\mathcal B$.
		\end{itemize}
		For such a $Y$, we write $\sigma=Y\cap \textrm{EnSub}_\kappa$.
		Then $\sigma\in\mathcal A$. If $u\in\sigma$ and $n<\omega$, then
		$u(n)\in Y$ by elementarity. Conversely, every
		$a\in N\mid\kappa$ belongs to the range of some $u\in \textrm{EnSub}_\kappa$,
		obtained by enumerating a countable elementary submodel of
		$N\mid\kappa$ containing $a$. Thus, if
		$a\in Y\cap(N\mid\kappa)$, elementarity gives such a
		$u\in Y\cap \textrm{EnSub}_\kappa$. Hence
		\[
		\bigcup\{\rang{u}:u\in\sigma\}=Y\cap(N\mid\kappa)\in\mathcal B.
		\]
		The set of such $Y$ is a quasi-club in $V$, so $\mathcal D$ contains
		the quasi-club consisting of the intersections $Y\cap \textrm{EnSub}_\kappa$.
		Hence $\mathcal D$ belongs to the supercompactness measure in $V$
		extending $\mu_{\textrm{EnSub}_\kappa}$. Since $\mathcal D\in N$, we have $\mathcal D\in\mu_{\textrm{EnSub}_\kappa}$.
	\end{proof}
	
	\Cref{lem:measure-projection} can also be proved using the
	uniqueness of supercompactness measures.
	
	The proof of \cref{thm:weak-comeager-uniformization} is similar
	to the proof of \cref{thm:real-weak-comeager-uniformization}.
	For $A$ and $F$ as in the statement, the injectivity of
	$\operatorname{Code}$ makes
	$F\circ(\operatorname{Code}\restriction A)^{-1}$ a function on
	$\operatorname{Code}[A]$. We use
	\cref{thm:measure-uniformization} to uniformize this function
	on a suitable subset of $\operatorname{Code}[A]$.
	
	\begin{theorem}\label{thm:weak-comeager-uniformization}
		Let $\mu_{\textrm{EnSub}_\kappa}\in N$ be the supercompactness measure on
		$\mcl P_{\omega_1}(\textrm{EnSub}_\kappa)$.
		Suppose that $A\in N$ is a subset of $(\textrm{EnSub}_\kappa)^\omega$
		that is $\mu_{\textrm{EnSub}_\kappa}$-weakly comeager below
		$p\in\col{\omega}{\textrm{EnSub}_\kappa}$.
		Let $F\colon A\to N$ be a function in $N$ with $F(f)\ne\emptyset$ for every $f\in A$.
		Then there are a set $B\subseteq A$ and a function $G\colon B\to N$
		such that $B$ is $\mu_{\textrm{EnSub}_\kappa}$-weakly comeager below some $q\leq p$
		and $G(f)\in F(f)$ for every $f\in B$.
	\end{theorem}
	\begin{proof} 
		\resetproofassertions
		By \cref{lem:directed-supports} and the observation
		following the definition of $\operatorname{Code}$, we may replace $A$
		by a $\mu_{\textrm{EnSub}_\kappa}$-weakly comeager subset below $p$ and assume that
		\begin{proofassertiongroup}[series=weak-comeager-domain]
			\item\label{clause:weak-comeager-uniformization:codes-in-ensub}
			$\operatorname{Code}[A]\subseteq \textrm{EnSub}_\kappa$.
		\end{proofassertiongroup}
		By fineness of $\mu_{\textrm{EnSub}_\kappa}$, we may choose $\mathcal A'\in\mu_{\textrm{EnSub}_\kappa}$ consisting
		of nonempty sets $\sigma$ with $\rang{p}\subseteq\sigma$ such that
		\begin{proofassertionitems}[resume*=weak-comeager-domain]
			\item\label{clause:weak-comeager-uniformization:comeager-fibers}
			for every $\sigma\in\mathcal A'$, the set
			$\{f\in A:f\colon\omega\to\sigma\}$ is comeager below $p$.
		\end{proofassertionitems}
		Set
		\[
		C\df\{\operatorname{Code}(f):f\in A\},
		\]
		and  define a function $F'$ on $C$ by
		\[
		F'(\operatorname{Code}(f))\df F(f).
		\]
		Let $\mu_{N\mid\kappa}\in N$ be the supercompactness measure on
		$\mcl P_{\omega_1}(N\mid\kappa)$.
		We have the following.
		\begin{proofassertiongroup}[series=weak-comeager-codes,beginpenalty=10000,midpenalty=10000]
			\item\label{clause:weak-comeager-uniformization:coded-function}
			Both $C$ and $F'$ are in $N$, and $F'\colon C\to N$.
			\item\label{clause:weak-comeager-uniformization:nonempty-values}
			$F'(c)\ne\emptyset$ for every $c\in C$.
			\item\label{clause:weak-comeager-uniformization:code-ranges-measure-one}
			$C^*\df \{\rang{c}:c\in C\}\in\mu_{N\mid\kappa}$.
		\end{proofassertiongroup}
		In \textup{(\ref{clause:weak-comeager-uniformization:coded-function})},
		$F'$ is well-defined because $\operatorname{Code}$ is injective,
		and $C,F'\in N$ because $A,F\in N$.
		Also, \textup{(\ref{clause:weak-comeager-uniformization:nonempty-values})}
		holds because $F(f)\ne\emptyset$ for every $f\in A$.
		
		To prove \textup{(\ref{clause:weak-comeager-uniformization:code-ranges-measure-one})},
		suppose otherwise. \Cref{lem:measure-projection} gives a nonempty
		$\sigma\in\mathcal A'$ such that
		$\bigcup\{\rang{u}:u\in\sigma\}\notin C^*$.
		The surjections from $\omega$ onto $\sigma$ form a comeager subset
		of $\sigma^\omega$, so there is $f\in A$ with $\rang{f}=\sigma$.
		Since
		\[
		\rang{\operatorname{Code}(f)}
		=\bigcup\{\rang{u}:u\in\sigma\}\in C^*,
		\]
		we get a contradiction.
		
		By \textup{(\ref{clause:weak-comeager-uniformization:code-ranges-measure-one})},
		we have $C\in\mathcal F_{N\mid\kappa}$.
		We choose $\kappa'<\xi$ such that
		$\bigcup\rang{F'}\subseteq N\mid\kappa'$ and set
		\[
		B'\df\{(c,z):c\in C,\ z\in F'(c)\}.
		\]
		\Cref{thm:measure-uniformization}, applied with $A=C$ and
		$R=B'$, gives $D\in\mu_{N\mid\kappa}$ and a function $H\in N$ such that, setting
		\[
		C'=\{\operatorname{Code}(f)\in C:
		\rang{\operatorname{Code}(f)}\in D\},
		\]
		both of the following hold:
		\begin{proofassertiongroup}
			\item\label{clause:weak-comeager-uniformization:restricted-domain-comeager}
			The set $\{f\in A:\operatorname{Code}(f)\in C'\}$ is
			$\mu_{\textrm{EnSub}_\kappa}$-weakly comeager below $p$.
			\item\label{clause:weak-comeager-uniformization:uniformizing-function}
			$H\colon C'\to N\mid\kappa'$ and
			$(\forall c\in C')~H(c)\in F'(c)$.
		\end{proofassertiongroup}
		To see \textup{(\ref{clause:weak-comeager-uniformization:restricted-domain-comeager})},
		first set \[
		B\df\{f\in A:\operatorname{Code}(f)\in C'\},
		\]
		and 
		\[
		U\df\{\sigma\in\mathcal A':(\exists f\in B)(\rang{f}=\sigma)\}.
		\]
		Notice now that for each $\sigma\in U$, we have that \[
		\{f\in B:\rang{f}=\sigma\}
		=\{f\in A:\rang{f}=\sigma\}.
		\]
		For $\mu_{\textrm{EnSub}_\kappa}$-almost every $\sigma$, the set
		$\{f\in A:\rang{f}=\sigma\}$ is comeager below $p$.
		Thus \textup{(\ref{clause:weak-comeager-uniformization:restricted-domain-comeager})}
		follows once we show that
		\begin{proofassertiongroup}
			\item\label{clause:thm:weak-comeager-uniformization:suffices-show-u-mu-have-mu-sigma} $U\in\mu_{\textrm{EnSub}_\kappa}$.
		\end{proofassertiongroup}
		By \cref{lem:measure-projection},
		\[
		\left\{\sigma\in\mathcal A':
		\bigcup\{\rang{u}:u\in\sigma\}\in D\right\}\in\mu_{\textrm{EnSub}_\kappa}.
		\]
		For each $\sigma\in\mathcal A'$, the surjections from $\omega$ onto
		$\sigma$ form a comeager subset of $\sigma^\omega$.
		Thus \textup{(\ref{clause:weak-comeager-uniformization:comeager-fibers})}
		gives $f\in A$ with $\rang{f}=\sigma$.
		For every $f\in A$ with $\rang{f}=\sigma$, we have
		\[
		\rang{\operatorname{Code}(f)}=\bigcup\{\rang{u}:u\in\sigma\}.
		\]
		Consequently,
		\[
		\begin{aligned}
			\sigma\in U
			&\Longleftrightarrow \sigma\in\mathcal A'\ \land\
			(\exists f\in B)~\rang{f}=\sigma\\
			&\Longleftrightarrow \sigma\in\mathcal A'\ \land\
			(\exists f\in A)~[\rang{f}=\sigma\ \land\
			\operatorname{Code}(f)\in C']\\
			&\Longleftrightarrow \sigma\in\mathcal A'\ \land\
			(\exists f\in A)~[\rang{f}=\sigma\ \land\
			\rang{\operatorname{Code}(f)}\in D]\\
			&\Longleftrightarrow \sigma\in\mathcal A'\ \land\
			\bigcup\{\rang{u}:u\in\sigma\}\in D.
		\end{aligned}
		\]
		Hence $U\in\mu_{\textrm{EnSub}_\kappa}$.
		
		For each $f\in B$, define
		\[
		G(f)\df H(\operatorname{Code}(f)).
		\]
		Clearly, $G$ uniformizes $F\restriction B$.
	\end{proof}

	We obtain the following corollary from
	\cref{thm:weak-comeager-uniformization}.

	\begin{corollary}\label[corollary]{thm:uniformization_in_minimal_N_model}
		Let $N$ be a minimal Nairian model.
		Let $\alpha\in N\cap\ord$. 
		Let $\mu\in N$ be a supercompactness
		measure on $\mcl P_{\omega_1}(\kappa_\alpha^\omega)$.
		Then
		$\mu$-Uniformization holds in $N$.
	\end{corollary}
	\begin{proof}
		It is enough to prove uniformization for relations
		\[
		R\subseteq(\kappa_\alpha^\omega)^\omega
		\times(N\mid\kappa'),
		\]
		where $\kappa'\geq\kappa_\alpha$ and
		$\{\rang{f}:f\in\dom{R}\}\in\mu$. We fix such an $R$ and also a surjection
		\[
		F:\kappa_\alpha^\omega\twoheadrightarrow N\mid\kappa_\alpha.
		\]
    For each $f\in\dom{R}$, we let $g_f$ be an enumeration of $F[\rang f]$ such that $f$ and $g_f$ can be decoded from each other using $F$ as an oracle. For example, letting $p_i$ be the $i$th odd prime number, we could enumerate members of the set $\{F(f(i)): i\in \omega\}$ on the even coordinates and enumerate $f(i)$ using the powers of $p_i$: more precisely, by letting $g_f(p_i^{k+1})=f(i)(k)$. 
    
Notice next that for $\mu_{\kappa_\alpha^\omega}$-almost every
		$\sigma\in\mcl P_{\omega_1}(\kappa_\alpha^\omega)$, we have
		$F[\sigma]\prec N\mid\kappa_\alpha$. To see this, we work in the parent
		model $V$. As in the proof of \cref{lem:measure-projection}, we choose
		$C\subseteq\R$ and $\beta<\Theta^V$ such that
		$F,\kappa_\alpha^\omega,N\mid\kappa_\alpha\in L_\beta(C,\R)$, and we set
		\[
		W\df(L_\beta(C,\R),F,\in),
		\]
		where $F$ is named by a constant. If $Y\prec W$ is countable in $V$,
		then $\dom{F},\rang{F}\in Y$. By elementarity and the surjectivity of $F$,
		\[
		F[Y\cap\kappa_\alpha^\omega]=Y\cap(N\mid\kappa_\alpha)
		\prec N\mid\kappa_\alpha.
		\]
		The intersections $Y\cap\kappa_\alpha^\omega$ form a quasi-club in $V$.
		Hence $F[\sigma]\prec N\mid\kappa_\alpha$ for
		$\widehat\mu_{\kappa_\alpha^\omega}$-almost every $\sigma$.
		The set of such $\sigma$ belongs to $N$, so it also belongs to
		$\mu_{\kappa_\alpha^\omega}$. The claim now follows from  \cref{thm:measure-real-uniformization}, which allows us to uniformize the relation $R'\subseteq (N|\kappa_\alpha)^\omega\times N|\kappa'$ given by $(g_f, x)\in R'$ if and only if $(f, x) \in R$. The key point here is that this uniformization shrinks the underlying measure one set but not the set of enumerations of its members. 
	\end{proof}

	\section{The main theorems}\label{sec:mainthms}
	We now prove our main theorems.
	
	\begin{theorem}[$V=L(\mcl P(\R))+\hpc+\treg$]\label{Mainthm1}
		Suppose that $N$ is a minimal Nairian model. Let $\alpha\in N\cap\ord$. Suppose that $\mu\in N$ is a supercompactness measure on $\mcl P_{\omega_1}(\kappa_\alpha^\omega)$.  Then $i_{\mu}(\omega_1)=\kappa_{\alpha+1}$, where  $i_{\mu}$  is the $\mu$-ultrapower embedding over $N$. 
	\end{theorem}
	\begin{proof}
		By \cref{thm:theta-sequence-properties,thm:uniformization_in_minimal_N_model}, we can apply \cref{thm:image-of-omega-one} to $N$ and obtain the result.
	\end{proof}
	\begin{theorem}[$V=L(\mcl P(\R))+\hpc+\treg$]\label{Mainthm2}
		Suppose that $N$ is a minimal Nairian model. Then $N$ has a class forcing extension in which $T_{\mathrm{IS}}+\gch$ holds.
	\end{theorem}
	\begin{proof}
		By \cref{thm:uniformization_in_minimal_N_model}, $\mu$-Uniformization holds in $N$
		for every ordinal $\alpha\in N$ and every supercompactness measure $\mu\in N$
		on $\mcl P_{\omega_1}(\kappa_\alpha^\omega)$.
		Thus we can apply
		\cref{thm:generic-supercompactness-and-gch} to find  the required
		class forcing extension.
	\end{proof}

	\begin{theorem}\label{Mainthm3}
		The consistency strength of $T_{\mathrm{IS}}+\gch$ is strictly below that of the existence of a Woodin cardinal which is a limit of Woodin cardinals.
	\end{theorem}
	\begin{proof}
		By an unpublished result of Sargsyan, the theory $T$ consisting of $\zfc$ together with the assertion that there exists a Woodin cardinal that is a limit of Woodin cardinals proves that it is consistent that there exists a transitive model $V$ satisfying $V=L(\mcl P(\R))+\hpc+\treg$ in which there are Nairian models. By \cref{Mainthm2}, $T$ proves the consistency of $T_{\mathrm{IS}}+\gch$. The strictness follows from Gödel's Incompleteness Theorem.
	\end{proof}
	
	\begin{remark}
		Under the same assumptions, the same argument shows that if one forces Choice over a minimal Nairian model, starting with $\mathbb P_{\rm max}$ as in \cite{NM2}, then one obtains the ideals with the same properties on every regular cardinal $\kappa\geq\omega_2$ in the generic extension.
	\end{remark}
	\section*{AI disclaimer}
	\addcontentsline{toc}{section}{AI disclaimer}
	The authors discovered the theorems and their proofs before the summer
	of 2026, without any assistance from AI. Many parts of the paper were
	written before advanced AI models became available. Since September
	2026, the authors have used AI for the following tasks:
	\begin{enumerate}[label=\arabic*.]
		\item To convert the second author's handwritten notes into \LaTeX.
		\item To find errors and typos in the paper.
		\item To improve the exposition.
	\end{enumerate}
	
	\printbibliography[heading=bibintoc]

@misc{GoldbergSargsyan2026Measures,
	title={Measures in Nairian models, I}, 
	author={Gabriel Goldberg and Grigor Sargsyan},
	year={2026},
	eprint={2610.19088},
	archivePrefix={arXiv},
	primaryClass={math.LO},
	url={https://arxiv.org/abs/2610.19088}, 
}

@article{ForemanGLC,
    author = {Foreman, Matthew},
    title = {Generic large cardinals: New axioms for mathematics?},
    journal = {Documenta Mathematica},
    volume = {Extra Volume ICM 1998, II},
    pages = {11--21},
    year = {1998},
}

@article{ASTWZ24,
    author = {Adolf, Dominik and Sargsyan, Grigor and Trang, Nam
              and Wilson, Trevor M. and Zeman, Martin},
    title = {Ideals and strong axioms of determinacy},
    journal = {J. Amer. Math. Soc.},
    volume = {37},
    number = {4},
    pages = {1203--1273},
    year = {2024},
    doi = {10.1090/jams/1041},
    eprint = {2111.06220v2},
    archivePrefix = {arXiv},
    primaryClass = {math.LO},
    url = {https://arxiv.org/abs/2111.06220v2},
    addendum = {Numbering follows the arXiv version},
}

@misc{EH24,
    author = {Eskew, Monroe and Hayut, Yair},
    title = {Dense ideals},
    year = {2024},
    note = {Preprint},
    eprint = {2410.14359},
    archivePrefix = {arXiv},
    url = {https://arxiv.org/abs/2410.14359},
}

@book {Wilson-thesis,
    AUTHOR = {Wilson, Trevor Miles},
     TITLE = {Contributions to {D}escriptive {I}nner {M}odel {T}heory},
      NOTE = {Thesis (Ph.D.)--University of California, Berkeley},
 PUBLISHER = {ProQuest LLC, Ann Arbor, MI},
      YEAR = {2012},
 PAGETOTAL = {102},
      ISBN = {978-1267-97617-8},
   MRCLASS = {99-05},
  MRNUMBER = {3152619},
       URL =
              {http://gateway.proquest.com/openurl?url_ver=Z39.88-2004&rft_val_fmt=info:ofi/fmt:kev:mtx:dissertation&res_dat=xri:pqm&rft_dat=xri:pqdiss:3555999},
}

@incollection{Sol,
	author = {Solovay, Robert M.},
	booktitle = {Large cardinals, determinacy and other topics. {T}he {C}abal {S}eminar. {V}ol. {IV}},
	pages = {66--95},
	publisher = {Cambridge Univ. Press, Cambridge},
	series = {Lect. Notes Log.},
	title = {The independence of {DC} from {AD}},
	volume = {49},
	year = {2021},
	}

@misc{NM,
      title={{Nairian} Models}, 
      author={Douglas Blue and Paul B. Larson and Grigor Sargsyan},
      year={2025},
      eprint={2501.18958},
      archivePrefix={arXiv},
      primaryClass={math.LO},
      url={https://arxiv.org/abs/2501.18958}, 
}

@misc{NM2,
      title={The failure of square at all uncountable cardinals is weaker than a {Woodin} limit of {Woodin} cardinals}, 
      author={Douglas Blue and Paul B. Larson and Grigor Sargsyan},
      year={2026},
      eprint={2602.13077},
      archivePrefix={arXiv},
      primaryClass={math.LO},
      url={https://arxiv.org/abs/2602.13077}, 
      note={Theorem numbers refer to the revised manuscript, also circulated as
            {Forcing the Axiom of Choice over Nairian Models}},
}

@book{Steel23,
    author = {Steel, John R.},
    title = {A Comparison Process for Mouse Pairs},
    series = {Lecture Notes in Logic},
    volume = {51},
    publisher = {Cambridge University Press},
    year = {2023},
    doi = {10.1017/9781108886840},
}

@incollection{KecWoo08inbookGene,
	author = {Kechris, Alexander S. and Woodin, W. Hugh},
	booktitle = {Games, Scales, and {Suslin} Cardinals: The {Cabal Seminar}, Volume {I}},
	editor = {Kechris, Alexander S. and L{\"o}we, Benedikt and Steel, John R.},
	series = {Lecture Notes in Logic},
	volume = {31},
	pages = {379--397},
	publisher = {Cambridge University Press},
	title = {Generic codes for uncountable ordinals, partition properties, and elementary embeddings},
	year = {2008},
	doi = {10.1017/CBO9780511546488.019},
}
\end{document}